\documentclass[12pt,a4paper]{amsart}
\usepackage{amsmath}
\usepackage{amscd}

\numberwithin{equation}{section}

\theoremstyle{plain} 
\newtheorem{thm}{Theorem}[section]
\newtheorem{cor}[thm]{Corollary}
\newtheorem{lem}[thm]{Lemma}
\newtheorem{prop}[thm]{Proposition}

\newtheorem{defn}[thm]{Definition}

\newtheorem{rem}[thm]{Remark}

\theoremstyle{remark}

\newcommand{\thmref}[1]{Theorem~\ref{#1}}
\newcommand{\propref}[1]{Proposition~\ref{#1}}

\newcommand{\lemref}[1]{Lemma~\ref{#1}}

\newcommand{\remref}[1]{Remark~\ref{#1}}
\newcommand{\defnref}[1]{Definition~\ref{#1}}

\newsavebox{\SmallMathBox}

\def\pdo{\psi{\rm do}}

\def\Ci{C^\infty}

\def\dd{\partial}

\def\Di{D\kern -.65em /}
\def\Dii{D\kern -.45em /}
\def\di{{\dd}\kern -.55em /}
\def\dii{{\dd}\kern -.40em /}

\def\noi{\noindent}

\def\ol{\overline}
\def\op{\oplus}

\def\to{\rightarrow}
\def\too{\longrightarrow}
\def\mto{\mapsto}
\def\mtoo{\longmapsto}

\def\ws{\widehat{\sigma}}
\def\wt{\widetilde}

\def\wR{\widehat{R}}

\def\re{{\rm Re}}

\def\ZZ{{\bf Z}}
\def\pp{{\bf p}}

\def\Aa{{\mathcal A}}

\def\Cc{{\mathcal C}}
\def\Dd{{\mathcal D}}
\def\Ee{{\mathcal E}}
\def\Ff{{\mathcal F}}

\def\Hh{{\mathcal H}}

\def\Kk{{\mathcal K}}

\def\Mm{{\mathcal M}}
\def\Nn{{\mathcal N}}

\def\Rr{{\mathcal R}}
\def\Ss{{\mathcal S}}

\def\={\cong}
\def\>{\supset}
\def\<{\subset}
\def\ii{^{-1}}
\def\si{^{-s}}
\def\st{^{\star}}
\def\pp{^{\perp}}
\def\12{\frac{1}{2}}
\def\2{\Dd}
\def\3{\Nn}
\def\4{\Rr}
\def\6{\cup}
\def\8{\otimes}
\def\0{^{\circ}}
\def\){\hfill{\ \qed}\enddemo}

\def\a{\alpha}

\def\b{\beta}
\def\R{\mathbb{R}}
\def\C{\mathbb{C}}
\def\Cn{\mathbb{C}^{n}}
\def\Crn{\mathbb{C}^{rn}}
\def\Cnr{\mathbb{C}^{rn}}

\def\d{\delta}
\def\D{\Delta}
\def\e{\varepsilon}
\def\ep{\epsilon}
\def\f{\varphi}

\def\g{\gamma}
\def\G{\Gamma}

\def\th{\theta}

\def\la{\lambda}

\def\o{\infty}
\def\s{\sigma}
\def\Si{\Sigma}

\def\z{\zeta}
\def\Z{\ZZ}

\def\Det{\mbox{\rm Det\,}}
\def\DET{\mbox{\rm DET\,}}

\def\dom{\mbox{\rm dom\,}}

\def\End{\mbox{\rm End}}

\def\GL{\mbox{\rm GL}}

\def\Gr2n{\mbox{${\rm Gr}(\Cn\oplus\Cn)$}}

\def\Grk2n{\mbox{${\rm Gr}_{k}(\Cn\oplus\Cn)$}}
\def\Grk{\mbox{${\rm Gr}_{k}$}}
\def\Gr{\mbox{${\rm Gr}$}}

\def\Hom{\mbox{\rm Hom}}

\def\Ind{\mbox{\rm Ind\,}}

\def\gr{\mbox{\rm graph}}

\def\LIM{\mbox{\rm LIM}}

\def\Ker{\mbox{\rm Ker}}
\def\ker{\mbox{\rm ker}}
\def\Cok{\mbox{\rm Coker}}

\def\ind{\mbox{\rm ind\,}}

\def\ran{\mbox{\rm ran\,}}

\def\Si{S\kern -.65em /}

\def\tr{\mbox{\rm tr\,}}
\def\Tr{\mbox{\rm Tr\,}}

\def\Pf{\mathbb{P}}

\def\wP{{\widehat P}}
\def\Sf{\mathbb{S}}

\def\wS{{\widehat S}}

\def\wg{{\widehat g}}
\def\wga{{\widehat \gamma}}

\def\wD{{\widehat \D}}
\def\wDla{{\widehat \D}_{\la}}

\def\wK{{\widehat K}}
\def\wQ{{\widehat Q}}

\def\wKk{{\widehat \Kk}}

\def\wM{{\widehat M}}
\def\wN{{\widehat N}}
\def\wMmla{{\widehat \Mm}_{\la}}

\def\Mmla{\Mm_{\la}}

\def\wk{{\widehat k}}
\def\Df{\mathbb{D}}

\def\wDd{{\widehat D}}
\def\wE{{\widehat \Ee}}

\def\wpsi{{\widehat \psi}}
\def\dla{\frac{\dd}{\dd\la}}

\def\dD{{\dot {D}}}
\def\dDD{{\dot {\Delta}}}

\def\wMm{{\widehat \Mm}}

\def\lrb{\left( }
\def\rrb{\right) }
\def\lsb{\left[ }
\def\rsb{\right] }

\def\t{\tilde}

\def\ord{\mbox{\rm ord}}

\def\bb{\bullet}

\begin{document}

\title[Determinants of global boundary problems]{{\large Zeta Determinants on Manifolds with Boundary}
\vskip 5mm {\Small {\rm Simon Scott}}}

\vskip 2cm

\begin{abstract}
We study the $\z$-determinant of global boundary problems of
APS-type through a general theory for relative spectral
invariants. In particular, we compute the $\z$-determinant for
Dirac-Laplacian boundary problems in terms of a scattering
Fredholm determinant over the boundary.
\end{abstract}

\maketitle

\tableofcontents


\section{Introduction}

The purpose of this paper is to study the $\z$-function
regularized determinant for the Dirac Laplacian of an
Atiyah-Patodi-Singer (APS)-type boundary problem. We do so within
a general framework for studying relative global spectral
invariants on manifolds with boundary. Despite the primacy of the
determinant of the Dirac Laplacian over closed manifolds,
relatively little has been known in the case of global boundary
problems of APS-type. Geometric index theory of such boundary
value problems began with the index formula \cite{AtPaSi75}
\begin{equation}\label{e:apsindex}
\ind(D_{\Pi_{\geq}}) = \int_X \omega(D) - \frac{\eta(\Dd_Y) +
\dim\Ker(\Dd_Y)}{2} \ .
\end{equation}
Here $D : \Ci(X,E^1)\too \Ci(X, E^2)$  is a first-order elliptic
differential operator of Dirac-type acting over a compact manifold
$X$ with boundary $\dd X= Y$. Near the boundary $D$ is assumed to
act in the tangential direction via a first-order self-adjoint
elliptic operator $\Dd_Y$ over  $Y$.  A boundary problem $D_{B} =
D$ is defined by restricting the domain of $D$ to those sections
whose boundary values lie in the kernel of a suitable order zero
pseudodifferential operator $B$ on the space of boundary fields.
APS-type boundary problems refer to the case where $B$ is, in a
suitable sense, `comparable' to the projection $\Pi_{\geq}$ onto
the eigenspaces of $\Dd_Y$ with non-negative eigenvalue. The
other ingredients in \eqref{e:apsindex} are the index density
$\omega(D)$ restricted from the closed double, and the
eta-invariant $\eta(\Dd_Y)$, defined as the meromorphically
continued value at $s=0$ of $\eta(\Dd_Y,s) = \Tr(\Dd_Y
|\Dd_Y|^{-s-1})$.

A striking consequence of \eqref{e:apsindex} is that, in contrast
to the case of closed manifolds, the index of APS-type boundary
problems is not a homotopy invariant. If, however, we restrict
the class of boundary conditions to a suitable classifying space
for even K-theory, then the index does become a homotopy
invariant of the boundary condition. An appropriate parameter
space is a restricted Grassmannian $Gr(D)$ of pseudodifferential
operator ($\pdo$) projections $P$ on the $L^2$ boundary fields
which are comparable to $\Pi_{\geq}$. Such a Grassmannian has
homotopy type $\mathbb{Z}\times BU$ and its connected components
$Gr^{(r)}(D)$ are labeled by the index $r = \ind(D_P)$. One moves
between the different components according to the {\em
relative-index formula}
\begin{equation}\label{e:index1}
 \ind (D_{P_1}) -  \ind (D_{P_2})= \ind (P_2,P_1) \ ,
\end{equation}
where $ (P_2,P_1) := P_1\circ P_2  :\ran(P_2)\to \ran(P_1)$ acts
between the ranges of the projections $P_1, P_2\in Gr(D)$.

The identity \eqref{e:index1} depends on two decisive properties
of global boundary problems over $Gr(D)$. The first is analytic:
the restriction to the boundary of the infinite-dimensional
solution space $\Ker(D)$, to the subspace $H(D):=\Ker(D)_{|Y}$ of
boundary sections, is a continuous bijection, with canonical left
inverse defined by the Poisson operator of $D$. The resulting
isomorphism between the finite-dimensional kernel of $D_P$ and
the kernel of the boundary operator $S(P) = P\circ P(D) : H(D)\to
\ran(P)$, where $P(D)$ is the Calderon projection, and similarly,
between the kernels of the adjoint operators, means that
\begin{equation}\label{e:index2}
  \ind (D_P) = \ind(S(P)) \ .
\end{equation}
The second property is geometric: $Gr(D)$ is a homogeneous
manifold, acted on transitively by an infinite-dimensional
restricted general linear group on the space of boundary fields,
resulting in the identity $\ind (P_1,P_2) + \ind (P_2,P_3) = \ind
(P_1,P_3)$ for any $P_1, P_2, P_3 \in Gr(D)$. Then
\eqref{e:index1} follows trivially from \eqref{e:index2}.

Although the relative index formula is quite classical, these two
properties resonate more powerfully when one turns to the harder
problem of computing for APS-type boundary problems the spectral
and differential geometric invariants familiar from closed
manifolds. A precise understanding of these invariants is crucial
for a direct approach to a geometric index theory of global
boundary problems of APS-type parallel to that for closed
manifolds \cite{BeGeVe92,BiFr85}. An important but rather
different perspective is provided by the $b$-calculus developed by
Melrose \cite{MePi97}.

One spectral invariant that certainly is well understood is the
$\eta$-invariant $\eta(D_P)$ for self-adjoint global boundary
problems over odd-dimensional manifolds. As a result of its
semi-local character the $\eta$-invariant obeys a strikingly
simple (to state) additivity property with respect to a partition
of a closed manifold \cite{BrLe98,KiLe00,Mu94,Wo99}. In this case
the homogeneous structure of the `self-adjoint' Grassmannian
takes a much simpler form: the range of any such $P$ occurs as
the graph of a unitary isomorphism $T:F^+ \to F^-$, where
$F^{\pm}$ denote the spaces of boundary chiral spinor fields
\cite{Sc95}.

The next invariant in the spectral hierarchy remains far more
mysterious. The spectral $\z$-function of the Laplacian boundary
problem $\D_P = (D_P)^* D_P$ is defined for $\re(s)>>0$ by the
operator trace
\begin{equation*}
  \z(\D_P,s) = \Tr(\D_P\si) \ ,
\end{equation*}
where we assume that $D_P$ is invertible. Recent results of Grubb,
following earlier joint work with Seeley \cite{Gr99,Gr99a,GrSe96},
show that for $P$ in the `smooth' Grassmannian $Gr_{\o}(D)$ (see
$\S 3$), $\z(\D_P,s)$ has a meromorphic continuation to $\C$ which
is regular at $s=0$. This means there is a well-defined
regularized $\z$-determinant of the Laplacian
\begin{equation}\label{e:zdet}
  {\rm det}_{\z}(\D_P) = \exp\lrb
  -\frac{d}{ds}_{|s=0}\z(\D_P,s)\rrb \ .
\end{equation}
Owing to its highly non-local nature ${\rm det}_{\z}(\D_P)$ is a
hopelessly difficult invariant to compute. There is, on the other
hand, a quite different but also completely canonical
regularization of the determinant of $\D_P$ as the Fredholm
determinant of the boundary `Laplacian' $$S(P)^* S(P): H(D) \too
H(D) \ .$$ The Fredholm determinant $\det_F$ is defined for
operators on a Hilbert space differing from the identity by an
operator of trace class and is the natural extension to
infinite-dimensional spaces of the usual determinant in finite
dimensions. Its analytical status, however, is essentially
opposite to that of \eqref{e:zdet}. More precisely, the
$\z$-determinant is {\em not\,} an extension of the Fredholm
determinant---operators with Fredholm determinants do not have
$\z$-determinants, operators with $\z$-determinants do not have
Fredholm determinants\footnote{We consider here
infinite-dimensional Hilbert spaces.}. A more subtle fact is
nevertheless true: the {\it relative} $\z$-determinant is a true
extension of the Fredholm determinant. Here a relative
regularized determinant means a regularization of the ratio $\det
A_1 / \det A_2$ for `comparable' operators $A_1, A_2$ (see $\S
2$). Thus, any pair of determinant class (= Id + Trace Class)
operators have a well-defined relative $\z$-determinant.
Moreover, there is a, roughly converse, `relativity principle for
determinants' which states that ratios of $\z$-determinants for
certain preferred classes of unbounded operators can be written
canonically in terms of Fredholm determinants.

Applied to global boundary problems the relativity principle for
determinants is a restatement of the fact that in order to define
a topologically meaningful Grassmannian one must do so relative
to a basepoint projection. (Dimension one is an exception since
we are reduced in this case to the usual finite-dimensional
Grassmannian, no basepoint is needed and for this reason explicit
formulas for the $\z$-determinant of ordinary differential
operators are possible. See $\S 5$.) This is familiar in physics
in the quantization of Fermions where the basepoint corresponds
to the Dirac sea splitting into positive and negative energy
modes (the APS splitting). The application to the determinant is
a well known, if imprecise, idea in physics folklore used
extensively in defining path integrals in QFT and String Theory.

In $\S2$ we prove a precise formulation of the relativity
principle for determinants adequate for our purposes here. The
main result in this paper is Theorem A in the following table. The
table summarizes relative formulas for the key spectral
invariants.

\vskip 5mm

\begin{center}
\begin{tabular}{|l|c|} \hline
 \hskip 23mm Invariant & Relative Formula \\ \hline
  &   \\
    Index &   $\ind (D_{P_1}) -  \ind (D_{P_2})= \ind (P_2,P_1)$ \\
     &  \\
  Eta-invariant: Odd-dimensions & $\eta(D_{P_1},0) - \eta(D_{P_2},0)
  = \frac{1}{\pi i}\log\det_F (T_2\ii T_1)$  \\
   &  \\
  Zeta-determinant: Odd-dimensions & $\frac{{\rm det}_{\z}
  D_{P_1}}{{\rm det}_{\z} D_{P_2}}
  =  \frac{{\rm
det}_F \lrb\frac{1}{2}(I + T_1\ii K)\rrb}{{\rm
det}_F\lrb\frac{1}{2}(I + T_2\ii K)\rrb} $ \\
   &  \\
   Laplacian Zeta-determinant  & {\bf Theorem A}: \ \
   $\frac{{\rm det}_{\z}(\D_{P_1})}{{\rm det}_{\z}(\D_{P_2})} =
\frac{{\rm det}_{F}\lrb S(P_1)^* S(P_1 )\rrb }{{\rm det}_{F}\lrb
S(P_2)^* S(P_2 )\rrb } $    \\
  &   \\
   \hline
\end{tabular}
\end{center}

The third formula in the table is Theorem (0.1) of \cite{ScWo99}
for self-adjoint global boundary problems $D_{P_1}, D_{P_2}$ over
a compact odd-dimensional manifold. The operators $T_i, K : F^+
\to F^-$ are the boundary unitary isomorphisms discussed earlier,
with $H(D) = {\rm graph}(K)$. Because, in this case, the
$\eta$-invariant is essentially the phase of the determinant, the
second formula, which holds $\mod{2\Z}$, is an easy corollary
when the operators are invertible \cite{KiLe00}. Theorem A holds
for general Dirac-type operators and in all dimensions. Notice,
furthermore, that it is stated invariantly, independently of the
choice of `coordinates' $T_i$---in $\S 3.4$ we explain how the
second and third formulas in the table are derived from the
invariant general formulas proved in $\S 2$ for the relative
$\eta$-invariant and $\z$-determinant.

The relative determinant formulas in the table encode a certain
spectral duality between the rapidly diverging eigenvalues of the
global boundary problems and the eigenvalues of the boundary
Laplacians, which converge rapidly to $1$.  The point being that,
for comparable global boundary problems, taking quotients
produces arithmetically similar behaviour.

This extends to a differential geometric duality between smooth
families of global boundary problems $(\Df,\Pf) = \{D^{b}_{P_b} \
| \ b\in B\}$ parameterized by a manifold $B$ and the
corresponding family of boundary operators $\Sf(\Pf) = \{S(P_b) =
P_b\circ P(D^b) \ | \ b\in B\}$. Each family has an associated
index bundle and determinant line bundle. Geometrically the
regularized determinants $\det_{\z}\D_P$ and $\det_{F}(S(P)^*
S(P))$ define the so-called Quillen metric \cite{BeGeVe92,BiFr85}
and canonical metric \cite{Sc99} on the respective determinant
line bundles. For example, if $D$ is held fixed and $P$ allowed
to vary in the parameter manifold $B = Gr_{\o}(D)$ the canonical
metric is just the `Fubini-Study' metric  of \cite{PrSe86} over
the restricted Grassmannian. From this view point, Theorem A
expresses the relative equality of these metrics with respect to
the obvious determinant line bundle isomorphism
$\DET(\Df,\Pf)\cong \DET(\Sf(\Pf))$. The following table lists
relative geometric index theory formulas for families of APS-type
boundary problems:

\vskip 5mm

\begin{center}
\begin{tabular}{|l|c|} \hline
 \hskip 13mm Invariant & Relative Formula \\ \hline
 &   \\
Index bundle & $\Ind(\Df,\Pf_1) - \Ind(\Df,\Pf_2) =
\Ind(\Pf_2,\Pf_1)$\\
   &  \\
Determinant line bundle & $\DET(\Df,\Pf_1) \otimes
\DET(\Df,\Pf_2)^{*} \cong \DET(\Pf_2,\Pf_1)$  \\
   &  \\
Zeta curvature & $\Omega^{\z}_1 - \Omega^{\z}_2  =
\Omega^{\Cc}_1 - \Omega^{\Cc}_2$  \\
  &   \\
\hline
\end{tabular}
\end{center}

The relative index bundle formula is taking place in $K^0 (B)$.
For a functional analytic proof  see \cite{DaZh97}. The
determinant line bundle isomorphism is explained in \cite{Sc99},
where, essentially, the definition is given of the canonical
curvature form $\Omega^{\Cc}$ on $\DET(\Sf(\Pf))$. For the
construction of the $\z$-connection on $\DET(\Df,\Pf)$ with
curvature form $\Omega^{\z}$ and proofs of all three identities,
see \cite{Sc01}.

Notice that by setting $P_2 = P(D)$, the relative formulas in the
tables may be re-expressed as an interior term and a boundary
correction term.

\vskip 2mm

There is an essentially immediate application of the methods here
to non-compact manifolds. For a closely related detailed study of
the Laplacian we refer to the seminal paper of Muller \cite{Mu98}.
For an account of how determinants of global boundary problems
fit into the framework of TQFT we refer to \cite{MiSc99}. The
results of this paper were announced in \cite{Sc00}.

\vskip 2mm

The paper is organized as follows.

\vskip 2mm

In $\S 2$ we prove a precise form of the relativity principle for
determinants using regularized limits of Fredholm scattering
determinants (\thmref{t:thmB}). In $\S 2.1$ we explain how this
is related to the heat operator regularization of the
determinant---the more usual scenario for studying scattering
determinants. The relative eta invariant for comparable
self-adjoint operators requires a somewhat different treatment.
In $\S 2.2$ we prove a general formula for the relative eta
invariant as the difference of two scattering determinant limits.
$\S 2.3$ is concerned with a general multiplicativity property
for the zeta determinant.

In  $\S 3$ we first review some analytic facts about first-order
global boundary problems which will be needed as we proceed.  We
 explain how the scattering determinant arises canonically in terms of
natural isomorphisms between the various determinant lines. In
\thmref{t:thmD} we prove an explicit formula for the relative
zeta determinant of first-order global boundary problems. As an
example, we use this formula to give a new derivation of Theorem
(0.1) of \cite{ScWo99}.

In $\S 4$ we use \thmref{t:thmB} to prove a formula for the
relative zeta determinant of the Dirac Laplacian in terms of an
equivalent first-order system. Methods from $\S 3$ then reduce
this to the equality in Theorem A.

In $\S 5$ we give a new proof using our methods of the results in
\cite{LeTo99} for differential operators in dimension one. In
this sense, the results of this paper may be regarded as the
extension of \cite{LeTo99} to all dimensions.

\vskip 2mm

I thank Gerd Grubb for helpful conversations.


\section{Regularized Limits and the Relative Zeta Determinant}

Let $A_1, A_2$ be invertible closed operators on a Hilbert space
$H$ with a common spectral cut $R_{\th} = \{re^{i\th} \ | \ r\geq
0\}, \ \th\in[0,2\pi)$. This supposes $\d, \varrho
>0$ such that the resolvents $(A_i - \la)\ii$ are holomorphic in
the sector
\begin{equation}\label{e:sector}
\Lambda_{\th} = \{z\in \C \ | \ |\arg(z)-\th| < \d \;\;{\rm
or}\;\; |z| < \varrho\}
\end{equation}
and such that the operator norms $\|(A_i - \la)\ii\|$ are
$O(|\la|\ii)$ as $\la\to\infty$ in $\Lambda_{\th}$. For
$\re(s)>0$ one then has the complex powers first studied by
Seeley \cite{Se67}
\begin{equation}\label{e:powers}
  A_{i}\si = \frac{i}{2\pi}\int_{C}\la\si\, (A_i - \la)\ii \,
  d\la \ ,
\end{equation}
where
\begin{equation}\label{e:branch}
\la\si = |\la|\si e^{-is\arg(\la)}, \hskip 8mm \th - 2\pi \leq
\arg(\la) < \th \ ,
\end{equation}
is the branch of $\la\si$ defined by the spectral cut $R_{\th}$,
and $C = C_{\th}$ is the negatively oriented contour
$C_{\th,\downarrow}\cup C_{\rho,\th,\th-2\pi} \cup C_{\th -
2\pi,\uparrow}$, with
\begin{equation}\label{e:contour}
C_{\phi,\downarrow} = \{\la = re^{i\phi} \ | \ \o > r \geq \rho\}
\ , \hskip 5mm C_{\rho,\phi,\phi^{'}} = \{\rho e^{i\th^{'}} \ | \
\phi \geq \th^{'} \geq \phi^{'}\} \ ,
\end{equation}
$$C_{\phi,\uparrow } = \{\la = re^{i\phi} \ | \ \rho \leq r < \o
\} \ ,$$  and $\rho < \varrho$. We assume there is a real $\a_0$
such that the operators
\begin{equation*}
\dd_{\la}^{m}(A_i - \la)\ii = m!(A_i - \la)^{-m-1}
\end{equation*}
are trace class for $m > -\a_0$, with asymptotic expansions as
$\la\to \o$ in $\Lambda_{\th}$
\begin{equation}\label{e:exp1}
  \Tr (\dd_{\la}^{m}(A_i - \la)\ii)   \sim
  \sum_{j=0}^{\infty}\sum_{k=0}^{1} a_{j,k}^{(i)}(m).(-\la)^{-\a_j
  - m}\log^k (-\la) \ ,
\end{equation}
where  $0< m + \a_0 < \ldots < m + \a_j \nearrow +\infty$. Here,
$\log = \log_{\th}$ is the branch of the logarithm specified by
\eqref{e:branch}; changing $\th$ may change the coefficients
$a_{j,k}^{(i)}(m)$. Since $\la^{k-s} \dd_{\la}^{k-1}(A_i - \la)\ii
\to 0 $ as $\la\to \o$ along $C$ for $\re(s)>0$, we can integrate
by parts in \eqref{e:powers} to obtain
\begin{equation}\label{e:A-by-parts}
  A_{i}\si = \frac{1}{(s-1)\ldots (s-m)}.\frac{i}{2\pi}\int_{C}\la^{m-s}\,
  \dd_{\la}^{m}(A_i - \la)\ii \, d\la \ .
\end{equation}
From \eqref{e:exp1} and \eqref{e:A-by-parts}, the operators
$A_{i}\si$ are trace class in the half-plane $\re(s) > 1 + m $ and
we can define there the spectral zeta functions of $A_1$, $A_2$
\begin{equation*}
  \z_{\th}(A_{i},s)  =  \Tr A_i\si, \hskip 15mm \re(s) > 1 + m \ .
\end{equation*}
Substituting the asymptotic expansion \eqref{e:exp1} in
\begin{equation}\label{e:zetabyparts}
  \z_{\th}(A_{i},s) = \frac{1}{(s-1)\ldots (s-m)}.\frac{i}{2\pi}\int_{C}\la^{m-s}\,
  \Tr(\dd_{\la}^{m}(A_i - \la)\ii) \, d\la  \ .
\end{equation}
yields the meromorphic continuation of $ \z_{\th}(A_{i},s)$ to
all of $\C$ with singularity structure
\begin{equation}\label{e:zetasing}
  \frac{\pi}{\sin(\pi s)}\z_{\th}(A_{i},s) \sim  \sum_{j=0}^{\infty}\sum_{k=0}^{1}
  \frac{\t{a}_{j,k}^{(i)}}{(s + \a_j - 1)^{k+1}} \ ,
\end{equation}
where, independently of m,
\begin{equation}\label{e:singcoeffs}
 \t{a}_{j,k}^{(i)} \approx \mu_k \G(\a_j)\G(\a_j + m )\ii
 a_{j,k}^{(i)}(m) \ ,
\end{equation}
$\mu_k = (1 + i(\th - \pi))^k$ and $\G(s)$ is the Gamma function
(see \cite{Se66}, \cite{GrSe96} Prop 2.9, here generalized to
arbitrary $\th$).

\vskip 2mm

\noi {\it Notation}: In equation \eqref{e:singcoeffs} $\approx$
indicates that
\begin{equation*}
\frac{\G(\a_j)\G(\a_j + m )\ii
 a_{j,0}^{(i)}(m)}{s+\a_j -1} + \frac{\mu_1 \G(\a_j)\G(\a_j + m )\ii
 a_{j,1}^{(i)}(m)}{(s+\a_j -1)^2}
\end{equation*}
gives the full pole structure at $s = 1-\a_j$. A function defined
in the sector $\Lambda_{\th}$ has an asymptotic expansion
\begin{equation*}
f(\la) \sim
  \sum_{j=0}^{\infty}\sum_{k=0}^{1} c_{j,k}(-\la)^{-\b_j}\log^k (-\la)
 + c(-\la)^{-\nu}
 \end{equation*}
as $\la\to\o$ with $ \b_j \nearrow +\o, \nu >0$ means that for
any $\ep >0$ and $N$ with $\b_N > \nu$,
\begin{equation*}
f(\la) =
  \sum_{j=0}^{N-1}\sum_{k=0}^{1} c_{j,k}(-\la)^{-\b_j}\log^k (-\la)
 + c(-\la)^{-\nu} + O(|\la|^{-\b_N + \ep}) \ ,
 \end{equation*}
for $\la$ sufficiently large, while a function $g$ on $\C$ has
singularity structure
\begin{equation*}
  g(s) \sim  \sum_{j=0}^{\infty}\sum_{k=0}^{1}
  \frac{d_{j,k}}{(s + \g_j - 1)^k}
\end{equation*}
means that
\begin{equation*}
  g(s) =   \sum_{j=0}^{N-1}\sum_{k=0}^{1}
  \frac{d_{j,k}}{(s + \g_j - 1)^k} + h_{N}(s) \ ,
\end{equation*}
with $h_N$ holomorphic for $1-\g_N < \re(s) < N + 1$.

\vskip 2mm

At any rate, \eqref{e:zetasing} implies that the term with
coefficient $a_{j,k}^{(i)}(m)$ in the resolvent trace expansion
\eqref{e:exp1} corresponds to a pole of $\frac{\pi}{\sin(\pi s)}
\z_{\th}(A_i,s)$ at $s = 1 - \a_j$ of order $k+1$. In particular,
since $\frac{\sin(\pi s)}{\pi} = s + O(s^3)$ around $s=0$, if
\begin{equation}\label{e:zeronopole}
\t{a}_{J,1}^{(1)} = 0 = \t{a}_{J,1}^{(2)} \ ,\hskip 10mm \a_J = 1
\ ,
\end{equation}
then \eqref{e:zetasing} implies the $\z_{\th}(A_i,s)$ have no
pole at $s=0$ and
\begin{equation}\label{e:zetazero1}
 \z_{\th}(A_i,0) = \t{a}_{J,0}^{(i)} = \frac{a_{J,0}^{(i)}(m)
 }{m!} \ .
\end{equation}
The regularity at $s=0$ means we can define the $\z$-determinants
\begin{equation*}
 {\rm det}_{\z,\th}A_1 = e^{-\z^{'}_{\th}(A_1,0)} \ , \hskip 15mm
 {\rm det}_{\z,\th}A_2 = e^{-\z^{'}_{\th}(A_2,0)} \ ,
\end{equation*}
where $\z^{'}_{\th} = d/ds(\z_{\th})$. If \eqref{e:zeronopole}
holds we refer to each of $A_1, A_2$ as $\z$-{\it admissible}.
Thus, for example, elliptic  $\pdo$s of order $d> 0$ over a closed
manifold of dimension $n$ are $\z$-admissible with
$a_{j,k}^{(i)}(m)$ locally determined, $\a_{j} = (j-n)/d$, so $J =
n+d$; for {\em differential} operators $a_{j,1}^{(i)}(m)=0$ (no
log terms). In the following we do not assume that the operators
$A_i$ are $\z$-admissible unless explicitly stated.

\vskip 4mm

\begin{defn}
We refer to a pair $(A_1, A_2)$ of invertible closed operators on
$H$ with spectral cut $R_{\th}$ as $\z$-{\it comparable} if for
$\la\in  \Lambda_{\th}$  :\newline \noi {\bf (I)} The relative
resolvent $(A_1 - \la)\ii - (A_2 - \la)\ii$ is trace class and
\begin{equation}\label{e:reltrace}
\Tr((A_1 - \la)\ii - (A_2 - \la)\ii) = -\frac{\dd}{\dd\la}\log{\rm
det}_{F}\Ss_{\la} \ .
\end{equation}
Here the `scattering' operator $\Ss_{\la} = \Ss_{\la}(A_1,A_2)$ is
an operator of the form $Id + W_{\la}$ on a Hilbert space
$\Hh_{\la}\subseteq H$ with $W_{\la}$ of trace class, so that
$\Ss_{\la}$ has a Fredholm determinant ${\rm det}_F\Ss_{\la} := 1
+ \sum_{k\geq 1}\Tr(\wedge^k W_{\la})$ taken on $\Hh_{\la}$.

 \noi
{\bf (II)} There is an asymptotic expansion as $\la\to\o$
\begin{equation}\label{e:exp3}
  \Tr ((A_1 - \la)\ii - (A_2 - \la)\ii)   \sim
  \sum_{\stackrel{j=0}{j\neq J}}^{\infty}\sum_{k=0}^{1}
  b_{j,k}(-\la)^{-\a_j}\log^k (-\la) +  b_{J,0}(-\la)^{-1} \ ,
\end{equation}
where  $0 < \a_0 < \ldots < \a_j \nearrow +\o$ and $\a_J = 1$.
\end{defn}

\begin{rem}\label{r:dz}

\vskip 2mm

\noi (1) If  $\Hh := \cup_{\la} \Hh_{\la}$ forms a (trivializable)
vector bundle, then the right-side of \eqref{e:reltrace} can be
written $-\Tr(\Ss_{\la}\ii \ d_{\la} \Ss_{\la}),$ where $d_{\la}$
is a covariant derivative on $\Hom(\Hh, H)$.  There is a canonical
choice for $d_{\la}$ induced from the covariant derivative
$\nabla_{\dd/\dd\la} = P(\la)\cdot(\dd/\dd\la)\cdot P(\la)$ on
$\Hh$, with $P(\la)$ the projection on $H$ with range $\Hh_{\la}$.

 (2) If an expansion \eqref{e:exp3} exists, then $\a_j > 0$
since $A_1\ii - A_2\ii$ is trace class.
\end{rem}

\vskip 3mm

If $A_1, A_2$ are $\z$-comparable, then $A_1\si - A_2\si$ is trace
class for $\re(s)> 1$. Hence we define the {\it relative spectral
$\z$-function} by
\begin{equation}\label{e:defrelz}
  \z_{\th}(A_{1},A_{2},s)  =  \Tr (A_1\si - A_2\si) \ , \hskip 15mm \re(s) >  1 \ .
\end{equation}
In view of \eqref{e:reltrace} we have
\begin{eqnarray}
  \z_{\th}(A_{1},A_{2},s) & = & \frac{i}{2\pi}\int_{C}\la\si\,
  \Tr((A_1 - \la)\ii - (A_2 - \la)\ii) \, d\la  \nonumber \\
  & = & -\frac{i}{2\pi}\int_{C}\la\si\,
  \frac{\dd}{\dd\la}\log {\rm det}_F \Ss_{\la} \, d\la \ . \label{e:relzetS}
\end{eqnarray}
$\z_{\th}(A_{1},A_{2},s)$  thus extends holomorphically to
$\re(s)> 1 -\a_0 $ while the asymptotic expansion \eqref{e:exp3}
defines the meromorphic continuation to $\C$ with singularity
structure
\begin{equation}\label{e:relzetasing}
  \G(s)\z_{\th}(A_1,A_2,s) \sim  \sum_{\stackrel{j=0}{j\neq J}}^{\infty}\sum_{k=0}^{1}
  \frac{\t{b}_{j,k}}{(s + \a_j - 1)^k} + \frac{\t{b}_{J,0}}{s}\ ,
\end{equation}
where
\begin{equation}\label{e:singcoeffs1}
  \t{b}_{j,k} = \mu_k\G(\a_j)\ii b_{j,k} , \hskip 10mm
  \t{b}_{J,0}= b_{J,0} = \z_{\th}(A_1,A_2,0) \ .
\end{equation}
Since $b_{J,1} = 0$ in \eqref{e:exp3} then
$\z_{\th}(A_{1},A_{2},s)$ is regular at $s=0$ and we can define
the {\it relative $\z$-determinant} by
\begin{equation}\label{e:relzdet}
{\rm det}_{\z,\th}(A_1,A_2) = e^{-\z^{'}_{\th}(A_1,A_2,0)} \ .
\end{equation}
No assumption is made on the existence or regularity of
$\z_{\th}(A_i ,s)$. If $A_1, A_2$ are $\z$-admissible and (I) of
Definition (2.1) holds, then \eqref{e:exp1}, \eqref{e:zeronopole}
imply  as $\la\to\o$ in $\Lambda_{\th}$
\begin{equation}\label{e:exp2}
  \dd_{\la}^{m}\Tr ((A_1 - \la)\ii - (A_2 - \la)\ii) \sim
\end{equation}
$$  \sum_{\stackrel{p=0}{p\neq J}}^{\infty}\sum_{k=0}^{1}
  (a_{p,k}^{(1)}(m) - a_{p,k}^{(2)}(m))(-\la)^{-\a_p
  - m}\log^k (-\la) + (a_{J,0}^{(1)}(m) - a_{J,0}^{(2)}(m))(-\la)^{-1-m} \ ,$$
  with, by $\z$-comparability,
$\a_p > 0$---for, the relative resolvent
 trace is then $O(|\la|^{-\e})$ as $\la\to\o$, some $\e>0$,
and hence $a_{j,k}^{(1)}(m) - a_{j,k}^{(2)}(m) = 0$ in
\eqref{e:exp1} for $\a_j\leq 0$, while in \eqref{e:exp2}
$p=j-\max\{j \ | \ \a_j\leq 0\} + 1$,
 resulting in the regularity of
$\z_{\th}(A_{1},A_{2},s)$ for $\re(s)> 1 - \e$. With $\a_p > 0$ we
can integrate \eqref{e:exp2} to obtain for $\la\to\o$ an
asymptotic expansion of the form \eqref{e:exp3}.

The $b_{j,k}$ are related to the coefficients in \eqref{e:exp2}
via universal constants; in \eqref{e:relzetasing}
\begin{equation}\label{e:singcoeffs1a}
  \t{b}_{j,k} \approx
  \t{a}_{j,k}^{(1)} - \t{a}_{j,k}^{(2)} \ ,
\end{equation}
and in particular
\begin{equation}\label{e:singcoeffs2}
 \z_{\th}(A_{1},A_{2},0) = \t{b}_{J,0}=  \t{a}_{J,0}^{(1)} -
  \t{a}_{J,0}^{(2)} =  \z_{\th}(A_{1},0) - \z_{\th}(A_{2},0) \ .
\end{equation}

\noi More generally:
\begin{lem}\label{lem:relz=zrel}
If $A_1, A_2$ are $\z$-admissible operators such that (I) of
Definition (2.1) holds, then $A_1, A_2$ are $\z$-comparable and as
meromorphic functions on $\C$
\begin{equation}\label{e:relz=zrel}
  \z_{\th}(A_{1},A_{2},s) =  \z_{\th}(A_{1},s) - \z_{\th}(A_{2},s) \ .
\end{equation}
\end{lem}
\begin{proof}
The first statement is proved above. For $\re(s)> 1 -\a_0$,
\eqref{e:relz=zrel} is obvious. Elsewhere, from
\eqref{e:singcoeffs1}, or \eqref{e:exp2}, \eqref{e:exp3}, the
left and right sides of \eqref{e:relz=zrel} have the same
singularity structure, hence $ \z_{\th}(A_{1},A_{2},s) -
\z_{\th}(A_{1},s) + \z_{\th}(A_{2},s)$ is a holomorphic
continuation of zero from $\re(s)> 1 -\a_0$ to all of $\C$, and is
therefore identically zero.
\end{proof}

To compute $\det_{\z,\th}(A_1,A_2)$ in terms of the scattering
matrix we need to know more about the asymptotic behaviour of
$\Ss_{\la}$.
\begin{lem}\label{lem:asympSla}
Let $f$ be a differentiable function in the sector
$\Lambda_{\th}$ with an asymptotic expansion as $\la\to \o$
\begin{equation}\label{e:deriv-exp}
  -\frac{\dd f}{\dd\la}   \sim
  \sum_{j=0}^{\infty}\sum_{k=0}^{1} c_{j,k}(-\la)^{-\b_j}\log^k (-\la)
 + c_0 (-\la)\ii
\end{equation}
with $\b_j\nearrow +\infty$ and $\b_j\neq 1$. Then
\begin{equation}\label{e:cf}
  c_f (\la) := f(\la) -
  \sum_{j=0}^{r}\sum_{k=0}^{1}
  \frac{c_{j,k}}{1-\b_j}\lrb (-\la)^{-\b_j +1}\log^k (-\la)
  -  \frac{k}{1-\b_j}(-\la)^{-\b_j +1}\rrb  - c_0 \log(-\la),
\end{equation}
where $r = \max\{k \ | \ \b_k < 1 \}$,  converges uniformly as
$\la\to\o$. Denoting this limit by $$c_1 := {\rm
lim}^{\th}_{\la\to\o} c_f (\la) \ ,$$ (the $\th$ indicating the
limit is taken in the sector $\Lambda_{\th}$), there is an
asymptotic expansion as $\la\to\o$ in $\Lambda_{\th}$
\begin{equation}\label{e:f-exp}
  f(\la) \sim
  \sum_{j=0}^{\o}\sum_{k=0}^{1}
  \frac{c_{j,k}}{1-\b_j}\lrb (-\la)^{-\b_j +1}\log^k (-\la)
  -  \frac{k}{1-\b_j}(-\la)^{-\b_j +1}\rrb  + c_0 \log(-\la) + c_1 \ .
\end{equation}
\end{lem}
\begin{proof}
Let $\la_0 \in R_{\varphi}$  with $\varphi=\arg(\la)$ and $|\la|
< |\la_0|$. Choosing  $|\la|$ sufficiently large so that
\eqref{e:deriv-exp} holds, then
\begin{equation*}
|c_f(\la) - c_f(\la_0)| \leq \int_{|\la_0|}^{|\la|}
\left|\frac{\dd c_f }{\dd\mu}\right| \ d\mu
 = \int_{|\la_0|}^{|\la|} |\frac{\dd f}{\dd\mu} +
 \sum_{j=0}^{r}\sum_{k=0}^{1} c_{j,k}(-\mu)^{-\b_j}\log^k (-\mu)
 + c_0 (-\mu)\ii|  \ d\mu
\end{equation*}
\begin{equation*}
 =
  \int_{|\la_0|}^{|\la|} O(\mu^{-\b_{r+1}}) \ d\mu \leq
  C|\la|^{-\b_{r+1}+ 1} \ .
\end{equation*}
Since $\b_{r+1} >1$, then $(c_f (\la))$ is convergent by the
Cauchy criterion.

Integrating \eqref{e:deriv-exp} between $\la$ and $\la_0$, we
obtain for large $|\la|,|\la_0|$ and any $\ep
>0$
\begin{equation*}
  f(\la) =
  \sum_{j=0}^{N-1}\sum_{k=0}^{1}
  \frac{c_{j,k}}{1-\b_j}\lrb (-\la)^{-\b_j +1}\log^k (-\la)
  -  \frac{k}{1-\b_j}(-\la)^{-\b_j +1}\rrb
\end{equation*}
$$
 + c_0 \log(-\la) +
  O(|\la|^{-\b_N + 1 + \ep}) + c_f (\la_0) +
  O(|\la_0|^{-\b_N + 1 + \ep}) \ .
$$ Letting $\la_0\to\o$ we reach the conclusion.
\end{proof}

Applying \lemref{lem:asympSla} to $\f(\la) = \log\det_F \Ss_{\la}$
and from \eqref{e:reltrace}, \eqref{e:exp3}, we see that as
$\la\to\o$ in $\Lambda_{\th}$ there is an asymptotic expansion
\begin{equation}\label{e:logdetSla}
  \log{\rm det}_F \Ss_{\la} \sim
  \sum_{\stackrel{j=0}{j\neq J
}}^{\o}\sum_{k=0}^{1}
  b^{'}_{j,k} (-\la)^{-\a_j +1}\log^k (-\la)
+ b_{J,0}\log(-\la) + c_{{\rm rel}} \ ,
\end{equation}
where $b^{'}_{j,0} = b_{j,0} (1-\a_j)\ii - b_{j,1} (1-\a_j)^{-2}$,
$b^{'}_{j,1} = b_{j,1} (1-\a_j)\ii$ ($j\neq J$), and the constant
term is
\begin{equation}\label{e:const}
c_{{\rm rel}}   = {\rm lim}^{\th}_{\la\to\o} [ \log{\rm det}_F
\Ss_{\la} - b_{J,0}\log(-\la) -
\sum_{j=0}^{J-1}\sum_{k=0}^{1}b^{'}_{j,k} (-\la)^{-\a_j +1}\log^k
(-\la)] \ .
\end{equation}

The regularized limit of a function in the sector $\Lambda_{\th}$
with an asymptotic expansion $f(\la)\sim
\sum_{j=0}^{\infty}\sum_{k=0}^{1}\textsf{c}_{jk}(-\la)^{-\b_j}\log^{k}(-\la)
 + \textsf{c}_0\log(-\la) + \textsf{c}_1$ as $\la\to\o$,
 where  $\b_j\nearrow +\infty$ and
$\b_j\neq 0$, picks out the constant term in the expansion
\begin{equation*}
 \LIM^{\th}_{\la\to\infty}f(\la) = \textsf{c}_{1} \ .
\end{equation*}
We have (with $\Ss := \Ss_0$):

\vskip 4mm

\begin{thm}\label{t:thmB} For $\z$-comparable operators $A_1,A_2$
\begin{equation}\label{e:thm3}
{\rm det}_{\z,\th}(A_1,A_2) = {\rm det}_F \Ss\, . \,
 e^{-{\rm LIM^{\th}_{\la\to\infty}}\log{\rm det}_F\Ss_{\la}} \; .
\end{equation}
With $\z_{{\rm rel}}(0) := \z_{\th}(A_1,A_2,0)$, one has
\begin{equation}\label{e:reglim}
{\rm LIM^{\th}}_{\la\to\infty}\log {\rm det}_F \Ss_{\la} = {\rm
lim}^{\th}_{\la\to\o} [ \log{\rm det}_F \Ss_{\la} - \z_{{\rm
rel}}(0)\log(-\la)
\end{equation}
$$ - \sum_{j=0}^{J-1}\sum_{k=0}^{1}b^{'}_{j,k} (-\la)^{-\a_j
+1}\log^k (-\la)] \ .$$ If $A_1,A_2$ are $\z$-admissible
\begin{equation}\label{e:reldet=ratio}
{\rm det}_{\z,\th}(A_1,A_2) = \frac{{\rm det}_{\z,\th}A_1}{{\rm
det}_{\z,\th}A_2} \ .
\end{equation}
\end{thm}
\vskip 3mm
\begin{proof}
The identity \eqref{e:reldet=ratio} is immediate from
\lemref{lem:relz=zrel}, while \eqref{e:reglim} follows from
\eqref{e:logdetSla}, \eqref{e:const}, and
\eqref{e:singcoeffs1},\eqref{e:singcoeffs2}.

To prove \eqref{e:thm3}, since $\la\si\log\det_F\Ss_{\la}\to 0$ at
the ends of $C$ for $\re(s)>1-\a_0$, we can integrate by parts in
\eqref{e:relzetS} to obtain
\begin{equation}\label{e:relzeta-g}
\z_{\th}(A_1,A_2,s) = sg(s)
\end{equation}
and hence that
\begin{equation}\label{e:relzetaderiv}
\z_{\th}^{'}(A_1,A_2,0) = \frac{d}{ds}_{|s=0}(sg(s)|^{{\rm mer}})
\ ,
\end{equation}
where, with $f(\la) = \frac{\log{\rm det}_F \Ss_{\la}}{(-\la)}$,
\begin{equation}\label{e:g(s)}
 g(s) = \frac{i}{2\pi}\int_{C}\la\si \, f(\la)  d\la \,
\end{equation}
has a simple pole at $s=0$ with residue $b_{J,0}$. The notation
$h(s)|^{{\rm mer}}$ indicates the meromorphically continued
function.

We carry out the meromorphic continuation of $\z_{\th}(A_1,A_2,s)$
along the lines of \cite{GrSe96} Prop 2.9 . First, $\log{\rm
det}_F \Ss_{\la}$ is regular near $\la = 0$ and so $f(\la)$ is
meromorphic there with Laurent expansion
\begin{equation}\label{e:around0}
  f(\la) = \frac{\log{\rm det}_F \Ss}{(-\la)} + \sum_{j=0}^{\o}
  b_j (-\la)^j \ .
\end{equation}
Since $\int_{C}\la^{-1-s} d\la = 0$ for $\re(s)>0$, then
\begin{equation}\label{e:g(s)2}
 g(s) = \frac{i}{2\pi}\int_{C}\la\si \, f_0(\la)  d\la  \ ,
 \hskip 10mm f_0 (\la) := f(\la) - \frac{\log{\rm det}_F
 \Ss}{(-\la)} \ .
\end{equation}
For $1-\a_0 < \re(s)<1$, the circular part $C_{\rho,\th,\th-2\pi}$
of the contour $C$ can now be shrunk to the origin, which reduces
$g(s)$ to
\begin{equation}\label{e:g(s)3}
 g(s) = \frac{\sin(\pi s)}{\pi}e^{i(\pi - \th)(s-1)}
 \int_{0}^{\o} r\si \, f_0(re^{i\th}) \
 dr \ ,
\end{equation}
using $e^{-is(\th - 2\pi)} - e^{-is\th} = 2i\sin(\pi s).e^{is(\pi
- \th)}$. On the other hand, from \eqref{e:logdetSla} there is an
asymptotic expansion as $\la\to\o$ along $R_{\th}$
\begin{equation}\label{e:f0-asymp}
  f_0 (\la) \sim c_0 \frac{\log(-\la)}{(-\la)} + \frac{c_1 - \log{\rm det}_F
 \Ss}{(-\la)} + \sum_{\stackrel{j=0}{j\neq J
}}^{\o}\sum_{k=0}^{1} b^{'}_{j,k}(-\la)^{-\a_j}\log^k (-\la)
\end{equation}
where $c_0 = b_{J,0} = \z_{rel}(0)$, $c_1 = c_{{\rm rel}}$.
Hence, since $-\la = -re^{i\th} = re^{i(\th-\pi)}$ with respect
to $R_{\th}$, for any $\ep > 0, N> J+1$ we have as $r\to\o$
\begin{eqnarray}
  f_0 (re^{i\th}) & = & c_0 \frac{\log(re^{i(\th-\pi)})}{re^{i(\th-\pi)}}
  + \frac{c_1 - \log{\rm det}_F
 \Ss}{re^{i(\th-\pi)}} \label{e:f0r-asymp} \\
& + &
 \sum_{\stackrel{j=0}{j\neq J
}}^{N-1}\sum_{k=0}^{1} b^{'}_{j,k} r^{-\a_j}e^{i(\th-\pi)\a_j}
 \log^k (re^{i(\th-\pi)}) + O(r^{-\a_N + \ep}). \nonumber
\end{eqnarray}

\noi Therefore
\begin{equation*}
  e^{-i(\th-\pi)}e^{-is(\pi-\th)}\frac{\pi}{\sin(\pi s)}g(s) =
  \int_{0}^{1}\lsb\sum_{j=0}^{N-1} b_j e^{i(\th-\pi)j}r^{j-s} +
  r^{-s}O(r^N ) \rsb \ dr
\end{equation*}
$$ + \int_{1}^{\o}\lsb e^{-i(\th-\pi)}c_0 r^{-s-1}\log(r) +
e^{-i(\th-\pi)}(c_1 + i(\th-\pi)c_0 - \log{\rm det}_F
 \Ss) r^{-s-1} \right.$$
$$ + \sum_{\stackrel{j=0}{j\neq J }}^{N-1}\sum_{k=0}^1 \left.
c_{j,k,\th} r^{-\a_{j}-s}\log^k (r) +  r^{-s}O(r^{-\a_N + \ep})
\rsb \ dr$$
\begin{equation*}
 = -\sum_{j=0}^{N-1} \frac{b_j e^{i(\th-\pi)j}}{s-j-1} +
+ \frac{e^{-i(\th-\pi)}c_0}{s^2} + \frac{e^{-i(\th-\pi)}(c_1 +
i(\th-\pi)c_0 - \log{\rm det}_F \Ss)}{s}
\end{equation*}
$$ + \sum_{\stackrel{j=0}{j\neq J }}^{N-1}\sum_{k=0}^1
\frac{c_{j,k,\th}}{(s+\a_j -1)^k} + h_N (s) \ ,$$ where $h_N$ is
holomorphic for $1- \a_N + \ep < \re(s) < N+1$. Here we use the
meromorphic extension to $\C$ of
\begin{equation}\label{e:int1}
 \int_{0}^{1} r^{j-s} \ dr = \frac{-1}{s-j-1} \ , \hskip 20mm \re(s)
 < j+1 \ ,
\end{equation}
\begin{equation}\label{e:int2}
 \int_{1}^{\o} r^{\b-s}\log^k (r) \ dr = \frac{1}{(s-\b -1)^{k+1}} \ ,
 \hskip 10mm k=0,1 \ , \hskip 5mm \re(s)
 > \b +1 \ .
\end{equation}

This implies the singularity structure
$$e^{-i(\pi-\th)s}\frac{\pi}{\sin(\pi s)}g(s) \sim \frac{c_0}{s^2}
+ \frac{c_1 + i(\th-\pi)c_0 - \log{\rm det}_F \Ss}{s}  -
\sum_{j=0}^{\o} \frac{b_{j}}{s-j-1}$$
\begin{equation}\label{e:poles}
 + \sum_{\stackrel{j=0}{j\neq J
}}^{\o}\sum_{k=0}^1 \frac{e^{i(\th-\pi)j} c_{j,k,\th}}{(s+\a_j
-1)^k}\ ,
\end{equation}
Around $s=0$ one has $\frac{\sin(\pi s)}{\pi} = s + O(s^3)$ and
hence
\begin{equation}\label{e:sgpoles}
s g(s) = e^{i(\pi-\th)s}(s^2 + O(s^4)) \lrb \frac{c_0}{s^2} +
\frac{c_1 + i(\th-\pi)c_0 - \log{\rm det}_F \Ss}{s}\rrb + s^2
p(s) \ ,
\end{equation}
where $p$ is meromorphic on $\C$ and holomorphic around $s=0$,
giving the pole structure in \eqref{e:poles} away from the
origin. We therefore have near $s=0$
\begin{equation*}
\frac{d}{ds}(s g(s)) = i(\pi - \th) e^{i(\pi-\th)s}\left(c_0 +
s(c_1 + i(\th-\pi)c_0 - \log{\rm det}_F \Ss)\right)
\end{equation*}
$$  + O(s^2) + e^{i(\pi-\th)s}(c_1 + i(\th-\pi)c_0 - \log{\rm
det}_F \Ss) + O(s) \ .$$ And hence from \eqref{e:relzetaderiv}
\begin{equation}\label{e:thm3a}
\z_{\th}^{'}(A_1,A_2,0) = c_1  - \log{\rm det}_F \Ss \
\end{equation}
and this is equation \eqref{e:thm3}.
\end{proof}

\begin{rem}\label{rem:freedom}
There is freedom in specifying $\log\det_F \Ss_{\la}$ up to the
addition of a constant, and hence in specifying $\Ss_{\la}$ up to
composition with an element of $\GL_{1}(\Hh_{\la}) = \{\Ee\in
GL(\Hh_{\la}) \ | \ \Ee-I\in L_{1}(\Hh_{\la})\}$, where
$L_{1}(\Hh_{\la})$ is the ideal of trace-class operators; that
is, $\Ss_{\la}\Ee$ is also a scattering operator for any $\Ee\in
GL_{1}(\Hh_{\la})$. However, since $\det_F : GL_{1}(\Hh_{\la})
\to \C^*$ is a group homomorphism and the regularized limit is
linear
\begin{equation}\label{e:linLIM}
\LIM^{\th}_{\la\to\infty}(g(\la) +  c \, . \, f(\la)) =
\LIM^{\th}_{\la\to\infty}(g(\la)) + c \, . \,
\LIM^{\th}_{\la\to\infty}(f(\la)) \ ,
\end{equation}
any  constant $c$, then \eqref{e:thm3} and \eqref{e:thm3a} are
unambiguous.
\end{rem}

\bigskip

With the regularized limit $\LIM_{z\to 0}(h(z))$ denoting the
constant term in the Laurent expansion of a function $h(z)$ around
$z=0$, we can recast \eqref{e:thm3} as follows:

\begin{prop}\label{p:detaroundzero} If $A_1, A_2$ are
$\z$-comparable, then
\begin{eqnarray}
\log{\rm det}_{\z,\th}(A_1,A_2) & = & -\LIM_{s\to
0}\left[\left.\frac{i}{2\pi}\int_{C}\la^{-s-1}\log{\rm det}_F
\Ss_{\la} \ d\la \right|^{mer} \right] \nonumber \\
& = & -\left[\left.\frac{i}{2\pi}\int_{C}\la^{-s-1}\log{\rm
det}_F \Ss_{\la} \ d\la  - \frac{
\z_{rel,\th}(0)}{s}\right]\right|^{mer}_{s=0}
\label{e:detaroundzero1}
\end{eqnarray}
Equivalently,
\begin{eqnarray}
\log{\rm det}_{\z,\th}(A_1,A_2) & = & -\LIM_{s\to 0}\left[\left.
\G(s) \z_{\th}(A_1,A_2,s)\right|^{mer} \right] + \g
\z_{rel,\th}(0) \nonumber\\
& = & -\left[\left. \G(s)\z_{\th}(A_1,A_2,s) - \frac{
\z_{rel,\th}(0)}{s} \right]\right|^{mer}_{s=0} + \g
\z_{rel,\th}(0)
 \ . \label{e:detaroundzero2}
\end{eqnarray}
\end{prop}
\begin{proof}
From \eqref{e:g(s)},\eqref{e:sgpoles}, around $s=0$ one has
\begin{eqnarray*}
 \frac{i}{2\pi}\int_{C}\la^{-s-1}\, \log{\rm
det}_F\Ss \ d\la  & = &
 e^{i(\pi-\th)s}\left(\frac{c_0}{s} + c_1 + i(\th -
\pi)c_0 - \log{\rm det}_F\Ss + O(s) \right) \\
& = &
 \frac{c_0}{s} + c_1 - \log{\rm det}_F\Ss + O(s) \ ,
\end{eqnarray*}
and hence \eqref{e:detaroundzero1} follows from \eqref{e:thm3a}.

On the other hand,  $\G(s) = s\ii + \g s + O(s)$ near $s=0$, so
from \eqref{e:relzeta-g}
\begin{eqnarray*}
\G(s) \z_{\th}(A_1,A_2,s) & = & \G(s)sg(s) \\
& = & e^{i(\pi-\th)s}(\frac{c_0}{s} + c_1 + i(\th-\pi)c_0 -
\log{\rm det}_F\Ss  + c_0\g + O(s)) \ ,
\end{eqnarray*}
and so \eqref{e:detaroundzero2} follows similarly.
\end{proof}

We also have:

\begin{prop}\label{p:detasymp} Let $A_1, A_2$ be
$\z$-comparable and $\z$-admissible, then there is an asymptotic
expansion as $\mu\to\o$ in $\Lambda_{\th}$
\begin{equation*}
\log{\rm det}_{\z,\th}(A_1 - \mu) - \log{\rm det}_{\z,\th}(A_2 -
\mu) \sim
\end{equation*}
$$\sum_{\stackrel{j=0}{j\neq J }}^{\o}\sum_{k=0}^{1}b^{'}_{j,k}
(-\la)^{-\a_j +1}\log^k (-\la) + b_{J,0}\log(-\la) \ .$$
  In particular, the constant term is zero : ${\rm LIM}^{\th}_{\mu\to
\o}\log\left(\frac{{\rm det}_{\z,\th}(A_1 - \mu)}{{\rm
det}_{\z,\th}(A_2 - \mu)}\right) = 0 \ .$
\end{prop}
\begin{proof}
For $\mu\in \Lambda_{\th}$ the operators $A_i - \mu$ are
$\z$-comparable, and hence by \thmref{t:thmB}
\begin{equation*}
\log{\rm det}_{\z,\th}(A_1 - \mu, A_2 - \mu) = \log{\rm
det}_{F}\Ss_{\mu} - {\rm LIM^{\th}_{\la\to\infty}}\log{\rm
det}_{F}\Ss_{\mu+\la} \ .
\end{equation*}
Since ${\rm LIM^{\th}_{\la\to\infty}}\log{\rm
det}_{F}\Ss_{\mu+\la} = {\rm LIM^{\th}_{\la\to\infty}}\log{\rm
det}_{F}\Ss_{\la}$ the conclusion is reached from
\eqref{e:logdetSla},\eqref{e:f0-asymp}.
\end{proof}

Finally, it is useful to note that a similar proof allows
\thmref{t:thmB} to be abstracted and generalized slightly:
\begin{prop}\label{p:thmB*}
Let $\Phi$ be a function on $\C$ which is meromorphic at $0$ with
Laurent expansion $\Phi(\la) = \sum_{j=-m}^{\o}b_j (-\la)^j$, and
holomorphic in a sector $\Lambda_{\th}$ with for some $r\in\Z$ a
uniform asymptotic expansion of $\dd\Phi/\dd \la$ as $\la\to\o$ in
$\Lambda_{\th}$
\begin{equation}\label{e:Phi-exp}
\frac{\dd\Phi}{\dd\la} \sim
\sum_{j=-r}^{\o}\sum_{k=0}^{1}a_{j,k}(\la)^{-\a_j}\log^k(-\la) \ ,
\end{equation}
where $-\o < \a_{-r} < \ldots < \a_{-1} \leq 0 < \a_0 < \ldots <
\a_{j}\nearrow \o.$ Then
\begin{equation*}
Z(s) = -\frac{i}{2\pi}\int_{C}\la\si \frac{\dd\Phi}{\dd\la} \ d\la
\end{equation*}
is holomorphic for $\re(s) > 1-\a_{-r}$ and has a meromorphic
continuation to all of $\C$ with poles determined by the
coefficients of \eqref{e:Phi-exp}. In particular, if $a_{J,1} =
0$ (with $\a_J = 1$) then $Z(s)$ is holomorphic around $s=0$ with
$Z(0) = a_{J,0}$, and
\begin{equation}\label{e:Zderiv}
  -Z^{'}(0) = \Phi(0) -\LIM^{\theta}_{\la\to\o}\Phi(\la) \ ,
\end{equation}
where $\Phi(0) := {\rm LIM}_{\la\to 0}\Phi(\la) = b_0$.
\end{prop}

\subsection{Relative Heat Kernel Regularization}

In this Section we derive \thmref{t:thmB} using the heat operator
trace for operators $A_1, A_2$ with spectrum contained in a
sector of the right-half plane. This applies, for example, to the
Dirac Laplacian on a manifold with boundary.

We assume that $R_{\pi}$ is a spectral cut for $A_1, A_2$, so
that $\|(A_i - \la)\ii\| = O(|\la|\ii)$ for large $\la$ in
$\Lambda_{\pi}$ with $\d > \pi/2$ in \eqref{e:sector}. Let $\Cc$
be a contour surrounding ${\rm sp}(A_1), {\rm sp}(A_2)$, coming in
on a ray with argument in $(0,\pi/2)$, encircling the origin, and
leaving on a ray with argument in $(-\pi/2,0)$. For $t>0$, one
then has the heat operators
\begin{equation*}\label{e:heat}
  e^{-tA_i} := \frac{i}{2\pi}\int_{\Cc} e^{-t\la}(A_i - \la)\ii \
  d\la  = \frac{i}{2\pi}\int_{\Cc} t^{-m} e^{-t\la}\dd_{\la}^{m}(A_i - \la)\ii  \
  d\la \ .
\end{equation*}
If we assume $\dd_{\la}^{m}(A_i - \la)\ii$ is trace class for some
$m$, then $e^{-tA_i}$ is trace class with
 $$  \Tr ( e^{-tA_i} )  = \frac{i }{2\pi}\int_{\Cc}
  t^{-m} e^{-t\la}\Tr(\dd_{\la}^m(A_i - \la)\ii) \
  d\la \ .$$  The resolvent trace expansions \eqref{e:exp1} thus imply
  heat trace expansions as $t \to 0$
\begin{equation}\label{e:heatrace}
  \Tr ( e^{-tA_i} )   \sim
  \sum_{j=0}^{\infty}\sum_{k=0}^{1} \t{a}_{j,k}^{(i)} t^{\a_j
  -1}\log^k t \ ,
\end{equation}
with coefficients differing from those in \eqref{e:exp1} by
universal constants, while $\Tr ( e^{-tA_i} ) = O(e^{-ct})$, some
$c>0$, as $t\to\o$. The heat representation of the power operators
\begin{equation}\label{e:heatpower}
  A_{i}\si = \frac{1}{\G(s)}\int_{0}^{\o}t^{s-1}e^{-tA_i} \ dt \ ,
  \hskip 10mm \re(s) > 0 \ ,
\end{equation}
then implies $\z(A_i,s) =
\G(s)\ii\int_{0}^{\o}t^{s-1}\Tr(e^{-tA_i}) \ dt$ for $\re(s)>
1-\a_0$,  with \eqref{e:heatrace} giving the pole structure of the
meromorphic extension to $\C$.

For positive $\z$-admissible operators the heat cut-off
regularization, defined by
\begin{equation}\label{e:heatcuttoff}
  \log{\rm det}_{heat}(A_i) :=
  \LIM_{\e\to 0}\int_{\e}^{\o} - \frac{1}{t} \Tr(e^{-tA_i}) \ dt \ ,
\end{equation}
picks out the constant term in the asymptotic expansion as $\e\to
0$ of
\begin{equation*}
f(\e):=\int_{\e}^{\o} - \frac{1}{t} \Tr(e^{-tA_i}) \ dt \sim
-\z(A_i,0)\log \e + \LIM_{\e\to 0}f(\e) -
\sum_{j=0}^{\o}\sum_{k=0}^{1}f_{j,k}\e^{\a_j-1}\log^{k}\e \ .
\end{equation*}
Since $\dd f/\dd\e=\e\ii\Tr(e^{-\e A_i})$, the existence of such
an expansion follows from \eqref{e:heatrace} and the small time
asymptotics analogue of \lemref{lem:asympSla}.

The definition in \eqref{e:heatcuttoff} is motivated by
$\det_{heat}(A) = \det_{F}(A)$ in finite dimensions. However, if
$H$ is infinite-dimensional and $A$ is determinant class, then
$e^{-tA}$ is not trace class and \eqref{e:heatcuttoff} is
undefined. There is, nevertheless, for any pair of
$\z$-comparable operators $A_1, A_2$  a well-defined relative
heat cut-off determinant
\begin{equation}\label{e:relheatcuttoff}
  \log{\rm det}_{heat}(A_1,A_2) :=
  \LIM_{\e\to 0}\int_{\e}^{\o} - \frac{1}{t} \Tr(e^{-tA_1} - e^{-tA_2}) \ dt \ .
\end{equation}
This includes both when $A_1, A_2$ are determinant class or
$\z$-admissible. In the former case, $e^{-tA_1} - e^{-tA_2}$ is
now trace class, $\log{\rm det}_{heat}(A_1,A_2) = \int_{0}^{\o} -
\frac{1}{t} \Tr(e^{-tA_1} - e^{-tA_2}) \ dt$ and ${\rm
det}_{heat}(A_1,A_2) = \det_F (A_1)/\det_F (A_2)$. This is the
$r$-integrated version of
\begin{equation}\label{e:why}
\Tr(A_r\ii \dot{A}_r) =
  \int_{0}^{\o} -t\ii \Tr(\frac{d}{dr}e^{-tA_r}) \ dt \ , \hskip
  15mm (A_r \;\; {\rm det \; class})
\end{equation}
where $A_1 \leq A_r\leq A_2$ is a smooth 1-parameter family of
 non-negative determinant class operators (to see \eqref{e:why}, set $s=1$ in
\eqref{e:heatpower} and use Duhamel's principle).

On the other hand, when the $A_i$ are $\z$-admissible
\begin{equation}\label{e:relheatmult}
  {\rm det}_{heat}(A_1,A_2) = \frac{{\rm det}_{heat}(A_1)}{{\rm
  det}_{heat}(A_2)} \ .
\end{equation}
By \eqref{e:zeta-heat}, below, \eqref{e:relheatmult} is a
restatement  of \eqref{e:reldet=ratio} and \eqref{e:singcoeffs2}.

The relative heat and $\z$-function $\z(A_1,A_2,s) =
\G(s)\ii\int_{0}^{\o}t^{s-1}\Tr(e^{-tA_1} - e^{-tA_2}) \ dt$
regularizations are related to the scattering determinant is as
follows.

\vskip 4mm

\begin{thm}\label{t:thmB*}  Let $A_1, A_2$ be positive $\z$-comparable
operators with $\th = \pi$, as above. Then as $\e\to 0 $ there is
an asymptotic expansion
\begin{equation}\label{e:heatcuttoffexp}
\int_{\e}^{\o} - \frac{1}{t} \Tr(e^{-tA_1} - e^{-tA_2}) \ dt \sim
\log{\rm det}_F \Ss - \LIM_{\mu\to\o}\log{\rm det}_F \Ss_{-\mu} -
\z_{rel}(0)\G^{'}(1)
\end{equation}
$$\hskip 40mm - \z_{rel}(0)\log \e + \sum_{j=0}^{\o}\sum_{k=0}^1
c^{'}_{j,k}\e^{\a_j - 1}\log^k \e \
  . $$
Hence
\begin{equation}\label{e:logrelheat}
\log{\rm det}_{heat}(A_1,A_2) =  \log{\rm det}_F \Ss -
\LIM_{\la\to\o}\log{\rm det}_F \Ss_{-\la} - \z_{rel}(0)\G^{'}(1) \
.
\end{equation}
One has
\begin{equation}\label{e:zeta-heat}
{\rm det}_{\z}(A_1,A_2) =  {\rm det}_{heat}(A_1,A_2)
e^{\z_{rel}(0)\G^{'}(1)} \ .
\end{equation}
\end{thm}
\begin{rem}
$\G^{'}(1) = -\g$.
\end{rem}

\vskip 3mm

\begin{proof}
Equation \eqref{e:zeta-heat} follows from \eqref{e:logrelheat}
and \eqref{e:thm3}. Alternatively, it is proved directly, without
reference to \thmref{t:thmB}, by an obvious modification of the
following proof of \eqref{e:heatcuttoffexp}.

From \eqref{e:reltrace}, \eqref{e:heat} we have

\begin{equation}\label{e:scatterheatrace}
\Tr(e^{-tA_1} - e^{-tA_2}) =  -\frac{i}{2\pi}\int_{\Cc}
e^{-t\la}\frac{\dd}{\dd\la}\log{\rm det}_F \Ss_{\la} \ d\la =
-t\frac{i}{2\pi}\int_{\Cc} e^{-t\la}\log{\rm det}_F \Ss_{\la} \
d\la \ .
\end{equation}
Hence $$\int_{\e}^{\o}-\frac{1}{t}\Tr(e^{-tA_1} - e^{-tA_2}) \ dt
= \frac{i}{2\pi}\int_{\Cc} \lim_{\omega\to
\o}\left.\frac{e^{-t\la}}{-\la}\right|^{\omega}_{t=\e}\log{\rm
det}_F \Ss_{\la} \ d\la$$
\begin{equation}\label{e:ecuttoff}
= \frac{i}{2\pi}\int_{\Cc} \lim_{\omega\to
\o}\frac{e^{-\mu}\log{\rm det}_F \Ss_{\mu/\omega}}{-\mu} \ d\mu -
 \frac{i}{2\pi}\int_{\Cc} \frac{e^{-\rho}\log{\rm det}_F \Ss_{\rho/\e}}{-\rho} \
 d\rho \ ,
\end{equation}
using $\mu = \omega\la, \rho = \e\la$ and homotopy invariance to
shift the contours $\omega\Cc, \e\Cc$ to $\Cc$.

Since $\lim_{\omega\to \o}\log{\rm det}_F \Ss_{\mu/\omega}
=\log{\rm det}_F \Ss$ and since the contour in the first term in
\eqref{e:ecuttoff} can be closed at $\o$, we have
\begin{equation}\label{e:first-term}
\frac{i}{2\pi}\int_{\Cc} \lim_{\omega\to
\o}\frac{e^{-\mu}\log{\rm det}_F \Ss_{\mu/\omega}}{-\mu} \ d\mu =
\log{\rm det}_F \Ss \ .
\end{equation}
Now from \eqref{e:logdetSla}, \eqref{e:f0-asymp}, for any $\d>0$,
as $\e\to 0$
\begin{equation}\label{e:e-logdet}
\frac{\log{\rm det}_F \Ss_{\rho/\e}}{-\rho} = c_1(-\rho)\ii +
c_0(-\rho)\ii\log(-\rho) + c_0(-\rho)\ii\log(\e)
\end{equation}
$$+ \sum_{j=J+1}^{N-1}\sum_{k=0}^1
c_{j,k}((-\rho)^{1-\a_j}\e^{1-\a_j}\left(\log(-\rho) +
\log(\e)\right)^k + O(|\rho\e|^{-\a_N + \d}) \ . $$ Hence, as
$\e\to 0$,
\begin{equation}\label{e:second-term}
\frac{i}{2\pi}\int_{\Cc} \frac{e^{-\rho}\log{\rm det}_F
\Ss_{\rho/\e}}{-\rho} \ d\rho = c_1 . \frac{i}{2\pi}\int_{\Cc}
e^{-\rho}(-\rho)\ii \ d\rho \end{equation} $$ + c_0 .
\frac{i}{2\pi}\int_{\Cc} e^{-\rho}(-\rho)\ii\log(-\rho) \ d\rho  +
c_0\log(\e).\frac{i}{2\pi}\int_{\Cc} e^{-\rho}(-\rho)\ii \ d\rho
$$ $$\sum_{j=J+1}^{N-1}\sum_{k=0}^1 c_{j,k}\e^{1-\a_j}.
\frac{i}{2\pi}\int_{\Cc}
e^{-\rho}(-\rho)^{1-\a_j}\left(\log(-\rho) + \log(\e)\right)^k +
O(|\e|^{\a_N + \d}) \ . $$ From the contour integral formula for
the Gamma function
\begin{equation*}
\G(s)\ii  = \frac{i}{2\pi}\int_{\Cc} e^{-\rho}(-\rho)\si \ d\rho
\ ,
\end{equation*}
we have
\begin{equation}\label{e:gamma1}
 \frac{i}{2\pi}\int_{\Cc} e^{-\rho}(-\rho)\ii \ d\rho = \G(1)\ii
 = 1 \ ,
\end{equation}
and
\begin{equation}\label{e:gamma2}
\frac{i}{2\pi}\int_{\Cc} e^{-\rho}(-\rho)\ii\log(-\rho) \ d\rho \
 = -\frac{d}{ds}_{|s=1}\left(\G(s)\ii\right)
  = \G^{'}(1) \ .
\end{equation}
From \eqref{e:ecuttoff}, \eqref{e:first-term},
\eqref{e:second-term}, \eqref{e:gamma1}, \eqref{e:gamma2}, we
reach the conclusion.
\end{proof}

For closely related formulae see $\S$3 of \cite{Mu98}.

\subsection{Relative Eta Invariants}

The dependence of the relative $\z$-determinant on the choice of
spectral cut $R_{\th}$ is measured by the regularized limit in
\eqref{e:thm3}:
\begin{lem}\label{lem:theta}
Let $A_1,A_2$ be $\z$-comparable operators for spectral cuts
$\theta, \phi\in [0,2\pi)$ with scattering matrices
$\Ss_{\la}^{\th}, \Ss_{\mu}^{\phi}$ chosen so that $\det_F
\Ss_{0}^{\th} = \det_F \Ss_{0}^{\phi}$. Then
\begin{equation*}\label{e:theta}
\frac{{\rm det}_{\z,\th}(A_1,A_2)}{{\rm det}_{\z,\phi}(A_1,A_2)} =
 \exp\lsb -{\rm LIM_{\a\to\infty}}\log\lrb
 \frac{{\rm det}_F\Ss^{\th}_{e^{i\th}\a}}{{\rm det}_F\Ss^{\phi}_{e^{i\phi}\a}} \rrb\rsb \;
 ,
\end{equation*}
where $\LIM = \LIM^0, \;\a\in\R_{+}.$
\end{lem}

By \remref{rem:freedom} the requirement $\det_F \Ss_{0}^{\th} =
\det_F \Ss_{0}^{\phi}$ can always be fulfilled, and so
\eqref{e:theta} follows from \thmref{t:thmB} and
\eqref{e:linLIM}. Notice that the exponent is defined only
$\mod(2\pi i\Z)$ due to the ambiguity in defining $\log$, and that
the Fredholm determinants are taken on $H_{e^{i\th}\a},
H_{e^{i\phi}\a}$, respectively.

\vskip 2mm

We consider this for self-adjoint $A_1, A_2$ with ${\rm
sp}(A_i)\cap\R_{\pm}\neq\emptyset$; for example, for operators of
Dirac-type. There are, then, up to homotopy, two choices for
$\th$:
\begin{equation*}
\th = \frac{\pi}{2} \ , \hskip 10mm  -\frac{3\pi}{2}\leq \arg(\la)
< \frac{\pi}{2} \ , \hskip 10mm (-1)^s = e^{-i\pi s} \ ,
\end{equation*}
or,
\begin{equation*}
\th = \frac{3\pi}{2} \ , \hskip 10mm  -\frac{\pi}{2}\leq \arg(\la)
< \frac{3\pi}{2} \ , \hskip 10mm (-1)^s = e^{i\pi s} \ ,
\end{equation*}
which we may indicate by $\th = +, \th = -$, respectively. We
assume that $A_1, A_2$ are $\z$-comparable so that for
$\mu\in\Lambda_{\pm}$
\begin{eqnarray}
\hskip 20mm \Tr((A_1 - \mu)\ii &-& (A_2 - \mu)\ii) =
-\frac{\dd}{\dd\mu}\log{\rm det}_{F}\Ss^{\pm}_{\mu}
\label{e:thetareltrace} \\
     & \sim &  \sum_{j=0}^{\infty}\sum_{k=0}^{1}
  a_{j,k}^{\pm}(-\mu)^{-\a_j}\log^k (-\mu) \hskip 10mm {\rm as}
  \;\;\mu\to\o\;\;{\rm in}\; \Lambda_{\pm} \ . \label{e:thetaexp}
\end{eqnarray}
We may omit the $\pm$ superscripts in the following, $\mu$
indicating the appropriate scattering operator.

Since taking conjugates switches between spectral cuts it is
enough to assume $\z$-comparability for just one of $\th = \pi/2$
or $3\pi/2$. Considering $\mu =\pm i\a \in \pm i\R_{+}$ gives
\begin{equation}\label{e:aplusminus}
a_{J,0}^{-} = \ol{a_{J,0}^{+}} \hskip 10mm
\end{equation}
and
\begin{equation}\label{e:dlogplusminus}
\dd_{\a}\log{\rm det}_{F}\Ss_{-i\a} = \dd_{\a}\log{\rm
det}_{F}\Ss_{i\a}^{*} \ .
\end{equation}
Corresponding to $C_{+} = C_{\pi/2}$ and $C_{-} = C_{3\pi/2}$ we
have two relative $\z$-functions $\z_{\pm}(A_1,A_2,s)$. Their
disparity at $s=0$ is measured to first order by
$\z(A_{1}^{2},A_{2}^{2},0)$ and the {\it relative
$\eta$-invariant}
\begin{equation*}\label{e:releta}
\eta(A_1,A_2) := {\rm LIM}_{s\to 0} \, \eta(A_1,A_2,s)|^{mer} \ ,
\end{equation*}
where for $\re(s) >> 0$,
\begin{eqnarray}
\hskip 10mm \eta(A_1,A_2,s) & = & \Tr( A_1|A_1|^{-s-1} -
A_2|A_2|^{-s-1}) \nonumber \\
 & = & \frac{i}{2\pi} \int_{C_{\pi}}\la^{-\frac{s+1}{2}}
 \Tr( A_1(A_{1}^{2} - \la)\ii - A_2(A_{2}^{2} - \la)\ii ) \ d\la \
 ,
\label{e:reletafunction}
\end{eqnarray}
and $C_{\pi}$ is a contour of type \eqref{e:contour} with $\th =
\pi$. (Here $A_{i}^{2} = A_{i}^{*}A_i$ have (dense) domains ${\rm
dom}(A_{i}^{2}) = \{\xi\in H \ | \ \xi, A_i\xi\in {\rm
dom}(A_i)\}$).

The existence of $\z_{\pi}(A_{1}^{2},A_{2}^{2},s)$ and
$\eta(A_1,A_2,s)$ for $\re(s) >> 0$, the justification of
\eqref{e:reletafunction}, and their meromorphic continuation to
$\C$ follow from the $\z$-comparability of $A_1, A_2$ via the
identities
\begin{equation}\label{e:I1}
  A_i (A_{i}^2 - \la)\ii = \frac{1}{2}\left[(A_i - \la^{1/2})\ii +
  (A_i + \la^{1/2})\ii\right]  \ ,
\end{equation}
\begin{equation}\label{e:I2}
   (A_{i}^2 - \la)\ii = \frac{1}{2\la^{1/2}}\left[(A_i - \la^{1/2})\ii
   -  (A_i + \la^{1/2})\ii \right] \ .
\end{equation}
Here $\la^{1/2}$ is uniquely specified by $R_{\pi}$. It is
important to observe that the transformation $\la \to \la^{1/2}$
opens $C_{\pi}$ out into a vertical contour
\begin{equation}\label{e:Ctransform}
   C_{\pi} \mtoo  C_{1/2} := C_{\frac{\pi}{2},\downarrow} \cup
C_{-\frac{\pi}{2},\uparrow} \cup
C_{\rho,\frac{\pi}{2},-\frac{\pi}{2}} \ .
\end{equation}

From \eqref{e:I1}, $\|A_i (A_{i}^2 - \la)\ii\| =
O(|\la^{-1/2}|)$, so $A_i|A_i|^{-\frac{s+1}{2}}$ is defined for
$\re(s) > 0$. Since  $A_1, A_2$ are $\z$-comparable, \eqref{e:I1}
implies $A_1(A_{1}^{2} - \la)\ii - A_2(A_{2}^{2} - \la)\ii $ is
trace class, and, from \eqref{e:thetaexp},  $\eta(A_1,A_2,s)$ is
defined for $\re(s)> 1-\a_0$ (see \propref{p:Aeta}).

If $A_1, A_2$ are individually $\z$-admissible then $\eta(A_i,s) =
\Tr(A_i|A_i|^{-\frac{s+1}{2}})$ are defined, and since
$\eta(A_1,A_2,s)$ and $\eta(A_1,s) - \eta(A_2,s)$ have the same
pole structure
\begin{equation}\label{e:etasequal}
   \eta(A_1,A_2,s) =  \eta(A_1,s) - \eta(A_2,s) \hskip 20mm
   (A_i \;\; {\rm \z-admissible}) \ .
\end{equation}

Likewise, from \eqref{e:I2}, $\|(A_{i}^2 - \la)\ii\| =
O(|\la\ii|)$ and $\z(A_{1}^{2},A_{2}^{2},s)$ is defined for
$\re(s)> 1-\a_0$. More precisely, setting $\log = \log_{\pi}$
from here on, we have:
\begin{prop}\label{p:Asq} For
$\la\in\Lambda_{\pi}$
\begin{equation}\label{e:reltrsq}
   \Tr((A_{1}^{2} - \la)\ii - (A_{2}^{2} -
\la)\ii ) = -\frac{\dd}{\dd\la}\log{\rm det}_{F}\Ss_{\la^{1/2}}
-\frac{\dd}{\dd\la}\log{\rm det}_{F}\Ss_{-\la^{1/2}} \ .
\end{equation}
With $-\la=\a\in\R_{+}$
\begin{equation}\label{e:reltrsq2}
   \Tr((A_{1}^{2} + \a)\ii - (A_{2}^{2} +
\a)\ii ) = \frac{\dd}{\dd\a}\log{\rm det}_{F}((\Ss_{\pm
i\sqrt{\a}})^{*} \Ss_{\pm i\sqrt{\a}}) \ .
\end{equation}
[This means either $(\Ss_{-i\sqrt{\a}})^{*} \Ss_{-i\sqrt{\a}}$ or
$(\Ss_{+i\sqrt{\a}})^{*} \Ss_{+i\sqrt{\a}}$.] As $\la\to\o$ in
$\Lambda_{\pi}$ there is an asymptotic expansion
\begin{equation}\label{e:asympsq1}
   \Tr((A_{1}^{2} - \la)\ii - (A_{2}^{2} -
\la)\ii ) \sim \sum_{j=0}^{\infty}\sum_{k=0}^{1}
  a^{'}_{j,k}(-\la)^{-\frac{\a_j + 1}{2}}\log^k (-\la)  \ .
\end{equation}
In particular,
\begin{equation}\label{e:asympsq2}
 a^{'}_{J,0} = \frac{a^{+}_{J,0}+a^{-}_{J,0}}{2} \ ,
 \hskip 10mm a^{'}_{J,1} = 0 \hskip 20mm (\a_J = 1) \ .
\end{equation}
\end{prop}
\begin{proof}
The first equation is a consequence of \eqref{e:thetareltrace},
\eqref{e:I2} and
\begin{equation}\label{e:Phi}
 \Tr((A_{1}^{2} - \la)\ii - (A_{2}^{2} -
\la)\ii ) = \Psi(\la^{1/2}) + \Psi(-\la^{1/2}) \ ,
\end{equation}
where $\Psi(\rho) = (2\rho)\ii \Tr((A_{1} - \rho)\ii - (A_{2} -
\rho)\ii )$. Here, one uses \eqref{e:Ctransform} in order to
track which sector $\la^{1/2}$ is in, and hence whether
$\Ss_{\la^{1/2}}$ is $\Ss^{+}_{\mu}$ or $\Ss^{-}_{\mu}$. The
change of branch of $\log$ between \eqref{e:thetareltrace} and
\eqref{e:I2} is unimportant.

If $-\la = \a\in\R_{+}$, so $\la\in C_{\pi,\uparrow}$, then
$\la^{1/2} = -i\sqrt{\a}$ with respect to $R_{\pi}$. Since
$\ol{\Psi(\rho)} = \Psi(\ol{\rho})$, the right-side of
\eqref{e:Phi} becomes $\Psi(-i\sqrt{\a}) +
\ol{\Psi(-i\sqrt{\a})}$, or, equivalently, $\ol{\Psi(i\sqrt{\a})}
+ \Psi(i\sqrt{\a})$, resulting in the $\pm$ alternatives in
\eqref{e:reltrsq2}, which now follows similarly to
\eqref{e:reltrsq} using $\ol{\dd_{\a}\log{\rm
det}_{F}\Ss_{\pm\a^{1/2}}} = \dd_{\a}\log{\rm
det}_{F}(\Ss_{\pm\a^{1/2}})^{*}$.

\eqref{e:asympsq1} follows from \eqref{e:thetaexp} and
\eqref{e:Phi}. It is built from both of the expansions
\eqref{e:thetaexp} as $\la^{1/2}\to\o$ in $\Lambda_{\pm}$, but
the coefficients in \eqref{e:asympsq1} will not in general be of
the simple form \eqref{e:asympsq2}  due to the change in spectral
cut. However, $\la\ii$ is unambiguously defined and
\eqref{e:asympsq2} follows by comparing \eqref{e:thetaexp},
\eqref{e:reltrsq}, \eqref{e:asympsq1}.
\end{proof}

The content of \propref{p:Asq} is that $A_{1}^{2},A_{2}^{2}$ are
$\z$-comparable:
\begin{thm}\label{t:Asq}
Let $A_1, A_2$ be self-adjoint $\z$-comparable operators, as
above. Then $\z(A_{1}^{2},A_{2}^{2},s)$ is regular around $s=0$,
and (with $\th=\pi$, $\Ss = \Ss_{0}$)
\begin{eqnarray}
{\rm det}_{\z}(A^{2}_1,A^{2}_2) & = & {\rm det}_{F}(\Ss^* \Ss)
 e^{-{\rm LIM}_{\a\to\infty}\log
 {\rm det}_F((\Ss_{\pm i\a})^{*} \Ss_{\pm i\a})}
 \label{e:reldetsq} \\
& = & |{\rm det}_{F}\Ss |^2
 e^{-{\rm LIM}_{\a\to\infty}\log
 |{\rm det}_F(\Ss_{\pm i\a})|^2} \nonumber \ .
\end{eqnarray}
If $A_1, A_2$  are individually $\z$-admissible, then so are
$A^{2}_1, A^{2}_2$ and, then, ${\rm det}_{\z}(A^{2}_1,A^{2}_2) =
{\rm det}_{\z}(A^{2}_1)/{\rm det}_{\z}(A^{2}_2)$.
\end{thm}
\begin{proof}
The first statement is equation \eqref{e:asympsq2}. From
\propref{p:thmB*} with $\Phi(\la) = \log{\rm
det}_{F}\Ss_{\la^{1/2}} + \log{\rm det}_{F}\Ss_{-\la^{1/2}}$ we
obtain
\begin{equation}\label{e:anytheta}
{\rm det}_{\z}(A^{2}_1,A^{2}_2) = {\rm det}_{F}(\Ss^* \Ss)
 e^{-{\rm LIM}^{\pi}_{\la\to\infty}\Phi(\la)} \ .
\end{equation}
Notice, though $\log{\rm det}_{F}\Ss^{\pm}$ may differ by a
constant, \eqref{e:anytheta} is unambiguous. The regularized
limit averages the limits in the sectors $\Lambda_{\pm}$ and so
is well-defined in $\Lambda_{\pi}$. Equation \eqref{e:reldetsq}
now follows from \eqref{e:reltrsq2} and \eqref{e:anytheta} by
computing the limit along $R_{-} = -R_{+}$, and noting ${\rm
LIM}_{\a\to\infty}f(\sqrt{\a}) =  {\rm LIM}_{\a\to\infty}f(\a)$.
The final statement follows on applying $\dd_{\la}^m$ to
\eqref{e:I2} for large $m$.
\end{proof}

The  analogue of \propref{p:Asq} for the eta-function is proved
similarly:
\begin{prop}\label{p:Aeta} $A_{1},A_{2}$ are $\eta$-comparable, in so far as, for
$\la\in\Lambda_{\pi}$,
\begin{equation}\label{e:reltreta}
   \Tr(A_1(A_{1}^{2} - \la)\ii - A_2(A_{2}^{2} -
\la)\ii ) = -\la^{1/2}\frac{\dd}{\dd\la}\log{\rm
det}_{F}\Ss_{\la^{1/2}} + \la^{1/2}\frac{\dd}{\dd\la}\log{\rm
det}_{F}\Ss_{-\la^{1/2}} \ ,
\end{equation}
and as $\la\to\o$ in $\Lambda_{\pi}$ there is an asymptotic
expansion
\begin{equation}\label{e:etaasympsq}
   \Tr(A_1(A_{1}^{2} - \la)\ii - A_2(A_{2}^{2} -
\la)\ii ) \sim \sum_{j=0}^{\infty}\sum_{k=0}^{1}
  a^{''}_{j,k}(-\la)^{-\frac{\a_j}{2}}\log^k (-\la)  \ .
\end{equation}
\end{prop}

\vskip 1mm

From \eqref{e:etaasympsq} we obtain the singularity structure
\begin{equation}\label{e:etasing}
\G(\frac{s+1}{2})\eta(A_1,A_2,s) \sim \sum_{\stackrel{j=0}{j\neq J
}}^{\infty}\sum_{k=0}^{1} \frac{2^k A_{j,k}}{(s+\a_j -1)^{k+1}} +
\frac{2A_{J,0}}{s} + \frac{4A_{J,1}}{s^2} \ ,
\end{equation}
where the $A_{j,k}$ differ from the $a_{j,k}^{''}$ by universal
constants. Since $A_{1},A_{2}$ are $\z$-comparable,  from
\eqref{e:reltreta} we find that $a^{''}_{J,1} = 0$, and
$A_{J,1}=0$. Hence, as $\G(s)$ is regular at $s=1/2$,
$\eta(A_1,A_2,s)$ can have at most a simple pole at $s=0$.

\vskip 2mm

Though \propref{p:Asq} and \propref{p:Aeta} are ostensibly
similar, the $\eta$-invariant has a quite different character. In
particular, it is not necessary for $A_i$ to be invertible in
order to define $\eta(A_i), \eta(A_1,A_2)$ -- we require only
that  at $\la=0$ the resolvents $(A_i-\la)\ii$ are
meromorphic\footnote{This has consequences topologically:
concretely, the eta-invariant provides a canonical generator for
$\pi_{1}(\Ff_{sa})$, where $\Ff_{sa}$ is the space of
self-adjoint Fredholm operators, and a transgression form in the
relative family's index \cite{MePi97,Sc01}.}. It is then
convenient to consider
\begin{equation*}
  \wt{\eta}(A_1,A_2) = \frac{\eta(A_1,A_2) + \dim\Ker(A_1) -
  \dim\Ker(A_2)}{2}\ ,
\end{equation*}
since $\wt{\eta}(A_1,A_2)\mod(\Z)$ varies smoothly with
1-parameter families \cite{BrLe98,KiLe00,Mu94,LeWo96}.

The topological nature of $\eta(A)$ originates in the following
difference of regularized limits in the sectors $\Lambda_{\pm}$:
\begin{thm}\label{t:releta}
If $A_1, A_2$ are self-adjoint $\z$-comparable operators, as
above, then
\begin{equation}\label{e:eta}
\wt{\eta}(A_1,A_2) = \frac{1}{2\pi i}{\rm
LIM}_{\a\to\infty}\left(\log {\rm det}_F\Ss_{-i\a} - \log{\rm
det}_F\Ss_{i\a}\right) + \frac{1}{2}\z(A_{1}^{2},A_{2}^{2},0)
\hskip 2mm \mod(\Z).
\end{equation}
\end{thm}
\vskip 4mm
\begin{proof}
With $\mu = \la^{1/2}$ we have from \eqref{e:reltreta},
\eqref{e:Ctransform} that
\begin{equation*}
\eta(A_1,A_2,s) = \frac{i}{2\pi}\int_{C_{1/2}}\mu\si[\Tr((A_1 -
\mu)\ii - (A_2 - \mu)\ii)
\end{equation*}
$$\hskip 50mm + \ \Tr((A_1 + \mu)\ii - (A_2 + \mu)\ii)] \ d\mu$$
\begin{equation}\label{e:sqroot}
\hskip 33mm = \frac{i}{2\pi}\int_{\re(s)=c}\mu\si
\left(\dd_{\mu}\log{\rm det}_F \Ss_{-\mu} - \dd_{\mu}\log{\rm
det}_F \Ss_{\mu} \right) \ d\mu \ ,
\end{equation}
 where $c$ is positive and less than the smallest
positive spectral value of $A_1$ or $A_2$.

Since $\mu\si\log\det_F\Ss_{\mu}\to 0$ for $\re(s)>1-\a_0$ as
$\mu\to\o$, integrating \eqref{e:sqroot} gives
\begin{equation}\label{e:releta-G}
\eta(A_1,A_2) = \LIM_{s\to 0}(sG(s)) \ ,
\end{equation}
where $ G(s) = \frac{i}{2\pi}\int_{\re(\mu)=c}\mu^{-s-1} \,
F(\mu) \, d\mu,$ and
\begin{equation}\label{e:F}
F(\mu) = \log{\rm det}_F \Ss_{-\mu} - \log{\rm det}_F \Ss_{\mu}
\hskip 4mm \mod(2\pi i \Z) \ .
\end{equation}
Here, since
\begin{equation}\label{e:int-mu-s-1}
\int_{\re(\mu)=c}\mu^{-s-1} \, d\mu = 0 \hskip 15mm \re(s)>0 \  ,
\end{equation}
the $\mod(2\pi i \Z)$ ambiguity in \eqref{e:F} is not seen in
$G(s)$.

At $\mu= 0$: though $r(\mu)=\Tr((A_1 - \mu)\ii - (A_2 - \mu)\ii)$
is meromorphic with residue $-\dim\Ker(A_1) + \dim\Ker(A_2)$,
$r(\mu) + r(-\mu) = \dd_{\mu}F(\mu)$ is regular, and hence
$(-\mu)\ii F(\mu)$ is meromorphic with a Laurent expansion
$(-\mu)\ii F(\mu)=\sum_{j=-1}^{\o}b_j(-\mu)^j$. Let $F_0(\mu) =
(-\mu)\ii(F(\mu)- b_{-1})$.

From \eqref{e:int-mu-s-1} and since $F_0(\mu)$ is regular at
$\mu=0$
\begin{equation*}
G(s) = \frac{i}{2\pi}\int_{\re(\mu)=c}\mu^{-s} \, F_0(\mu) \, d\mu
= \frac{i}{2\pi}\int_{i\R}\mu^{-s} \, F_0(\mu) \, d\mu \ .
\end{equation*}
We have $\mu = re^{\pm i\pi/2}$ on $\pm i\R$ and since the
orientation goes from $+i\o$ to $-i\o$,
\begin{equation}\label{e:G2}
G(s) = -\frac{e^{-\frac{i\pi s}{2}}}{2\pi}\int_{\o}^0 r\si \,
F_0(re^{\frac{i\pi s}{2}}) \, dr  + \frac{e^{\frac{i\pi
s}{2}}}{2\pi}\int_{0}^{\o} r\si \, F_0(re^{\frac{-i\pi s}{2}}) \,
dr \ .
\end{equation}
By the argument in \thmref{t:thmB}, if $h$ is holomorphic in
$\Lambda_{\th}$ with an asymptotic expansion $h(\la) \sim
h_1(-\la)\ii + h_0 (-\la)\ii\log(-\la) +
\sum_{j=0}^{\o}\sum_{k=0}^1 h_{j,k}(-\la)^{-\a_j}\log^k (-\la)$,
with $0<\a_j\nearrow \o$, as $\la\to\o$, then $h(s) =
\int_{0}^{\o} r\si f(re^{i\th}) \, dr$ defined for $\re(s) >>0$
extends meromorphically to $\C$ with singularity structure around
$s=0$
\begin{equation}\label{e:singatzero}
  h(s) = e^{i(\pi-\th)}\left(\frac{h_0}{s^2} +
  \frac{h_1 + i(\th-\pi)h_0}{s}\right) + p(s) \ ,
\end{equation}
where $p(s)$ is meromorphic on $\C$ with poles at $s=1-\a_j\neq
0$.

Since $A_1, A_2$ are $\z$-comparable,  as in \eqref{e:f0-asymp} we
obtain asymptotic expansions
\begin{equation}\label{e:asympatzero}
F_0 (\mu) \sim a^{\pm}_1(-\mu)\ii + a^{\pm}_0 (-\mu)\ii\log(-\mu)
+ \sum_{j=0}^{\o}\sum_{k=0}^1 c^{\pm}_{j,k}(-\la)^{-\a_j}\log^k
(-\la) \ ,
\end{equation}
as $\mu\too\pm i\o$, where $0<\a_j\nearrow \o$. From \eqref{e:F}
we compute
\begin{equation}\label{e:a1}
a^{\pm}_1 = {\rm LIM}_{\mu\to \pm i\infty}\left(\log {\rm
det}_F\Ss_{\mu} - \log{\rm det}_F\Ss_{-\mu}\right) - b_{-1} \mp
i\pi a_{J,0}^{\pm} \hskip 3mm \mod(2\pi i \Z) ,
\end{equation}
and
\begin{equation}\label{e:a0}
a^{\pm}_0 = \pm(a_{J,0}^{\pm} + a_{J,0}^{\mp}) \ ,
\end{equation}
cf. \eqref{e:f0-asymp} ( $i\pi a_{J,0}^{\pm}$ in \eqref{e:a1}
arises from $\log(\mu) = \log(-\mu) + i\pi)$.

From  \eqref{e:G2}, \eqref{e:singatzero}, \eqref{e:asympatzero}
\begin{equation}\label{e:sG-singatzero}
  sG(s) = -\frac{e^{-\frac{i\pi s}{2}}}{2\pi i}\left(\frac{a^{+}_0}{s} +
  a^{+}_1 - \frac{\pi i}{2} a^{+}_0\right)
+ \frac{e^{\frac{i\pi s}{2}}}{2\pi i}\left(\frac{a^{-}_0}{s} +
  a^{-}_1 - \frac{3\pi i}{2} a^{-}_0\right) + O(s)
\end{equation}
$$ =  \frac{(a^{-}_0 - a^{+}_0)/2\pi i}{s} + \frac{(a^{+}_0 -
a^{-}_0)}{2} + \frac{(a^{-}_1 - a^{+}_1)}{2\pi i} + O(s) \ .$$
Hence from \eqref{e:releta-G}, \eqref{e:a1}, \eqref{e:a0}
\begin{eqnarray*}
\eta(A_1,A_2) & = & \frac{(a^{+}_0 - a^{-}_0)}{2} + \frac{(a^{-}_1
- a^{+}_1)}{2\pi i} \\
& = &\frac{1}{\pi i }{\rm LIM}_{\a\to\infty}\left(\log {\rm
det}_F\Ss_{-i\a} - \log{\rm det}_F\Ss_{i\a}\right) +
\frac{a_{J,0}^{+} + a_{J,0}^{-}}{2} \hskip 3mm \mod(2\Z) \\
& = & \frac{1}{\pi i }{\rm LIM}_{\a\to\infty}\left(\log {\rm
det}_F\Ss_{-i\a} - \log{\rm det}_F\Ss_{i\a}\right) + a^{'}_{J,0}
\hskip 3mm \mod(2\Z),
\end{eqnarray*}
the final equality using \eqref{e:asympsq2}. Since $\z(A_{1}^2,
A_{2}^2,0) = a^{'}_{J,0} - \dim\Ker(A_{1}^2) +
\dim\Ker(A_{1}^2)$, and, by self-adjointness,
$\dim\Ker(A_{i}^2)=\dim\Ker(A_i)$ we reach the conclusion.
\end{proof}

\begin{rem}
\noi (1) Similar regularized winding numbers to \eqref{e:eta} for
suspended $\pdo$s have been studied in \cite{Me95}.

\noi (2) The methods of $\S$ 2.1 can be used to obtain
\eqref{e:eta} through the heat formula $$\eta(A_1,A_2,s) =
\G(\frac{s+1}{2})\ii\int_{0}^{\o}t^{\frac{s-1}{2}}\Tr(A_1
e^{-tA_{1}^{2}} - A_2 e^{-tA_{2}^{2}}) \ dt \ .$$

\noi (3) From \eqref{e:dlogplusminus}, the regularized limit in
\eqref{e:eta} is pure imaginary, while from \eqref{e:aplusminus}
and \eqref{e:asympsq2} $\z(A_{1}^2, A_{2}^2,0)$ is real. When
$A_1, A_2$ are invertible this corresponds to the the role of
these invariants in defining the phase of the determinant
\begin{equation}\label{e:sadet}
  {\rm det}_{\z,\pm}(A_1,A_2) = e^{\mp i\frac{\pi}{2}(\eta(A_2,A_2) - \z(A_{1}^2,
A_{2}^2,0))}{\rm det}_{\z,\pi}(|A_1|,|A_2|) \ .
\end{equation}
\end{rem}

\vskip 2mm

\subsection{A Multiplicativity Property}
We refer to $\z$-comparable operators $A_1, A_2$ as {\it strongly
$\z$-comparable} if $a_{j,k}=0$ for $j\leq J$ in \eqref{e:exp3};
in particular, $\z_{rel}(0)=0$.

From \eqref{e:reglim}, if $A_1,A_2$ are strongly $\z$-comparable
then $\LIM = \lim$ and
\begin{equation}\label{e:thm3zeta0}
{\rm det}_{\z,\th}(A_1,A_2) = {\rm det}_F \Ss\, . \,
 e^{-\lim^{\th}_{\la\to\infty}\log{\rm det}_F\Ss_{\la}} \; .
\end{equation}
More precisely, $\z_{\th}(A_{1},A_{2},s)$ is then holomorphic for
$\re(s)> 1-\a_{J+1}$ and so at $s=0$, without continuation, and
\eqref{e:thm3zeta0} follows from
\begin{equation*}
\z_{\th}(A_{1},A_{2},s) = -\frac{\sin(\pi s)}{\pi}e^{i(\pi -
\th)s} \int_{0}^{\o} r\si \, \frac{\dd}{\dd r}\log{\rm det}_F
(\Ss_{re^{i\th}}) \ dr \ .
\end{equation*}
This applies to the following multiplicativity property:

\vskip 2mm

\begin{thm}\label{t:thmC}
Let $A:H\to H$ be closed and invertible with spectral cut
$R_{\th}$ with $\|(A-\la)\ii\| = O(|\la|\ii)$ as $\la\to\o$ in
$\Lambda_{\th}$ and let $Q = I + W: H\to H$ with $W$ of trace
class. If the operator $AW$ is trace class, then $(AQ,A)$ are
strongly $\z$-comparable and
\begin{equation}\label{e:mult1}
  {\rm det}_{\z,\th}(AQ,A) = {\rm det}_{F}Q \ .
\end{equation}
If $A$ is $\z$-admissible, then so is $AQ$ and
\begin{equation}\label{e:mult2}
{\rm det}_{\z,\th}(AQ) =  {\rm det}_{\z,\th}(A).{\rm det}_{F}Q \ .
\end{equation}
\end{thm}

\vskip 2mm

Let us point out some immediate consequences. First, we have:
\begin{prop}\label{p:multzeta}
Let $M$ be a closed $n$-manifold and $\textsf{D}: H^s (M,E)\to
H^{s-d} (M,E)$ an elliptic $\pdo$ of order $d>0$ acting on
sections of a vector bundle $E$ over $M$. Let $Q= I + W$ where
$W$ is a $\pdo$ on $L^2(M,E)$ of order $\ord(W) < -n-d$. Then
\begin{equation*}
{\rm det}_{\z,\th}(\textsf{D}Q) =  {\rm
det}_{\z,\th}(\textsf{D}).{\rm det}_{F}Q \ ,
\end{equation*}
In particular, this holds if $W$ is a smoothing operator.
\end{prop}
\begin{proof}
It is well known that in this case $\textsf{D}$ is
$\z$-admissible. On the other hand, $\ord(\textsf{D} W) < - n$
and hence $\textsf{D} W$ is trace class.
\end{proof}

This generalizes Lemma(2.1) of \cite{KoVi93}. On the other hand,
using the multiplicativity of the Fredholm determinant, setting
$A=Q_2$ and $Q = Q_2\ii Q_1$ we have:

\begin{prop}\label{p:zetaFred}
 \thmref{t:thmC}, \eqref{e:mult1}, applies to any bounded operator
$A$ on $H$. In particular, $(Q,I)$ are strongly $\z$-comparable,
with $I$ the identity operator, and
\begin{equation}\label{e:zetaFred1}
  {\rm det}_{\z,\th}(Q,I) = {\rm det}_{F}Q \ .
\end{equation}
Equivalently, if $Q_1, Q_2$ are determinant class, they are
strongly $\z$-comparable and
\begin{equation}\label{e:zetaFred2}
  {\rm det}_{\z,\th}(Q_1,Q_2) = {\rm det}_{F}(Q_2\ii Q_1)
   =  \frac{{\rm det}_{F}Q_1}{{\rm det}_{F}Q_2} \ .
\end{equation}
\end{prop}

\bigskip

Thus, although $\z_{\th}(Q_i,s)$ is undefined if $Q_i$ is of
determinant class for any $s$ (with $H$ infinite-dimensional),
$\z_{\th}(Q_{1},Q_{2},s)$ is defined and holomorphic for $\re(s)>
-1$ (see below). Since the contour can be closed at $\o$ this
extends to all $s$---equivalently, \eqref{e:zetaFred1} is
independent of $\th$, providing one perspective on the following.

First, note that for self-adjoint strongly $\z$-comparable
operators, $\eta(A_1,A_2,s)$ is holomorphic for $\re(s)>
1-\a_{J+1}$, and \eqref{e:eta} is immediate with $\LIM = \lim$ on
setting $s=0$ in \eqref{e:sqroot}.
\begin{prop}
Let $Q_1,Q_2$ be self-adjoint and determinant class. Then
$\eta(Q_1,Q_2,s)$ defined for $\re(s)>-2$ extends to an entire
function and $$\widetilde{\eta}(Q_1,Q_2)\in\Z \ .$$
\end{prop}
\begin{proof}
The first statement follows from \eqref{e:zet0=0}. The identity
is immediate from \eqref{e:eta} and \eqref{e:limlogdetclass}.
\end{proof}
An equivalent view point is that ${\rm det}_{\z,\th}(Q_1,Q_2)$ is
real by \eqref{e:zetaFred2}, and so the phase in \eqref{e:sadet}
must be real. Matters are quite different for differential
operators ($\S$ 3.4).

 The proof \thmref{t:thmC} is as follows.
\begin{proof}
We have $AQ-\la = A-\la + S$, where $S=AW$ is trace class. Hence
$(AQ-\la)\ii - (A-\la)\ii$ is trace class, and for $\la$ large
\begin{equation}\label{e:mult-resolvents}
(AQ-\la)\ii = (A-\la)\ii + \sum_{k\geq 1}(A-\la)\ii
(S(A-\la)\ii)^k \ .
\end{equation}
From the trace norm estimate
\begin{equation*}
\|(A-\la)\ii (S(A-\la)\ii)^k\|_{Tr} \leq \|(A-\la)\ii\|^{k+1}(
\|S\|_{Tr})^k < C|\la|^{-k-1}\ ,
\end{equation*}
as $\la\to\o$, we have $\|(AQ-\la)\ii - (A-\la)\ii\|_{Tr} =
O(|\la|^{-2})$ and hence that
\begin{equation}\label{e:zet0=0}
|\Tr\lrb(AQ-\la)\ii - (A-\la)\ii\rrb|  = O(|\la|^{-2}) \ .
\end{equation}
Therefore $\z_{\th}(AQ,A,s) = \Tr((AQ)\si - A\si)$ is holomorphic
for $\re(s)> -1$ and $\z_{\th}(AQ,A,0) = 0$.  Using the symmetry
of the trace we find
\begin{equation*}
\Tr\lrb (AQ-\la)\ii - (A-\la)\ii \rrb  =
-\frac{\dd}{\dd\la}\log{\rm det}_F((A-\la)\ii (AQ-\la)),
\end{equation*}
so the scattering matrix is $\Ss_{\la} = (A-\la)\ii (AQ-\la)$.
Hence from \eqref{e:thm3zeta0}
\begin{equation*}
  {\rm det}_{\z,\th}(AQ,A) = {\rm det}_{F}Q \, . \,
  e^{\lim_{\la\to\o}^{\th} \log{\rm det}_F \Ss_{\la}} \ .
\end{equation*}
Finally, it is easy to see that $\det_F\Ss_{\la} = 1 +
O(|\la|\ii)$ for large $\la$ and hence that
\begin{equation}\label{e:limlogdetclass}
{\rm lim}_{\la\to\o}^{\th} \log{\rm det}_F \Ss_{\la} = 0 \hskip
5mm {\rm mod}(2\pi i\mathbb{Z}),
\end{equation}
proving \eqref{e:mult1}.

If $A$ is $\z$-admissible, then applying $\dd_{\la}^{m}$ to
\eqref{e:mult-resolvents} for large enough $m$, we find
$\dd_{\la}^{m}(AQ-\la)\ii$ is trace class with an expansion
\eqref{e:exp1} with no term $(-\la)\ii\log(-\la)$. Hence
\eqref{e:mult2} follows from \eqref{e:reldet=ratio}.
\end{proof}

\begin{rem}
In fact, ${\rm det}_{\z,\th}(AQ,A) = {\rm det}_{F}Q .
e^{-\LIM_{\la\to\o}^{\th}\log{\rm
  det}_{F}\Ss_{\la}}$ for any determinant class $Q$---
though expected, the vanishing of the regularized limit is
unresolved.
\end{rem}


\section{An Application to Global Boundary Problems of Dirac-type}

We turn now to the application of \thmref{t:thmB} to elliptic
differential operators on manifolds with boundary.

\subsection{Analytic Preliminaries}

Let $X$ be a compact Riemannian manifold with (closed) boundary
manifold $\dd X = Y$. Let $E^1,E^2$ be Hermitian vector bundles
over $X$ and let $A: \Ci(X,E^1)\to \Ci(X,E^2)$ be a first-order
elliptic differential operator. We assume a collar neighborhood $U
= [0,1)\times Y$ of the boundary such that
\begin{equation}\label{e:collar}
  A_{|U} = \sigma\left(\frac{\dd}{\dd u} + \Aa + R\right),
\end{equation}
where $\Aa$ is a first-order self-adjoint elliptic operator on
$\Ci(Y,E_{|Y}^1)$, $R$ is an operator of order $0$, and $\sigma :
E_{|U}^1 \to E_{|U}^2$ a unitary bundle isomorphism constant in
$u$. When \eqref{e:collar} holds, then $A$ is of {\it Dirac-type}.
The case when $R=0$ is called the {\it product case}.

We define the space of interior solutions of $A$
\begin{equation*}
\Ker(A,s) := \{\psi\in H^s (X,E^1) \ | \ A\psi = 0 \ \ {\rm in} \
\ X\backslash Y \} \ ,
\end{equation*}
and its restriction to $Y$
\begin{equation*}
H(A,s) := \g_0\Ker(A,s) \ ,
\end{equation*}
where, for each real $s > 1/2$, $\g_0 :H^s (X,E^1)\to
H^{s-1/2}(Y,E_{|Y}^1)$ is the continuous operator restricting
sections of $E^1$ in the $s^{th}$ Sobolev completion to the
boundary. Because of the {\it Unique Continuation Property}
$\g_0:\Ker(A,s)\to H(A,s)$ is a bijection, while the {\it Poisson
operator}
\begin{equation}\label{e:poisson}
  \Kk_{A} := r\tilde{A}\ii\tilde{\g}^*\sigma : H^{s-1/2}(Y,E_{|Y}^1)
  \too \Ker(A,s) \subset H^s (X,E^1) \ .
\end{equation}
defines a canonical left inverse to $\g_0$. Here, using the
doubling construction \cite{BoWo93} for example, we consider $X$
as embedded in a closed manifold $\t{X}$ with Hermitian bundles
$\t{E}_i$ such that $(\t{E}_i)_{|X} = E_i$, and such that $A$
extends to an invertible elliptic operator $\tilde{A}:H^s
(\tilde{X},\tilde{E}^1)\to H^{s-1} (\tilde{X},\tilde{E}^2)$, and
where $\tilde{\gamma} :H^s (\tilde{X},\tilde{E}^1)\to
H^{s-1/2}(Y,E_{|Y}^1), \ r :H^s (\tilde{X},\tilde{E}^2)\to
H^{s}(X,E^2)$ are the continuous restriction operators.

\begin{prop}\label{p:calderon}\cite{Se66,Se69,BoWo93,Gr96,Gr99}
The restriction
\begin{equation}\label{e:Calderon}
P(A) := \tilde{\gamma}\Kk_A \,
\end{equation}
of $\Kk_A$ to $Y$ is a $\pdo$ projection of order $0$ (the {\it
Calderon projection}) on the space of boundary sections
$H^{s-1/2}(Y,E_{|Y}^1)$ with range $H(A,s)$.

For $(y,\xi)\in T^* Y\backslash \{0\}$, the principal symbol
$\s[P(A)](y,\xi):E_{y}^1 \to E_{y}^1$ is the orthogonal projection
with range $N_{+}(y,\xi)$ equal to the direct sum of eigenspaces
of the principal symbol of $\Aa$ with positive eigenvalue.
Therefore $\s(P(A))$ is independent of the operator $R$.

In general, $P(A)$ is only a projector (an indempotent), but
\begin{equation}\label{e:ort}
  P(A)_{{\rm ort}} := P(A)^* P(A)(P(A)P(A)^* +
  (I-P(A)^*)(I-P(A)))\ii ,
\end{equation}
is the $\pdo$ projection (unique self-adjoint indempotent) with
range
\begin{equation}\label{e:Qrange}
\ran(P(A)_{{\rm ort}}) = \ran(P(A)) = \ran(P(A)^*) \ .
\end{equation} and principal symbol
$\s[P(A)] =\s[P(A)_{{\rm ort}}]  \ . $
\end{prop}

\vskip 2mm

The Calderon projection provides a natural basepoint with which
to define global boundary problems:

\begin{defn}\label{defn:wellposed}\cite{Se69,Gr99}
A classical $\pdo$ $B$ of order $0$ acting in $H^s(Y,E_{|Y}^1)$
with principal symbol $\s[B]$ defines a boundary condition for $A$
which is {\it well-posed} if:

(i) $B$ has closed range for each real $s$;

(ii) for $(y,\xi)\in T^* Y\backslash\ \{0\}$, $\s[B](y,\xi)$ maps
$N_{+}(y,\xi)$ injectively onto the range of $ \s[B](y,\xi)$ in
$\C^N$.
\end{defn}

\begin{defn}\label{defn:admissible}
A well-posed boundary condition $B$ for $A$ is {\it admissible} if
the $\pdo$ $B-P(A)$ is a $\pdo$ of order $< -n$.
\end{defn}

\begin{rem}\label{rem:defns} In the following we shall for clarity
and brevity assume that if $B$ is admissible then $B-P(A)$ is a
smoothing operator. The modifications needed for the general case
are straightforward.
\end{rem}

For each well-posed boundary condition $B$ for $A$ the global
boundary problem
\begin{equation}\label{e:ebvp}
 A_B = A : \dom(A_B) \to L^2 (X,E^2)
\end{equation}
with domain $$\dom(A_B) = \{\psi\in H^1 (X,E^1) \ | \ B\g_0 \psi =
0\}$$ is a closed operator from $L^2 (X,E^1)$ to $L^2 (X,E^2)$. An
equivalent global boundary problem is obtained by replacing $B$
by the $\pdo$ projection $P[B] := P_{\Ker(B)\pp}$ so that
\begin{equation}\label{e:PB}
A_B = A_{P[B]} \ ,
\end{equation}
where for any closed subspace $W\subset H_Y$, $P_W$ denotes the
(orthogonal) projection with range $W$.

The following preferred sub-class of well-posed boundary
conditions is of special interest. By an {\it APS-type boundary
condition} we mean a $\pdo$ projection $P$ of order $0$ on $H_Y$
such that $P(A)-P$ is a $\pdo$ of order $-1$. The {\it
pseudodifferential Grassmannian} $Gr_{1}(A)$ is the
infinite-dimensional manifold parameterizing such projections,
each such $P\in Gr_{1}(A)$ defines a global boundary problem $A_P
: \dom(A_P) \too L^2 (X,E^2)$. In particular, $Gr_{1}(A)$
contains the APS projection $\Pi_{\geq}$. This property is quite
crude in so far as it follows trivially from the equality
$\s[P(A)]=\s[\Pi_{\geq}]$. If $R = 0$, the flow over the collar
leads to the following harder result:
\begin{prop}\label{p:differbysmoothing} \cite{Sc95,Gr99} In the product case
\begin{equation*}
  P(A) - \Pi_{\geq} \hskip 7mm {\rm and} \hskip 7mm  P(A)^* - \Pi_{\geq}
\end{equation*}
 are $\pdo$s of order $-\o$ (smoothing operators).
\end{prop}
Tailored to the product case we therefore also consider the dense
submanifold $Gr_{\o}(A)$ of $Gr_{1}(A)$ parameterizing those $P$
such that $P-P(A)$ is a smoothing operator.

Clearly, any $P\in Gr_{\o}(A)$ is admissible. The following facts
will be useful later:
\begin{lem}\label{lem:adsmoothing}
Let $B_1, B_2$ be admissible boundary conditions for $A$. Then
each of $B_1 - B_2$, $P[B_1] - P[B_2]$, $B_1 P[B_2]\pp$ are
smoothing operators, where $P\pp := I-P$, any projection $P$. In
particular, if $B$ is admissible, then $$ P[B]\in Gr_{\o}(A) \ .$$
\end{lem}
\begin{proof}
See \remref{rem:defns}. The first statement is obvious from
\defnref{defn:admissible}. The second follows easily from  $P[B_i] =
(i/2\pi)\int_{\Gamma}(B^{*}_i B_i - \la)\ii \, d\la$, with
$\Gamma$ a contour surrounding the origin and not enclosing any
eigenvalues of $B^{*}_i B_i$. For the third, one has $B_1
P[B_2]\pp = B_1(P[B_2]\pp - P_{\Ker(B_1)}) = B_1(P[B_2]\pp -
P[B_1]\pp).$
\end{proof}

The existence of the Poisson operator reduces the construction of
a parametrix for $A_B$ to the construction of a parametrix for the
operator on boundary sections
\begin{equation}\label{e:S(P)}
  S_A (B) = B\circ P(A) : H(A) \too W = \ran(B) \ .
\end{equation}
$S(B) =  S_A (B)$ is a Fredholm operator with kernel and cokernel
consisting of smooth sections. The corresponding properties for
$A_B$ follow from canonical isomorphisms
\begin{equation}\label{e:kercoker}
\Ker (S(B))\cong \Ker (A_B) \ , \hskip 15 mm \Cok (S(B))\cong \Cok
(A_B) \ .
\end{equation}
The first is defined by the Poisson operator. The second follows
in the same way from $\Cok (A_B) = \Ker (A^{*}_{B})$ and $\Cok (S
(B)) = \Ker (S_{A^*}(B\st))$. Here, the operator $$A^{*}_B :=
(A_B)^* = A^{*}_{B\st} : \dom(A^{*}_{B\st}) \to L^2 (X,E^1)$$ is
the adjoint realization of $A_B$ with $\dom(A^{*}_B) = \{\phi\in
H^1 (X,E^2) \ | \ B\st \g_0 \phi =0 \}$. $A^*$ is the formal
adjoint of $A$ and the adjoint boundary condition on $L^2
(Y,E^{2}_{|Y})$ is
\begin{equation}\label{e:adjP}
B\st = \s P[B]\pp \s\ii \ .
\end{equation}
This follows from Green's formula
\begin{equation}\label{e:Green's}
  <A\psi,\phi>_2 - <\psi,A^* \phi>_1 = -<\s\g_0\psi,\g_0\phi>_Y  \ ,
\end{equation}
which takes the distributional form on $H^s(X,E^i)$
\begin{equation}\label{e:AXA}
  A_{P(A)}\ii A = I - \Kk_A\g_0 \ .
\end{equation}
In the collar neighborhood $U$ of $Y$
\begin{equation}\label{e:collar*}
  A^{*}_{|U} = -\sigma\ii\left(\frac{\dd}{\dd u} + \s\Aa\s\ii + \s
  R\s\ii\right)\ .
\end{equation}
Hence $A^*$ is of Dirac-type with Poisson operator $\Kk_{A^*} :
H^{s-1/2}(Y,E_{|Y}^2)\to H^s (X,E^2)$ and Calderon projection
$P(A^*) = \g_0\Kk_{A^*}$ having range $H(A^*)=\g_0 \Ker(A^*)$.
There is an obvious diffeomorphism $Gr_{\o}(A^*)\cong Gr_{\o}(A)$,
$P\st \leftrightarrow P$, which, in view of
\begin{equation}\label{e:calderon*}
  P(A^*) = \s P(A)\pp \s\ii = P(A)^* \ , \hskip 15 mm H(A^*) = \s(H(A)\pp) \ ,
\end{equation}
is base point preserving. As in \eqref{e:S(P)}, $A^{*}_{B}$ is
modelled by the boundary operator
\begin{equation}\label{e:S(P^*)}
  S_{A^*}(B\st) = B\st\circ P(A^*)  = \s B\pp P(A)\pp \s\ii: \s H(A)\pp \too
  \s W\pp \ ,
\end{equation}
where $W\pp :=\ran(P[B]\pp).$

\begin{rem}
With the identifications \eqref{e:kercoker} at hand, elementary
arguments \cite{BoWo93} yield the identities \eqref{e:index1} and
\eqref{e:index2}. For an alternative proof using functorial
methods see \cite{Sc01}. Details on the above facts can be
accessed in \cite{Se66,Se69,Gr96,Gr99,BoWo93,ScWo99}.
\end{rem}

\vskip 2mm

\noi $\bb$ {\bf Construction of a relative inverse from $S
(B)\ii$}

\vskip 2mm

From \eqref{e:kercoker}, if $A_B$ is invertible then so is $S(B)$
and we can define the {\em Poisson operator of the global
boundary problem} $A_B$ by
\begin{equation*}
  \Kk_{A}(B) := \Kk_A S(B)\ii P_{W} : H^{s-1/2}(Y,E_{|Y}^1)
  \too H^s (X,E^1) \ ,
\end{equation*}
where $W = \ran(B)$. This restricts to an isomorphism
\begin{equation}\label{e:KAPisom}
\Kk_{A}(B)_{|W} : W \too \Ker(A)
\end{equation}
with inverse
\begin{equation}\label{e:KAPinverse}
\lrb\Kk_{A}(B)_{|W}\rrb\ii = (B\g_0)_{|\Ker(A)} \ .
\end{equation}

\begin{prop}\label{p:relinv}
Let $B, B_1, B_2$ be well-posed for $A$ such that the global
boundary problems $A_B, A_{B_1}, A_{B_1}$ are invertible. Then
one has
\begin{equation}\label{e:master}
  A_B\ii A = I - \Kk_A (B) B\g_0 \ ,
\end{equation}
and hence
\begin{equation}\label{e:relinv}
A_{B_1}\ii = A_{B_2} \ii   - \Kk_{A}(B_1)B_1\g  A_{B_2} \ii \ .
\end{equation}
\end{prop}
\begin{proof}
From \eqref{e:AXA}
\begin{equation}\label{e:first}
A_{P(A)}\ii = A_{P(A)}\ii A A_{B}\ii = (A_{P(A)}\ii A)A_{B}\ii =
(I - \Kk\g_0) A_{B}\ii = A_B\ii   - \Kk_{A}\g_0  A_B \ii \ .
\end{equation}
And so
\begin{equation*}
B\g_0 A_{P(A)}\ii = - B\g_0 \Kk_{A}\g_0 A_B \ii  = -S(B)P(A)\g_0
A_B\ii \ .
\end{equation*}
Applying $\Kk_A(B)$ to both sides, we have $\Kk_A (B)B\g_0
A_{P(A)}\ii = -\Kk_A \g_0 A_B\ii$. Substituting in
\eqref{e:first} yields
\begin{equation}\label{e:second}
 A_B\ii  = A_{P(A)}\ii - \Kk_A (B)B\g_0 A_{P(A)}\ii  \ .
\end{equation}
Hence, since $(\Kk_A (B)B\g_0) \Kk_A \g_0 = \Kk_A \g_0$,
\begin{equation*}
 A_B\ii A  = (I - \Kk_A \g_0) - \Kk_A (B)B\g_0 (I - \Kk_A \g_0) = I - \Kk_A (B)B\g_0 \ .
\end{equation*}
Hence
\begin{equation*}
A_{B_1}\ii = A_{B_1}\ii A_{B_2}A_{B_2} \ii = (A_{B_1}\ii A)A_{B_2}
\ii  = A_{B_2} \ii   - \Kk_{A}(B_1)B_1\g  A_{B_2} \ii \ .
\end{equation*}
\end{proof}

\begin{rem}
The relative-inverse formula appears in various forms in the
literature. We refer in particular to \cite{Fo87,Gr99,ScWo99}.
\end{rem}

\subsection{The Relative Abstract Determinant}
The scattering determinant for $\z$-comparable global boundary
problems arises canonically at the level of determinant lines.

\vskip 2mm

\noi $\bb$ {\bf Determinant lines}

\vskip 2mm

The determinant of a Fredholm operator $E:H^1\to H^2$ exists
abstractly not as a number but as an element $\det E$ of a complex
line $\Det(E)$. Elements of the {\em determinant line} $\Det(E)$
are equivalence classes $[\Ee,\la]$ of pairs $(\Ee,\la)$, where
$\Ee:H^1\to H^2$ such that $\Ee - E$ is trace class\footnote{A
bounded operator $T:H^1\to H^2$  is trace class if $\|T\|_{Tr} :=
\Tr(T^* T)^{1/2}<\o$, where here $Tr$ means the sum of the
eigenvalues, but $T$ does not have a trace unless $H^1 = H^2$.}
and relative to the equivalence relation $(\Ee q,\la) \sim
(\Ee,\det_F (q)\la)$ for $q:H^1 \to H^1$ of determinant class.
Complex multiplication on $\Det(E)$ is defined by
\begin{equation}\label{e:Cmult}
\mu . [\Ee,\la] = [\Ee,\mu\la] \ .
\end{equation}
 The {\em abstract determinant}
$\det E := [E,1]$ is non-zero if and only if $E$ is invertible,
and there is a canonical isomorphism
\begin{equation}\label{e:detline}
   \Det(E) \cong \wedge^{max}\Ker(E)^*
\otimes \wedge^{max}\Cok(E) \ .
\end{equation}
(Clearly, any two complex lines are isomorphic, the issue, here
and below, is whether there is a canonical choice of isomorphism.)

Taking quotients of abstract determinants in $\Det(E)$ coincides
with the (relative) Fredholm determinant:
\begin{lem}\label{lem:detfred}
Let $E_1 :H^1\to H^2, E_2 :H^1\to H^2$ be Fredholm operators such
that $E_i - E$ are trace class. Then provided $E_2$ is invertible
\begin{equation}\label{e:ratio}
  \frac{{\rm det}(E_1)}{{\rm det}(E_2)} =  {\rm det}_F (E_2\ii
  E_1) \ ,
\end{equation}
where the quotient on the left side is taken in $\Det(E)$.
\end{lem}
\begin{proof}
The left side of \eqref{e:ratio} is the ratio
\begin{equation*}
  \frac{[E_1 , 1]}{[E_2 , 1]}   =    \frac{[E_2 .(E_2\ii
  E_1) , 1]}{[E_2 , 1]}  =  \frac{[E_2 , {\rm det}_F (E_2\ii E_1)]}{[E_2 , 1]}
\end{equation*}
which from \eqref{e:Cmult} is equal to the asserted determinant.
\end{proof}
For example, in \propref{p:multzeta}, one has
$\det_{\z,\th}(\textsf{D}Q)/\det_{\z,\th}(\textsf{D}) =
\det(\textsf{D}Q)/\det(\textsf{D})$.

\vskip 2mm

\noi $\bb$ {\bf Relative determinant lines for global boundary
problems}

\vskip 2mm

For well-posed boundary conditions $B_1, B_2$ for $A$, the global
boundary problems $A_{B_1}, A_{B_2}$ have different domains and
hence the abstract determinants live in different complex lines
$$\det(A_{B_1})\in \Det(A_{B_1}) \ , \hskip 10mm \det(A_{B_2})\in
\Det(A_{B_2}) \ .$$ (We assume here that the $A_{B_i}$ are
invertible.) This means that the {\it relative abstract
determinant} $$\det(A_{B_1}, A_{B_2}) :=
\det(A_{B_1})/\det(A_{B_2})$$  is undefined as a complex number.
(Equivalently, although of the form identity plus a smoothing
operator, the operators $ D_{B_1}\ii D_{B_2}$ and $ D_{B_1}
D_{B_2}\ii$ do not have Fredholm determinants). Rather
$\det(A_{B_1}, A_{B_2})$ is a canonical element of the $${\it
relative \ determinant \ line} :=\Det(B_2,B_1)$$ of the boundary
Fredholm operator\footnote{This is the origin of gauge anomalies
on manifolds with boundary.} $$(B_2,B_1):= B_1\circ
P[B_2]:\ran(B_2) \too \ran(B_1) \ . $$ More precisely, there is a
canonical isomorphism
\begin{equation}\label{e:isom1}
\Det(A_{B_1})\cong\Det(A_{B_2})\otimes\Det(B_2,B_1)
\end{equation}
and a canonical isomorphism
\begin{equation}\label{e:isom2}
\Det(A_{B})\cong \Det(S (B)) \ , \hskip 10mm \det (A_{B})
\longleftrightarrow \det(S (B)) \ .
\end{equation}
 Formally, these follow from \eqref{e:kercoker} and
\eqref{e:detline} and are the determinant line analogues of
\eqref{e:index1} and \eqref{e:index2}, for precise constructions
see \cite{Sc99,Sc01}. The isomorphism \eqref{e:isom1} says that to
define $\det(A_{B_1}, A_{B_2})$ as a complex number requires a
non-zero element of $\Det(B_2,B_1)$. By an {\it auxiliary
operator} for $A_{B_1}, A_{B_2}$ we mean an invertible operator
$\Ee:\ran(B_2)\to \ran(B_1)$ such that $\Ee - (B_2,B_1)$ is
trace-class. Such an operator defines the non-zero element
$\det(\Ee)\in\Det(B_2,B_1)$ and we hence obtain a regularized
relative determinant, defined as the quotient taken in
$\Det(A_{B_1})$ via \eqref{e:isom1}
\begin{equation*}
{\rm det}_{\Ee}(A_{B_1},A_{B_2}) = \frac{{\rm det}A_{B_1}}{{\rm
det}A_{B_2}\otimes {\rm det}(\Ee)} \ .
\end{equation*}
In particular, if $(B_2,B_1)$ is invertible then we obtain the
{\it Relative Canonical Determinant} ${\rm
det}_{\Cc}(A_{B_1},A_{B_2}):={\rm
det}_{(B_2,B_1)}(A_{B_1},A_{B_2})$.

\begin{prop}\label{p:relcdet}
Let $\Ee = \Ee(B_1,B_2)$ be an auxiliary operator for the global
boundary problems $A_{B_1} ,A_{B_2}.$ If $B_1 , B_2$ are
admissible boundary conditions, then
\begin{equation}\label{e:relEdet}
{\rm det}_{\Ee}(A_{B_1},A_{B_2}) = {\rm det}_{F}\left(\frac{S
(B_1)}{\Ee S (B_2)}\right),
\end{equation}
where the Fredholm determinant is taken on $H(A)$, and the
operator quotient means $\left(\Ee S (B_{2})\right)^{-1}S
(B_1):H(A)\too H(A)$. For two choices of auxiliary operator $\Ee,
\Ee^{'}$, one has
\begin{equation}\label{e:Edependence}
{\rm det}_{\Ee}(A_{B_1},A_{B_2}) = {\rm det}_F (\Ee\ii\Ee^{'}) .
{\rm det}_{\Ee^{'}}(A_{B_1},A_{B_2}) \ .
\end{equation}
If $(B_2,B_1)$ is invertible, one has
\begin{equation}\label{e:relCdet}
{\rm det}_{\Cc}(A_{B_1},A_{B_2}) = {\rm det}_{F}\left(\frac{S
(B_1)}{B_1 S (B_2)}\right) \ .
\end{equation}
\end{prop}
\begin{proof}
The identities \eqref{e:Edependence} and \eqref{e:relCdet} are
obvious from  \eqref{e:relEdet}. From
\begin{equation*}
  S (B_1) = \Ee B_2 P(A) + (B_1 P(A) - \Ee B_2 P(A)) \ ,
\end{equation*}
and since $B_i - P(A)$ have smooth kernels and $\Ee - (B_2,B_1)$
is trace-class, it is readily verified that $\left(\Ee S
(B_{2})\right)^{-1}S (B_1)$ is of determinant class on $H(A)$.

From \eqref{e:isom1}, \eqref{e:isom2}  we obtain a commutative
diagram of canonical isomorphisms
 $$\begin{CD} \Det(A_{B_1})
@>{\simeq}>> \Det(A_{B_1})\otimes\Det (B_{2},B_{1})\\ @VV{\simeq}V
@VV{\simeq}V \\ \Det(S (B_1))  @>{\simeq}>> \Det(S
(B_2))\otimes\Det (B_{2},B_{1})
\end{CD}$$
\vskip 1mm \noi in which the vertical maps take the abstract
determinant elements to each other, while in the bottom map
$\det(\Ee S (B_{2}))\longleftrightarrow \det(S (B_{2}))\otimes
\det\Ee$. (See \cite{Sc99}.) By construction, we therefore have
\begin{equation*}
\frac{{\rm det}A_{B_1}}{{\rm det}A_{B_2}\otimes {\rm det}(\Ee)} =
\frac{\det(S (B_{1}))}{\det(\Ee S (B_{2}))} \ ,
\end{equation*}
and by \eqref{e:ratio} this is the right side of
\eqref{e:relEdet}.
\end{proof}

\subsection{Relative $\z$-Determinant of First-Order Global Boundary Problems}
To see that ${\rm det}_{\Ee}((A-\la)_{B_1},(A-\la)_{B_2})$ is a
scattering determinant for $A_{B_1},A_{B_2}$ and to compute the
regularized limit, we study the zeta determinant under variation
of the operator and the boundary conditions.

First, we study the operator variation, with fixed boundary
conditions:
\begin{prop}\label{p:variation}
Let $A_z : \Ci(X,E^1)\to\Ci(X,E^2)$ be a 1-parameter family of
Dirac-type operators depending smoothly on a complex parameter
$z$ such that in the collar $U$
\begin{equation}\label{e:collarz}
  A_{z|U} = \sigma\left(\frac{\dd}{\dd u} + \Aa_z + R_z\right) \ ,
\end{equation}
where $\sigma$ and the principal symbol $\s(\Aa_z)$ of $\Aa_z$ are
independent of $z$.

Let $\dot{A}_z = (d/dz)A_z$ and let $B_1,B_2$ be admissible for
$A_z$ such that $A_{z,B_1},A_{z,B_2}$ are invertible for each
$z$. Then $\dot{A}_z (A_{z,B_1}\ii - A_{z,B_2}\ii)$ is a trace
class operator on $L^{2}(X,E_2)$ with
\begin{equation}\label{e:variation1}
\Tr(\dot{A}_z (A_{z,B_1}\ii - A_{z,B_2}\ii))  =  \Tr\lrb
S_{z}(B_1)\ii d_z S_{z}(B_1) - S_{z}(B_2)\ii d_z S_{z}(B_2) \rrb
\ ,
\end{equation}
where $S_z := S_{A_z}$ and $d_z$ is defined by the covariant
derivative $P(A_z)\cdot\frac{d}{dz}\cdot P(A_z)$ on $\Hh = \cup_z
H(A_z)$ (see \remref{r:dz}). Relative to a choice of auxiliary
operator $\Ee:\ran(B_2)\to\ran(B_1)$ one has \vskip 1mm
\begin{equation}\label{e:variation2}
\Tr(\dot{A}_z (A_{z,B_1}\ii - A_{z,B_2}\ii))  =  \frac{d}{d
z}\log{\rm det}_{F}\left(\frac{S_z (B_1)}{\Ee S_z (B_2)}\right) \
.
\end{equation}
$\Tr$ and $\det_F$ on the right-side of \eqref{e:variation1},
\eqref{e:variation2} are taken on $H(A_z)$ \vskip 1mm
\end{prop}
\begin{proof}
From \eqref{e:relinv} and \eqref{e:KAPinverse} we
compute that on ${\rm dom}(A_{z,B_1})$
\begin{eqnarray}
B_1\g_0 A_{z,B_2}\ii A_{z,B_1} & = & B_1 \g_0 ( A_{z,B_1}\ii
-\Kk_{z}(B_2)\g_0 A_{z,B_1} \ii ) A_{z,B_1} \nonumber \\ & = &
-B_1 \g_0 \Kk_{z}(B_2)\g_0 \nonumber \\ & = &  -\Kk_{z}(B_1)\ii
\Kk_{z}(B_2)B_2\g_0 \nonumber
\\ & = & -S_{z}(B_1)S_{z}(B_2)\ii B_2\g_0 \label{e:fact1} \ .
\end{eqnarray}
The vector bundle structure on $ \cup_z H(A_z)\cong \cup_z
\Ker(A_z)$ follows from the smooth dependence of the operators on
$z$ \cite{Sc01}. Let $d_z$ be the induced operator covariant
derivative.  Since $\Kk_{z}(B_1)$ has range in $\Ker(A_z)$ then
$A_z \Kk_{z}(B_1) = 0$, and hence $\dot{A}_{z}.\Kk_{z}(B_1) = -
A_{z}\dot{\Kk}_{z}(B_1).$ Since $B_1 \g_0 \Kk_z (B_1) = B_1
\g_0$, then $B_1 \g_0 d_z (\Kk_z (B_1)) = d_z (B_1 \g_0) =0$ so
that
\begin{equation}\label{e:fact2}
\frac{d}{d z}(A_{z}).\Kk_{z}(B_1) = - A_{z,B_1}d_z \Kk_{z}(B_1) \
,
\end{equation}
We shall also need the identity
\begin{eqnarray}
B_2d_z\g_0 \Kk_{z}(B_1) & = &  \frac{d}{d z}\left(B_2\g_0
\Kk_{z}(B_1)\right) \nonumber \\ & = &  \frac{d}{d z}\left(B_2
P(A_z) S_{z}(B_1)\ii B_1\right) \nonumber \\ & = & \frac{d}{d
z}\left(S_z(B_2)S_{z}(B_1)\ii B_1\right) \ . \label{e:fact3}
\end{eqnarray}
From \lemref{lem:adsmoothing} and \eqref{e:relinv} we have that
\begin{equation}\label{e:zrelinv}
A_{z,B_1}\ii - A_{z,B_2} \ii =  -\Kk_{z}(B_1)B_1 \g_0 A_{z,B_2}
\ii = -\Kk_{z}(B_1)B_1 P[B_2]\pp \g_0 A_{z,B_2} \ii \ ,
\end{equation}
has a smooth kernel. Hence $\dot{A}_z (A_{z,B_1}\ii -
A_{z,B_2}\ii)$ is trace class and
\begin{eqnarray}
\Tr_{L^2}(\dot{A}_z (A_{z,B_1}\ii - A_{z,B_2}\ii))  & = &
-\Tr_{L^2}(\dot{A}_z \Kk_{z}(B_1)B_1 P[B_2]\pp  P[B_2]\pp\g_0
A_{z,B_2} \ii) \nonumber \\  & = & -\Tr_{W_2\pp}(P[B_2]\pp\g_0
A_{z,B_2} \ii\dot{A}_z \Kk_{z}(B_1)B_1 P[B_2]\pp ) \label{e:tr1}
\end{eqnarray}
using the fact that $P[B_2] \g_0 A_{z,B_2} \ii = 0$ for the first
equality and that $\dot{A}_z \Kk_{z}(B_1)B_1 P[B_2] \pp$ has a
smooth kernel and $P[B_2]\pp\g_0 A_{z,B_2} \ii$ is bounded for the
second. Since the operator $B_2\pp\g_0 A_{z,B_2}\ii\dot{A}_z
\Kk_{z}(B_1)B_1$ is a $\pdo$ of order $0$, and thus bounded, and
$B_1 P[B_2] \pp$ is smoothing ($\lemref{lem:adsmoothing}$), the
trace \eqref{e:tr1} is equal to
\begin{eqnarray}
& & \Tr_{W_1}(B_1 \g_0 A_{z,B_2}\ii \frac{d}{d
z}(A_{z})\Kk_{z}(B_1)B_1) \nonumber \\
& = & \Tr_{W_1}\lrb B_1 \g_0 A_{z,B_2}\ii
A_{z,B_1}d_z\Kk_{z}(B_1)B_1\rrb \hskip 42mm \text{by} \
\eqref{e:fact2} \nonumber \\ & = & -\Tr_{W_1}\lrb
S_{z}(B_1)S_{z}(B_2)\ii B_2\g_0 d_z\Kk_{z}(B_1)B_1 \rrb \hskip
29mm \text{by}\ \eqref{e:fact1} \nonumber  \\ & = &
-\Tr_{W_1}\lrb S_{z}(B_1)S_{z}(B_2)\ii \frac{d}{d z}\lrb
S_{z}(B_2)S_{z}(B_1)\ii B_1\rrb \rrb  \hskip 14mm \text{by}\
\eqref{e:fact3}  \nonumber
\\ & = &\Tr_{H(A_z)}\lrb S_{z}(B_1)\ii d_z S_{z}(B_1) -
S_{z}(B_2)\ii d_z S_{z}(B_2) \rrb \ \label{e:lhs} .
\end{eqnarray}
Here we use the fact that the expression inside the trace in the
final equality has a smooth kernel, and so is trace class, in
order to swap the order of the operators from the previous line.
For $R\neq 0$ in \eqref{e:collar} it is only the difference in
\eqref{e:lhs} that is smoothing. In the product case, though,
each term is individually trace class.

On the other hand, the right side of \eqref{e:variation2} is equal
to $$
  \frac{d}{d
z}\log{\rm det}_{F}\lrb(S_z (B_1)(\Ee S_z (B_2))\ii\rrb $$
\begin{eqnarray*} & = &
\Tr_{H(A_z)}\lrb (\Ee S_z (B_2)) S_{z}(B_1)\ii \frac{d}{d z}
(S_{z}(B_1)(\Ee S_z (B_2))\ii) \rrb \\
 & = &
\Tr_{W_1}\lrb \Ee S_z (B_2) S_{z}(B_1)\ii \frac{d}{dz}
(S_{z}(B_1)S_z (B_2)\ii)\Ee\ii \rrb \ ,
\end{eqnarray*}
and this is clearly equal to \eqref{e:lhs} by symmetry of the
trace.
\end{proof}

We now have:

\begin{thm}\label{t:thmD} Let $A$ be a first-order elliptic operator of
Dirac-type and let $E^1 = E^2$. Let $B_1, B_2$ be admissible
boundary conditions for $A$ and such that $A_{B_1}, A_{B_2}$ are
invertible with common spectral cut $R_{\th}$, and such that
\eqref{e:exp3} holds. Then $A_{B_1}, A_{B_2}$ are $\z$-comparable
and, relative to a choice of auxiliary operator $\Ee$, have
scattering operator determinant ${\rm det}_{\Ee}(A_{B_1}-\la,
A_{B_2}-\la)$. One has
\begin{equation}\label{e:firstreldet}
{\rm det}_{\z,\th}(A_{B_1}, A_{B_2})  = {\rm det}_{\Ee}(A_{B_1},
A_{B_2}) \, . \,
 e^{-{\rm LIM}^{\th}_{\la\to\infty}\log{\rm det}_{\Ee}(A_{B_1}-\la,
A_{B_2}-\la)} \ .
\end{equation}
If $A_{B_1}, A_{B_2}$ are $\z$-admissible, then ( in terms of
\eqref{e:relEdet})
\begin{equation}\label{e:Eversion}
 \frac{{\rm
det}_{\z,\th}(A_{B_1})}{{\rm det}_{\z,\th}(A_{B_2})} = {\rm
det}_{F}\left(\frac{S(B_1)}{\Ee S (B_2)}\right) \, . \,
 e^{-\LIM^{\th}_{\la\to\infty}\log\det_F ((\Ee S_{\la} (B_2))\ii S_{\la}
 (B_1))} \ .
\end{equation}
\end{thm}
\begin{proof}
From \propref{p:calderon} we have that $B_1, B_2 $ are well-posed
and  admissible for $A_{\la} := A - \la$, while
\eqref{e:variation2} becomes
\begin{equation}\label{e:relresolventtrace}
\Tr((A_{B_1}-\la)\ii - (A_{B_2}-\la)\ii)  =  -\frac{\dd}{\dd
\la}\log{\rm det}_{F}\left(\frac{S_{\la} (B_1)}{\Ee S_{\la}
(B_2)}\right) \ ,
\end{equation}
Hence, with the stated assumptions, $A_{B_1}, A_{B_2}$ are
$\z$-comparable and \eqref{e:firstreldet} is immediate from
\thmref{t:thmB} and \eqref{e:relEdet}. Finally,
\eqref{e:Eversion} follows from \lemref{lem:relz=zrel} and
\eqref{e:reldet=ratio}.
\end{proof}

That the right-sides of \eqref{e:firstreldet}, \eqref{e:Eversion}
are independent of the choice of $\Ee$ is clear from
\eqref{e:Edependence} and \eqref{e:linLIM}. More precise
knowledge of the dependence of the regularized limit on the
operators $A$ and $B_i$ is obtained as follows.

\begin{prop}\label{cor:relvar}
With the conditions of \propref{p:variation},
\begin{equation}\label{e:relvar}
  \frac{d}{dz}\log{\rm det}_{\z,\th}(A_{z,B_1}, A_{z,B_2}) =
   \frac{d}{dz}\log{\rm det}_{\Ee}(A_{z,B_1}, A_{z,B_2}) \ ,
\end{equation}
and hence is independent of $\th$. Moreover, with
$\zeta_{z,rel}(0) := \z_{\th}(A_{z,B_1},A_{z,B_2},0)$
\begin{equation}\label{e:relzetazerovar}
  \frac{d}{dz}\zeta_{z,rel}(0) = 0 \ .
\end{equation}
The regularized limit term in \thmref{t:thmD} is independent of
the operator $A$, and depends only on the pseudodifferential
boundary conditions $B_1,B_2, P(A)$.
\end{prop}
\begin{proof} $(A_{z,B_1}-\la)\ii
- (A_{z,B_2}-\la)\ii $ has a smooth kernel, and hence the
operator $J(\la)= {\dot {A_z}}((A_{z,B_1}-\la)\ii -
(A_{z,B_2}-\la)\ii)$ is trace class. A well known argument
\cite{Fo87} gives $(d/dz)(\Tr((A_{z,B_1}-\la)\ii -
(A_{z,B_2}-\la)\ii)) = -\dd_{\la}\Phi(\la) \ ,$ where $\Phi(\la)
= \Tr(J(\la))$. Hence
\begin{equation}\label{e:dzeta}
-\frac{d}{dz} \z_{\th}(A_{z,B_1},A_{z,B_2},s) = \frac{i}{2\pi
}\int_{C}\la^{-s}\dd_{\la}\Phi(\la) \ d\la \ ,
\end{equation}
and so from \propref{p:thmB*}, $-\frac{d}{dz}
\z_{\th}^{'}(A_{z,B_1},A_{z,B_2},0)  =   \Tr(J(0)) - {\rm
LIM}^{(\th)}_{\la\to\o}\Tr(J(\la).$ By \eqref{e:variation2},
then, \eqref{e:relvar} is equivalent to ${\rm
LIM}^{(\th)}_{\la\to\o}\Tr(J(\la))=0.$ To see that, for $\re(s) >
1-\a_0$ we can integrate by parts in \eqref{e:dzeta} to obtain
\begin{equation}\label{e:dzeta2}
-\frac{d}{dz}\z_{\th}(A_{z,B_1},A_{z,B_2},s) = s\Tr({\dot {A_z}}
(A_{z,B_1}^{-s-1}- A_{z,B_2}^{-s-1})) \ .
\end{equation}
But $A_{z,B_1}^{-s-1}- A_{z,B_2}^{-s-1}$ has a smooth kernel for
$\re(s)> -1$ and hence \eqref{e:dzeta2} holds in that larger
half-plane. Setting $s=0$ in \eqref{e:dzeta2} therefore proves
\eqref{e:relzetazerovar}, while differentiating and setting $s=0$
we obtain
\begin{equation}\label{e:dzeta3}
-\frac{d}{dz}\z^{'}_{\th}(A_{z,B_1},A_{z,B_2},0) = \Tr({\dot
{A_z}} (A_{z,B_1}^{-1}- A_{z,B_2}^{-1})) \ ,
\end{equation}
which is equation \eqref{e:relvar}.

For the final statement, let $A_r , -\ep <  r < \ep $ be a smooth
path of Dirac-type operators, as in \propref{p:variation}, with
$A_{0} = A$. The variation near $\dd X$ is at most order $0$ and
so from \propref{p:calderon} $B_1, B_2$ are admissible for each
$A_r$. For small enough $\ep$ we can apply \eqref{e:relvar} and
comparing with \eqref{e:firstreldet} we reach the conclusion.
\end{proof}
If $A_r , \ 0 \leq r \leq t,$ is a smooth 1-parameter family
satisfying \eqref{e:collarz} with $A_{r,B_i}$ invertible, then
the final statement of \propref{cor:relvar} can equivalently be
expressed by the integrated version of \eqref{e:relvar}:
\begin{equation*}
\frac{{\rm det}_{\z,\th}(A_{t,B_1}, A_{t,B_2})} {{\rm
det}_{\z,\th}(A_{0,B_1}, A_{0,B_2})} = {\rm
det}_{F}\left(\frac{S_t(B_1)}{\Ee S_t(B_2)}\right) / {\rm
det}_{F}\left(\frac{S_0(B_1)}{\Ee S_0(B_2)}\right) \ .
\end{equation*}

Next, we compute the variation of the $\z$-determinant with
respect to the boundary condition. The following formula gives a
general direct variational formula\footnote{The usual approach to
computing the boundary variation is to try to `gauge transform'
the variation into an equivalent operator variation, see
\cite{DoWo91,ScWo99,Mu94} and also $\S 4.1$.}.
\begin{prop}\label{p:bvar}
Let $\{B_r \ | \ -\ep < r < \ep\}$ be a smooth 1-parameter family
of $\pdo$s on $L^2 (Y, E_{|Y}^1)$ such that $B_r - P(A)$ has a
smooth kernel and such that $A_{B_r}$ is invertible for each $r$
and $\z$-admissible. Then setting $S_{\la}(B_r) =
S_{(A-\la)}(B_r)$, one has
\begin{equation}\label{e:bvar}
\frac{d}{dr}\log{\rm det}_{\z,\th}(A_{B_r}) = \Tr\lrb S(B_r)\ii
\frac{d}{dr}S(B_r)\rrb  - \LIM^{\th}_{\la\to\o}\Tr\lrb
S_{\la}(B_r)\ii \frac{d}{dr}S_{\la}(B_r)\rrb \ .
\end{equation}
\end{prop}
\vskip 2mm
\begin{proof} Each $B_r$ is an admissible
boundary condition for the Dirac-type operator $A_{\la} := A -
\la$ and $\dot{B_r} = d/dr(B_r)$ is a smoothing operator on $L^2
(Y, E_{|Y}^1)$. We have
\begin{equation}
\frac{d}{dr}\frac{\dd}{\dd \la}\log{\rm
det}_F\left(\frac{S_{\la}(B_r)}{\Ee S_{\la}(B_0)}\right) =
\frac{\dd}{\dd \la}\Tr_{H(A_{\la})}(S_{\la}(B_r)\ii
\frac{d}{dr}S_{\la}(B_r)) \ .
\end{equation}
Hence from \eqref{e:relresolventtrace}
$$-\frac{d}{dr}\z_{\th}(A_{B_r},s) =
-\frac{d}{dr}\z_{\th}(A_{B_r},A_{B_0},s)
  =
  -\frac{i}{2\pi }\int_{C}\la^{-s} \frac{\dd}{\dd\la} \Tr(S_{\la}(B_{r})\ii
\frac{d}{dr}S_{\la}(B_{r})) \
  d\la \ ,$$
and so the result follows from \eqref{e:Zderiv}.
\end{proof}

\subsection{Local Coordinates and an Odd-Dimensional Example}

The identity \eqref{e:Eversion} can be given a more familiar form
if we work in local coordinates on $Gr_{\o}(D)$.

 \vskip 2mm

\noi $\bb$ {\bf The relative zeta determinant in Stiefel
Coordinates }

\vskip 2mm

To be concrete,  let $X$ be a compact Riemannian spin manifold and
consider a compatible Dirac operator $A=D:\Ci(X,E^1) \to
\Ci(X,E^2)$ acting between Clifford bundles in the product case
($R=0$). First, observe that to each basepoint $\Pi\in
Gr^{(0)}_{\o}(D)$ there is a dense open subset of
$Gr^{(0)}_{\o}(D)$. Setting
\begin{equation}
E = \ran(\Pi) \ , \hskip 10mm W = \ran(P) \ , \hskip 10mm W_i =
\ran(P_i) \ ,
\end{equation}
it is defined by
\begin{equation}\label{e:Uw1}
 U_E  =  \{P\in Gr^{(0)}_{\o}(D) \ | \ (\Pi, P) := P\circ \Pi :
  E \to W \; {\rm invertible} \} \ .
\end{equation}
Equivalently,
\begin{equation}\label{e:Uw2}
P\in U_E \Longleftrightarrow \ran(P)   = \gr(T:E\to E\pp) \hskip
15mm T\in\Hom_{\o}(E,E^{\perp}) \ ,
\end{equation}
where $\Hom_{k}(E,E^{\perp})= \{\Pi^{\perp}\textsf{Z}\Pi : E \to
E^{\perp} \ | \ \textsf{Z}\in \Psi_{k}(H_Y)\}$ and
$\Psi_{k}(H_Y)$ is the space of $\pdo$s on $H_Y = L^2
(Y,E_{|Y}^1)$ of order $-k\in\R_+\cup\{\o\}$. Equivalently,
\begin{equation}\label{e:Uw3}
P\in U_E \Longleftrightarrow P = P_T  = \begin{pmatrix}
    Q_{T}^{-1} & Q_{T}^{-1}T^{*} \\
    TQ_{T}^{-1} & TQ_{T}^{-1}T^{*}
  \end{pmatrix}, \hskip 10mm Q_T := I + T^* T \ ,
\end{equation}
$T\in\Hom_{\o}(E,E^{\perp})$. In this way a atlas for
$Gr^{(0)}_{\o}(D)$ can be constructed with respect to a countable
set of basepoint spectral projections in a similar way to
\cite{PrSe86}.

It is always possible to arrange for $P(D), P_1, P_2$ to lie in a
single coordinate patch $U_{E} \subset Gr^{(0)}_{\o}(D)$, so that,
\begin{equation}\label{e:graph}
  H(D) = \gr(K:E\to E\pp) \ , \hskip 8mm W_i = \gr(T_i:E\to
  E\pp) \ ,
\end{equation}
where $K, T_i\in \Hom_{\o}(E,E^{\perp})$. Since the operators
$S(P_i) = P_i P(D)$ are invertible, one such choice is $E=H(D)$,
in which case $K=0$. We may further assume, by perturbing $\th$
slightly if necessary, that the global boundary problem $D_{\Pi}$
has no eigenvalue along $R_{\th}$, and hence that $\Pi \circ
P(D-\la) : H(D-\la) \to \ran(\Pi)$ is invertible. This means that
\begin{equation}\label{e:graph2}
H(D-\la) = \gr(K_{\la} :E\to E\pp) \ , \hskip 8mm P(D-\la) =
P_{K_{\la}} \ ,
\end{equation}
for some unique $K_{\la}\in\Hom_{1}(E,E^{\perp})$ (the space of
restrictions of $\pdo$s of order $-1$). Recall that $P(D-\la)$ is
an element of $Gr_{1}(D)$, though not of $Gr_{\o}(D)$ if $\la\neq
0$.

We then have:

\begin{prop}\label{p:stiefel}
\begin{equation}\label{e:graphscattering}
{\rm det}_{F}\lrb\frac{S_{\la}(P_1)}{\Ee S_{\la}(P_2)}\rrb
 = {\rm det}_F \lrb\frac{I +
T^{*}_1 K_{\la}}{I + T^{*}_2 K_{\la}}\rrb \, . \, {\rm det}_F
(\Phi(\Ee,T_1,T_2)) \ ,
\end{equation}
where the Fredholm determinants are taken on $E$. The operator
$\Phi(\Ee,T_1,T_2): E\to E$ is invertible and independent of
$K_{\la}$. If $(P_2,P_1)$ is invertible, then
\begin{equation}\label{e:graphscattering2}
{\rm det}_{F}\lrb\frac{S_{\la}(P_1)}{P_1 S_{\la}(P_2)}\rrb
 = {\rm det}_F \lrb\frac{I +
T^{*}_1 K_{\la}}{I + T^{*}_2 K_{\la}}\rrb  \, . \, {\rm det}_F
(Q_2(I + T^{*}_1 T_2)\ii) \ .
\end{equation}
\end{prop}
\begin{proof}
Let $P,\t{P}\in Gr_{1}(D)$ with $\ran(P) = \gr(T),\, \ran(\t{P}) =
\gr(\t{T})$, where $T,\t{T} \in \Hom_{1}(E,E^{\perp})$. Then any
linear operator $R: \ran(P)\to \ran(\t{P})$ acts by $$(\xi, T\xi)
\mtoo (\Phi(R)\xi,\t{T}\Phi(R)\xi) \ ,$$ for some $\Phi(R)\in
\End(E)$. $\Phi$ respects operator composition: if $\t{R}:
\ran(\t{P})\to \ran(P)$
\begin{equation}\label{e:phihom}
\Phi(\t{R}R)=\Phi(\t{R})\Phi(R) \ .
\end{equation}
Moreover, if $\t{R}R:
\ran(P)\to \ran(P)$ is determinant class, then so is
$\Phi(\t{R}R)$ and
\begin{equation}\label{e:phidet}
  {\rm det}_F(\t{R}R) = {\rm det}_F(\Phi(\t{R}R))  \ ,
\end{equation}
where the left-side is taken on $\ran(P)$ and the right-side on
$E$.

In particular, the auxiliary operator $\Ee$ acts via
$\Phi(\Ee)\in\End(E)$. Similarly, $(P_2, P_1)$ acts via
$\Phi(P_2, P_1)\in\End(E)$ and  with $Q_i := Q_{T_i}$ one has
using \eqref{e:Uw3}
\begin{equation*}
(P_2, P_1) := P_1 \circ P_2 \begin{pmatrix}
  \xi \\
  T_2 \xi
\end{pmatrix} = \begin{pmatrix}
   Q_1\ii (I+T^{*}_1 T_2)\xi \\
  T_1 Q_1\ii (I+T^{*}_1 T_2)  \xi \end{pmatrix}
\end{equation*}
and so
\begin{equation}\label{e:case12}
  \Phi(P_2,P_1) = Q_1\ii (I+T^{*}_1 T_2) \ .
\end{equation}
Hence $(P_2, P_1)$ is invertible when $-1\notin {\rm sp}(T^{*}_1
T_2)$.  On the other hand, since $\Ee$ is invertible so is
$\Phi(\Ee)$, and because $\Ee - (P_2,P_1)$ is trace-class then
$\Phi(\Ee)- Q_i\ii (I + T^{*}_1 T_2)$ is also trace-class. Hence
$\Phi(\Ee)$ is of determinant class and $\det_F (\Phi(\Ee)) \neq
0$. It is easy to compute that $$\Phi(S_{\la}(P_1)) = Q_1\ii (I +
T^{*}_1 K_\la) \ ,\hskip 10mm \Phi(\Ee S_{\la}(P_2)) = \Phi(\Ee)
Q_2\ii (I + T^{*}_2 K_\la) \ .$$ From \eqref{e:phihom},
\eqref{e:phidet}, and the symmetry and multiplicativity of
$\det_F$, then by setting $\Phi(\Ee,T_1,T_2) = Q_2\Phi(\Ee)\ii
Q_1\ii$ we reach the conclusion.

Alternatively, since $(\Ee S_{\la} (P_2))\ii S_{\la}(P_1)) =
S_{\la} (P_2)\ii\Ee\ii S_{\la} (P_1)$, the computation can be
carried through by observing that for any invertible operator $R$
as above, one has (relative to graph coordinates)
\begin{equation}\label{e:Rinverse}
  P R\ii\t{P} =
\begin{pmatrix}
  \Phi(R)\ii & \Phi(R)\ii \t{T}^* \\
  T\Phi(R)\ii & T\Phi(R)\ii \t{T}^*
\end{pmatrix} \ .
\end{equation}
\end{proof}

We can now restate \thmref{t:thmD} as follows:

\vskip 2mm

\begin{thm}\label{t:thmD*} Let $D$ be a first-order Dirac-type operator in
the product case and let $E^1 = E^2$. Let $P_1, P_2\in
Gr_{\o}(D)$ such that $D_{P_1}, D_{P_2}$ are invertible with
common spectral cut $R_{\th}$, and such that \eqref{e:exp3}
holds. Then $D_{P_1}, D_{P_2}$ are $\z$-comparable and in local
Stiefel (graph) coordinates, as above, one has
\begin{equation}\label{e:stiefelversion}
 \frac{{\rm det}_{\z,\th}(D_{P_1})}{{\rm det}_{\z,\th}(D_{P_2})} =  \frac{{\rm
det}_F (I + T^{*}_1 K)}{{\rm det}_F(I + T^{*}_2 K)} \, .  \,
\exp\lsb -\LIM^{\th}_{\la\to\infty}\log {\rm det}_F \lrb\frac{I +
T^{*}_1 K_{\la}}{I + T^{*}_2 K_{\la}}\rrb \rsb \ .
\end{equation}
\end{thm}
\begin{proof} This is immediate from \eqref{e:Eversion}, \eqref{e:linLIM},
\eqref{e:graphscattering}. Because $P(D)\in Gr_{\o}(D)$ we have
replaced the determinant of the quotient by the quotient of the
determinants in the first term on the right-side of
\eqref{e:stiefelversion}.
\end{proof}

\begin{rem}\label{rem:stiefel}
More generally, Stiefel coordinates on $Gr^{(0)}_{\o}(D)$ refer to
an operator $[ M \; N] \in \Hom(E\oplus E\pp,E)$, where $M$ is
Fredholm with $\ind(M) = 0$, and $N\in \Hom_{\o}(E\pp,E)$. This
defines a point in the principal Stiefel frame bundle $ST_E \too
Gr^{(0)}_{\o}(D)$ (based at $E$), with bundle projection map
\begin{equation}\label{e:Pstiefel}
[ M \; N] \mtoo  P :=
\begin{pmatrix}
    M^{*}\Mm\ii M & M^{*}\Mm\ii N \\
    N^{*}\Mm\ii M & N^{*}\Mm\ii N
  \end{pmatrix} \ ,
\end{equation}
where $\Mm = M M^{*} + NN^{*}.$ In particular, graph coordinates
correspond to the canonical section $P_T \mto[ I \; T^*] $ of
$ST_E$ over  $U_E$. Stiefel coordinates $[ M_i \; N_i] $ for
$P_i$ modify \eqref{e:stiefelversion} by replacing $I+T^{*}_i
K_{\la}$ by  $M_i + N_i K_{\la} : E \to E$.
\end{rem}
\vskip 4mm

\noi $\bb$ {\bf Example: odd-dimensions revisited}

\vskip 2mm

To illustrate these formulae, we explain how they work for a
Dirac operator with $X$ odd-dimensional. In this case
\eqref{e:collar} takes the form
\begin{equation}\label{e:odd}
\begin{pmatrix}
i & 0 \\ 0 & -i
\end{pmatrix}
\left(
\partial_u +
\begin{pmatrix}
0 & \textsf{D}^{-}_Y \\ \textsf{D}^{+}_Y & 0
\end{pmatrix}
\right) \,
\end{equation}
\noi with respect to the decomposition $ H_Y = F^+ \oplus F^-$
into chiral spinor fields, where $\textsf{D}^{+}_Y$ is the chiral
Dirac operator, which, for brevity, we shall assume invertible.
The projection $P_{F^+}$ onto $F^+$ is not an element of
$Gr_{1}(D)$. It does, however, define a true (local) elliptic
boundary condition and is related to $\Pi_{\geq}$ in the
following precise way.

The involution defining the grading of $H_Y$ into positive and
negative energy (the APS condition) is the operator
\begin{equation*}
\Dd_Y\ii|\Dd_Y| =
  \begin{pmatrix}
     0 & (\textsf{D}^{+}_Y)\ii(\textsf{D}^{+}_Y \textsf{D}^{-}_Y)^{1/2} \\
     (\textsf{D}^{+}_Y \textsf{D}^{-}_Y)^{-1/2}\textsf{D}^{+}_Y & 0
   \end{pmatrix} \ .
\end{equation*}
Hence, defining $g_{+}$ to be the unitary isomorphism
$(\textsf{D}^{+}_Y \textsf{D}^{-}_Y)^{-1/2}\textsf{D}^{+}_Y  :
F^+ \to F^-$, we have
\begin{equation*}
\Pi_{\ge} = \frac{1}{2}(I + \Dd_Y\ii|\Dd_Y|) =
\frac{1}{2}\left(\begin{array}{cc}
           I & g_{+}^{-1} \\
           g_{+} & I
         \end{array}\right) \ .
\end{equation*}
The global boundary problem $D_{\Pi_{\geq}}$ is self-adjoint and,
more generally, a boundary condition $P\in Gr_{\o}(D)$ such that
$D_{P}$ is self-adjoint is characterized by having range equal to
the graph of an $L^2$-unitary isomorphism $T :F^+ \to F^-$ such
that $ T - g_+ $ has a smooth kernel \cite{Sc95}. Thus each
self-adjoint boundary condition $P=P_T$ defines a point $\det (T)$
of the determinant line
\begin{equation}\label{e:boundarydetline}
\Det(T) = \Det(g_+) \cong \Det((\textsf{D}^{+}_Y
\textsf{D}^{-}_Y)^{-1/2})\otimes \Det(\textsf{D}^{+}_Y) \cong
\Det(\textsf{D}^{+}_Y)
\end{equation}
of the boundary chiral Dirac operator $\textsf{D}^{+}_Y$. The
first isomorphism in \eqref{e:boundarydetline} is a general
functorial property of determinant lines under composition of
Fredholm operators \cite{Sc00}, while the second is defined
through the $\z$-determinant $${\rm det}_{\z} :
\Det(\textsf{D}^{+}_Y \textsf{D}^{-}_Y) \too \C \ .$$ That this
map is a {\em linear} isomorphism is a consequence of
\propref{p:multzeta}. A consequence of \eqref{e:boundarydetline}
and \thmref{t:releta} is that the $\eta$-invariant defines a
canonical linear isomorphism
\begin{equation}\label{e:etaboundarydetline}
e^{2\pi i\wt{\eta}(D)} :\Det(\textsf{D}^{+}_Y) \too \C \ ,
\end{equation}
via \eqref{e:boundarydetline} and the assignment $T\mtoo e^{2\pi
i\wt{\eta}(D_{P_T})}$. That is:
\begin{prop}
Absolutely---without a choice of boundary condition---the
exponentiated eta-invariant is a canonical element of the dual
determinant line of the boundary Dirac operator\footnote{The was
first pointed out by Segal \cite{Se90}.}
\begin{equation}\label{e:etaindual}
e^{2\pi i\wt{\eta}(D)}\in\Det(\textsf{D}^{+}_Y)^{*} \ .
\end{equation}
\end{prop}

To see this, first observe since $T$ is unitary that
\eqref{e:Uw3} becomes
\begin{equation*}
P_T  = \frac{1}{2}\begin{pmatrix}
    I & T\ii \\
    T & I
  \end{pmatrix} \ .
\end{equation*}
In particular, this holds for $P(D)$ for some unique unitary
$K:F^+ \to F^-$. It does not quite hold for $D-\la$ since the
operator is not of product type, but it is still true that
$H(D-\la)$ is the graph of a $\pdo$ operator $K_{\la}:F^+ \to F^-$
of order $0$, though not that $K_{\la}$ is an isometry or that
$K_{\la} - g_+$ is smoothing.

Consider two `self-adjoint' boundary conditions $P_1 = P_{T_1},
P_2 = P_{T_2} \in Gr_{\o}(D)$. The spectrum of the operators
$D_{P_i}$ is real and unbounded and, as in $\S$ 2.2, we denote
the two choices for $\th$ by $\pm$.
\begin{thm}\cite{ScWo99} For self-adjoint global boundary
problems $D_{P_1}, D_{P_2}$ for the Dirac operator over an
odd-dimensional spin manifold
\begin{equation}\label{e:SW}
 \frac{{\rm det}_{\z,\pm}(D_{P_1})}{{\rm det}_{\z,\pm}(D_{P_2})}
  =  \frac{{\rm
det}_F \lrb\frac{1}{2}(I + (T_1\ii K)^{\mp 1})\rrb}{{\rm
det}_F\lrb\frac{1}{2}(I + (T_2\ii K)^{\mp 1})\rrb} \ .
\end{equation}
Equivalently, if $P=P_T$
\begin{equation}
{\rm det}_{\z,\pm}(D_{P}) = {\rm det}_{\z,\pm}(D_{P(D)})\, . \,
{\rm det}_F\lsb\frac{1}{2}(I + (T\ii K)^{\mp 1})\rsb \ .
\end{equation}
\end{thm}
\begin{proof}
The equality \eqref{e:stiefelversion} becomes
\begin{eqnarray}
 \frac{{\rm det}_{\z,\pm }(D_{P_1})}{{\rm det}_{\z,\pm }(D_{P_2})}
 & = &   {\rm det}_F \lrb\frac{I +
T_1\ii K}{I + T_2\ii K}\rrb \, . \, \exp\lsb -{\rm
lim}^{\pm}_{\la\to\infty}\log {\rm det}_F \lrb\frac{I + T_1\ii
K_{\la}}{I + T_2\ii K_{\la}}\rrb \rsb \nonumber \\
 & = & \frac{{\rm
det}_F \lrb\frac{1}{2}(I + T_1\ii K)\rrb}{{\rm
det}_F\lrb\frac{1}{2}(I + T_2\ii K)\rrb} \, . \, \exp\lsb -{\rm
lim}^{\pm }_{\la\to\infty}\log {\rm det}_F \lrb\frac{I + T_1\ii
K_{\la}}{I + T_2\ii K_{\la}}\rrb \rsb \ \label{e:relsadet}.
\end{eqnarray}
The only extra subtlety introduced by $P_{F^+}\notin Gr_{\o}(D)$
 is that it is only the quotient of
operators $(I + T_1\ii K)/(I + T_2\ii K )$ which has a Fredholm
determinant. But since $T_i\ii K$ is of determinant class then so
is $(1/2)(I + T_1\ii K)$.  From \cite{Wo99} we have that
$D_{P_1}, D_{P_2}$  are strongly $\z$-comparable  and hence ${\rm
LIM}$ becomes the usual $\lim $ ($\S$ 2.4).

Finally,  either directly, using $g_{\la}\ii (D-\la)_{|U}g_{\la}
= D_{|U}$ with $g_{\la} = e^{-iu\la}\oplus e^{iu\la}$ in the
collar $U$, or using the symmetry argument of \cite{ScWo99}, one
has
\begin{equation}\label{e:Kla}
K_{\la} \to 0  \;\;  \text{as} \;\; \la\to \o \;\; \text{on} \;\;
R_{-} \ , \hskip 10mm K_{\la}^{-1} \to 0  \;\;  \text{as} \;\;
\la\to \o \;\; \text{on} \;\; R_{+} \ .
\end{equation}
The conclusion then follows from \eqref{e:relsadet}.
\end{proof}

From \eqref{e:SW}, switching the spectral cut conjugates the
relative zeta determinant. This corresponds to the equivalent
descriptions of $\ran(P_T)$ as $\gr(T :F^+ \to F^-)$ or $\gr(T\ii
:F^- \to F^+)$\footnote{Like the $\z$-determinant, the `canonical
determinant' of \cite{Sc95,ScWo99} is therefore not quite
canonical. The only completely canonical boundary determinant is
the quotient \eqref{e:relCdet}, which, like the relative
$\z$-determinant, has no `parity anomaly'.}. More generally,
allowing $D_{P_i}$ to be non-invertible, this disparity derives
from the relative eta-invariant:
\begin{thm}
\begin{equation}\label{e:relsaeta}
\wt{\eta}(D_{P_1}) - \wt{\eta}(D_{P_2})
  = \frac{1}{2\pi i}\log{\rm det}_F (T_2\ii T_1) \hskip 3mm \mod(\Z) \ .
\end{equation}
\end{thm}
\begin{proof}
From \cite{Wo99}, $D_{P_1}, D_{P_2}$  are strongly $\z$-comparable
and $\z(D_{P_1}^2,D_{P_2}^2,0)=0$. Hence from \eqref{e:eta},
\eqref{e:graphscattering} we have $\mod(2\pi i\Z)$
\begin{eqnarray*}
2\pi i \wt{\eta}(D_{P_1},D_{P_2}) & = & {\rm lim}_{\a\to
+\infty}\left(\log {\rm det}_F\Ss_{-i\a} - \log{\rm
det}_F\Ss_{i\a}\right)  \\  & = & {\rm lim}_{\a\to
+\infty}\left[\log {\rm det}_F \lrb\frac{I + T_1\ii K_{-i\a}}{I +
T_2\ii K_{-i\a}}\rrb - \log {\rm det}_F \lrb\frac{I + T_1\ii
K_{i\a}}{I + T_2\ii
K_{i\a}}\rrb\right]  \\
& = & \log {\rm det}_F \lrb T_2\ii T_1\rrb \;\;\mod(2\pi i\Z)
\end{eqnarray*}
where the final equality follows from \eqref{e:Kla}. Since
$D_{P_1}, D_{P_2}$ are $\z$-admissible, then \eqref{e:etasequal}
completes the proof.
\end{proof} The identity \eqref{e:relsaeta} is deduced in
\cite{KiLe00} from \eqref{e:SW}. Notice that \eqref{e:relsaeta}
is independent of $K:F^+ \to F^-$. More precisely, from
\eqref{e:ratio}, equation \eqref{e:relsaeta} is the assertion
that \eqref{e:etaboundarydetline} is linear.

\begin{rem}
(1) \ For a smooth family of $\z$-admissible operators
$\det_{\z}$ defines a section of the dual determinant line
bundle, that is, an element of Fock space. For a family of of
self-adjoint boundary problems $D_P$, $\det_{\z}$ defines a an
element of the Fock space associated to  $X$ (this is
\eqref{e:SW}), while the exponentiated $\eta$-invariant defines
an element of the boundary Fock space (this is
\eqref{e:etaindual},\eqref{e:relsaeta}).

\noi (2) \ The extension to the case where $\textsf{D}^{+}_Y$ is
non-invertible is easily done by augmenting $\textsf{D}^{+}_Y$ by
a unitary isomorphism $\s : \Ker(\textsf{D}^{+}_Y) \to
\Ker(\textsf{D}^{- }_Y)$. In particular, Thm(2.21) of
\cite{Mu94}, Thm(3.1) of \cite{LeWo96}, are special cases of
\eqref{e:relsaeta}.
\end{rem}

It is worth pointing out that, since $\ran(P_T\pp) = \gr(-T\ii :
F^{-} \to F^{+})$,
 if $M = X\cup X^{'}$ is a closed manifold with Dirac operator $A$
with $A_{|X}=D$ and $A_{|X^{'}}:=D^{'}$, then an easy corollary
of \eqref{e:relsaeta} is the

\vskip 1mm

\noi  {\it Weak Splitting Theorem : $\eta(D_{P_T}) +
\eta(D^{'}_{P_T^{\perp}})$  is constant as $P_T$ varies.}

\vskip 1mm

\noi The hard splitting Theorem \cite{BrLe98,KiLe00,Mu94,Wo99}
asserts this constant is precisely $\eta(A)$.

Another way of viewing the conjugation of the relative zeta
determinant on taking the conjugate spectral cut is through the
following formula for the relative Laplacian:
\begin{prop}
\begin{equation*}
\frac{{\rm det}_{\z,\pi}(D^2_{P_1})}{{\rm
det}_{\z,\pi}(D^2_{P_2})} = \left| \frac{{\rm
det}_{\z,\pm}(D_{P_1})}{{\rm det}_{\z,\pm}(D_{P_2})}
   \right|^2 =   \left|\frac{{\rm
det}_F \lrb\frac{1}{2}(I + T_1\ii K\rrb}{{\rm
det}_F\lrb\frac{1}{2}(I + T_2\ii K)\rrb} \right|^2 \ .
\end{equation*}
\end{prop}
\begin{proof}
Immediate from \eqref{e:reldetsq}, \eqref{e:graphscattering},
\eqref{e:relsaeta}.
\end{proof}
This formula is a special case of Theorem A to which we turn next.
(See also \remref{r:Aremark}(2)).


\section{An Application to the Laplacian on a Manifold with Boundary}

Let $X$ be an $n$-dimensional $\Ci$ compact Riemannian manifold
with boundary $Y$ and let $D : \Ci(X,E^1)\too \Ci(X, E^2)$ be a
Dirac-type operator with product case geometry, so that
\begin{equation}\label{e:collarD}
  D_{|U} = \sigma\left(\frac{\dd}{\dd u} + \Dd_Y\right),
\end{equation}
in a collar $U = [0,1)\times Y$ of the boundary, with notation as
in \eqref{e:collar}.

For each well-posed boundary condition $B$ for $D$ the associated
Dirac Laplacian
\begin{equation*}
\D_B = D_{B}^* D_B = D^* D : \dom(\D_B)\to L^2 (X,E^1)
\end{equation*}
with domain
\begin{equation*}
\dom(\D_B) = \{\psi\in H^2 (X,E^1) \ | \ B\g_0 \psi =0, B\st\g_0
D\psi = 0 \}
\end{equation*}
is a closed self-adjoint and positive operator on $L^2 (X, E^1)$
with discrete non-negative real spectrum.

The following result of Grubb allows us to define the zeta
determinant of $\D_B$.
\begin{prop}\cite{Gr99,Gr99a}
If $B$ is an admissible well-posed boundary condition for $D$,
then $\D_B$ is $\z$-admissible with spectral cut $R_{\pi}$.
\end{prop}
More precisely, Grubb proves that there is a resolvent trace
expansion for $m>n/2$ as $\la \to\o$  in closed subsectors of
$\C\backslash\ol{\R}_{+}$
\begin{equation}\label{e:Deltaexpansion}
\Tr\lrb\dd_{\la}^{m}(\D_B - \la)\ii\rrb \sim \sum_{j=-n}^{-1}a_j
(-\la)^{-j/2 -m - 1} + \sum_{j=0}^{\o}(a_{j,k} \log (-\la) + c_j)
(-\la)^{-j/2 -m - 1}
\end{equation}
and hence that the $\zeta$-function $\z(\D_B , s)$ defined by the
standard trace $\Tr(\D_B\si)$ for $\re(s) > n/2$ extends
meromorphically to all of $\C$ with the singularity structure
\begin{equation}\label{e:Deltapoles}
\G(s)\z(\D_B , s) \sim \sum_{j=-n}^{-1}\frac{\t{a}_j}{ s+ k/2} +
\frac{\dim\ker(\D_B)}{s} + \sum_{j=0}^{\o}\lrb\frac{\t{a}_j}{ (s+
k/2)^2} + \frac{\t{c}_j}{ s+ k/2}\rrb \ ,
\end{equation}
where $\th=\pi$ and the coefficients in \eqref{e:Deltapoles}
differ from those in \eqref{e:Deltaexpansion} by universal
constants. If $B-P(A)$ is a $\pdo$ of order $\leq - n$ then the
coefficient $\t{a}_0$ vanishes and so $\z(\D_B,s)$ is then
regular at $s=0$ and
\begin{equation*}
{\rm det}_{\z}\D_B = \exp(-\z^{'}(\D_B , 0))
\end{equation*}
is well-defined. In particular, $\det_{\z}\D_P$ exists for all
$P\in Gr_{\o}(D)$.

On the other hand, setting $S(P) := S_{D}(P)$, we have from
\propref{p:differbysmoothing} that the boundary `Laplacian'
\begin{equation*}
S(P)^* S(P) = P(D)^* \cdot P \cdot P(D) : H(D) \too H(D)
\end{equation*}
is of determinant class for all $P\in Gr_{\o}(D)$.

The main purpose of this section is to prove the following
Theorem.

\begin{thm}\label{t:thmA} Let $B_1, B_2$ be admissible well-posed
boundary conditions for a Dirac-type operator
$D:\Ci(X,E^1)\too\Ci(X,E^2)$. Then, with $P_i = P[B_i]$, one has
\begin{equation}\label{e:metrics1}
\frac{{\rm det}_{\z}(\D_{B_1})}{{\rm det}_{\z}(\D_{B_2})} =
\frac{{\rm det}_{F}\lrb S(P_1)^* S(P_1 )\rrb }{{\rm det}_{F}\lrb
S(P_2)^* S(P_2 )\rrb } \ .
\end{equation}
Or, from \eqref{e:zetaFred2},
\begin{equation}\label{e:metrics2}
{\rm det}_{\z}(\D_{B_1},\D_{B_2}) = {\rm det}_{\z}( S(P_1)^*
S(P_1) \  ,\  S(P_2)^* S(P_2) ) \ .
\end{equation}
 Equivalently, since $S(P(D)) = Id$,
 \begin{equation}\label{e:metrics3}
{\rm det}_{\z}(\D_{B}) = {\rm det}_{\z}(\D_{P(D)}).{\rm
det}_{F}\lrb S(P[B])^* S(P[B])\rrb \ .
\end{equation}
\end{thm}

\vskip 5mm

\begin{rem}\label{r:Aremark}
\noi (1) Because of \lemref{lem:adsmoothing} and \eqref{e:PB} it
is sufficient to assume that $B_i = P_i \in Gr_{\o}(D)$, and from
here on that is what we shall do.

\noi (2) In Stiefel graph coordinates \eqref{e:metrics1} has the
form
\begin{equation*}
\frac{{\rm det}_{\z}(\D_{P_1})}{{\rm det}_{\z}(\D_{P_2})}  =
\frac{{\rm det}_{F}\lrb Q_1\ii \lrb I + T^{*}_1 K\rrb \lrb I +
K^{*}T_1\rrb \rrb }{{\rm det}_{F}\lrb Q_2\ii \lrb I + T^{*}_2
K\rrb \lrb I + K^{*}T_2\rrb \rrb } =  e^{\log{\rm det}_F (Q_1\ii
Q_2)}\left|\frac{{\rm det}_{F}\lrb I + T^{*}_1 K\rrb}{{\rm
det}_{F}\lrb I + T^{*}_2 K\rrb} \right|^2 .
\end{equation*}

\noi (3) The simple form of \eqref{e:metrics1} depends on the
homogeneous structure of $Gr_{\o}(D)$, it does not persist to
more general classes of well-posed boundary conditions.

\noi (4)  We may replace $P(D)$ by $P(D)_{{\rm ort}}$ (cf.
\eqref{e:ort}) in \thmref{t:thmA}). This follows from the
invertibility of $P(D)(P(D)P(D)^* + (I-P(D)^*)(I-P(D)))\ii P(D)^*
P(D)$ on $H(D)$, which is therefore not detected in the quotient
on the right-side of \eqref{e:metrics1}.

 \noi (5) Implicit in \thmref{t:thmA} is the invertibility of
$D_{B_i}$. We assume invertibility when obviously required
without further mention. The identity \eqref{e:metrics3} is
globally defined.

\noi (6) The identifications hold for $P_1 - P_2$ differing just
by a $\pdo$ of order  $ < -n$.
\end{rem}

\subsection{Proof of \thmref{t:thmA}}

To identify the scattering operator we use a canonical
identification of the solution space of $\D_P$ with that of an
associated first-order elliptic system.

\vskip 2mm

$\bb$ {\bf An equivalent first-order elliptic system}

\vskip 2mm

We analyze $\D_P = D^{*}_P D_P$ through the first-order elliptic
operator acting on sections of $E^1 \oplus E^2$
\begin{equation*}
  \wD = \begin{pmatrix}
    0 & D^* \\
    D & -I \
  \end{pmatrix} : H^1 (X; E^1 \oplus E^2 )\to L^2 (X; E^1 \oplus
  E^2
  ) \ .
\end{equation*}
$$\wD(s_1,s_2) = (D^* s_2, Ds_1 - s_2) \ .$$ From \eqref{e:collar}
and \eqref{e:collar*} we find that $\wD$ is of Dirac-type with
\begin{equation*}
  \wD_{|U} = \ws\left(\frac{\dd}{\dd u} + \widehat{\Aa}_Y + \wR\right) ,
\end{equation*}
where
\begin{equation}\label{e:ws}
  \ws = \begin{pmatrix}
    0 & -\sigma\ii \\
    \sigma & 0 \
  \end{pmatrix} , \;\;\;\;\;\widehat{\Aa}_Y = \begin{pmatrix}
    \Dd_Y &  0\\
    0 &  -\sigma \Dd_Y\sigma\ii \
  \end{pmatrix}  ,
  \;\;\;\;\;\wR = \begin{pmatrix}
    0 & -\sigma\ii \\
    0 & 0 \
  \end{pmatrix}  ,
\end{equation}
satisfying the relations
\begin{equation}\label{e:relations}
\ws^2 = - I \ , \hskip 7mm \ws^* = - \ws \ , \hskip 7mm
\ws\widehat{\Aa}_Y  + \widehat{\Aa}_Y\ws = 0 \ , \hskip 7mm \ws\wR
+ \wR\ws = -I \ .
\end{equation}
Green's Theorem for the formally self-adjoint operator $\wD$ now
states that
\begin{equation}\label{e:wGreens}
  <\wD s_1 , s_2 > - <s_1, \wD s_2> = <-\ws\g_0 s_1, \g_0 s_2 > \ ,
\end{equation}
where here $\g_0 (\psi,\phi) = (\g_0\psi,\g_0\phi)$.

Setting $A=\wD$ in our discussion in $\S 4$, we have a Poisson
operator for $\wD$
\begin{equation}\label{e:wPoisson}
  \wKk := \Kk_{\wD} : H^{s-1/2}(Y, (E^1 \oplus E^2)_{|Y}) \too
  \Ker(\wD,s) \subset H^{s}(X, E^1 \oplus E^2),
\end{equation}
and Calderon projector
\begin{equation*}
P(\wD) = \g_0\wKk \ .
\end{equation*}
We can compute $P(\wD)$ quite explicitly:
\begin{lem}\label{lem:wCalderon}
\begin{equation}\label{e:wCalderon}
P(\wD) = \begin{pmatrix}
  P(D) & \g D_{P(D)}\ii \Kk_{*} \\
  0 & P(D^*)
\end{pmatrix} \ ,
\end{equation}
where $\Kk_{*}$ is the Poisson operator for $D^*$.
\end{lem}
\vskip 2mm \noi We postpone the proof for the moment. Notice,
however, since $\wR\neq 0$, that $P(\wD) - P(D)\oplus P(D^*)$ is
only a $\pdo$ of order $-1$ and not smoothing due to the
off-diagonal term. Further, it is a projector but not a projection
(cf. \remref{r:Aremark}(4)).

Since $\wD$ is of Dirac-type, we have a $\pdo$ Grassmannian
$Gr_{1}(\wD)$ of global boundary conditions for $\wD$, and for
each $Q\in Gr_{1}(\wD)$ a first-order global boundary problem
\begin{equation*}
\wD_{Q} = \wD :\dom(\wD_{Q}) \to L^2(X,E^1 \oplus E^2 ) \ .
\end{equation*}

We recover the resolvent $(\D_P -\la)\ii$ in the following way.
First, we have a canonical map
\begin{equation*}
  Gr_{\o}^{(r)}(D) \too Gr_{1}^{(0)}(\wD) \ ,
\end{equation*}
\begin{equation*}
  P\mtoo \wP := P \op P\st \ .
\end{equation*}
To see that $\wP$ is in the index zero component
$Gr^{(0)}_{1}(\wD)$, observe that the identity $I : \widehat{H}_Y
\to \widehat{H}_Y$ acting between the block decompositions
$H(D)\oplus H(D)\pp$, $W\oplus W\pp$, where $W=\ran(P)$, of
$\widehat{H}_Y = L^2 (Y, (E^1 \oplus E^2 )_{|Y})$ is
\begin{equation}\label{e:identity}
  I = \begin{pmatrix}
    S(P) & S\pp(P) \\
    S(P\pp) & S\pp(P\pp) \
  \end{pmatrix} \ .
\end{equation}
For any $\pdo$ $B$ on $H_Y$ set
\begin{equation*}
S\pp(B) = B\circ P(D)\pp : H(D)\pp \too \ran(B)  \ ,
\end{equation*}
and for $P\in Gr_{\o}(D)$ note that
\begin{equation}\label{e:star}
 S\st(P\st) := P\st\circ
P(D^*) = \s S\pp (P\pp) \s\ii  : \s H(D)\pp \too \s W\ ,
\end{equation}
and let
\begin{equation*}
S(\wP) = \wP\circ P(\wD) : H(\wD) \to \ran(\wP) = W\oplus W\pp \ ,
\end{equation*}
where we use \eqref{e:S(P^*)}. Then from \eqref{e:wCalderon} and
\eqref{e:identity}
\begin{equation*}
\ind(S(\wP)) =  \ind ( S(P)) + \ind(S\pp(P\pp))
   =  \ind \lrb I - \begin{pmatrix}
    0 & S\pp(P) \\
    S(P\pp) & 0 \\
  \end{pmatrix}\rrb
\end{equation*} which is zero, since the matrix operator is a
$\pdo$ of order $-1$ and hence compact. Alternatively, this fact
follows from $\ind(S(\wP))  = \ind(\wD_{\wP})$ and the identity
\begin{equation*}
\s \wP\pp \s\ii = \wP \ ,
\end{equation*}
which along with  \eqref{e:relations}, \eqref{e:wGreens} implies
that $\wD_{\wP}$ is self-adjoint considered as a closed operator
on $L^2 (X, E^1 \oplus E^2)$.

Next, we have a canonical inclusion defined by $D$
\begin{equation*}
  \widehat{i}: H^1 (X,E^1 ) \too L^2(X; E^1 \oplus E^2 ) \ ,
\hskip 8mm \widehat{i}(\psi) = (\psi, D\psi) \ .
\end{equation*}
Setting for $\la\in\C$
\begin{equation*}
\wDla = \begin{pmatrix}
    -\la & D^* \\
    D & -I \
  \end{pmatrix} : H^1 (X; E^1 \oplus E^2 )\to L^2 (X; E^1 \oplus
  E^2
  ) \ ,
\end{equation*}
the inclusion $\widehat{i}$ restricts to an isomorphism
\begin{equation*}
\widehat{i}_{|\Ker} : \Ker(\D - \la) \stackrel{\simeq}{\too}
\Ker(\wDla) \subset H^2 (X,E^1) \oplus H^1 (X, E^2) \ ,
\end{equation*}
with inverse $(s_1,s_2)\mto s_1$, where $\D - \la, \wDla$ are
acting in $H^2 (X,E^1), H^1 (X,E^1\oplus E^2)$, respectively. That
$\widehat{i}_{|\Ker}$ is injective with range $\Ker(\wDla)$
follows from the identity
\begin{equation}\label{e:kernels}
\wDla \begin{pmatrix}
    \psi \\
    D\psi \
\end{pmatrix} =
  \begin{pmatrix}
    (\D - \la)\psi \\
    0 \
\end{pmatrix} \ .
\end{equation}
On the other hand, if $(s_1,s_2)\in \Ker(\wDla)$, then $D^* s_2 =
\la s_1$ and $s_2 = Ds_1$ and hence $s_1 \in \Ker(\D-\la)$. In
particular, setting $s_i = \widehat{i}(\psi_i)$ we can extract
Green's formula for $\D$ from \eqref{e:wGreens} and
\eqref{e:kernels} (with $\la=0$) :
\begin{equation*}
  <\D\psi_1,\psi_2> - <\psi_1,\D\psi_2> =
  \left< -\ws\g_0 \begin{pmatrix}
    \psi_1 \\
    D\psi_1 \
\end{pmatrix}, \g_0 \begin{pmatrix}
    \psi_2 \\
    D\psi_2 \
\end{pmatrix}  \right> \ .
\end{equation*}

The operator $\widehat{i}$ also restricts to a canonical inclusion
\begin{equation}\label{e:domains}
  \widehat{i} : \dom(\D_P) \too \dom(\wD_{\wP}) \ .
\end{equation}
From \eqref{e:kernels} and \eqref{e:domains} we have for
$\la\in\C\backslash\ol{\R}_{+}$
\begin{equation}\label{e:resolvents}
(\D_P - \la)\ii = \left[  \wD_{\la,\wP}\ii \right]_{(1,1)} :
L^{2}(X,E^1) \too \dom(\D_P-\la) \ ,
\end{equation}
where for an operator $C = \begin{pmatrix}
  S & T \\
  U & V
\end{pmatrix}$ on $L^2 (X,E^1 \oplus E^2)$, we define $[C]_{(1,1)}
= S$. Equivalently,
\begin{equation}\label{e:resolvents2}
(\D_P - \la)\ii = \begin{pmatrix}
  I & 0
\end{pmatrix} \wD_{\la,\wP}\ii \begin{pmatrix}
    I \\
    0 \
\end{pmatrix}  \ .
\end{equation}
A precise formula for $\wD_{\la,\wP}\ii $ is given in
\eqref{e:id4}.

\vskip 5mm

\noi $\bb$ \ {\bf The scattering determinant}

\vskip 2mm

Let $P_1, P_2 \in Gr_{\o}(D)$ and for
$\mu\in\C\setminus\ol{\R}_{+}$ set
\begin{equation*}
S_{\mu}(\wP_i) := \wP_i\circ P(\wD_{\mu}) : H(\wD_{\mu}) \too
\ran(\wP_i) \ .
\end{equation*}
Let $\Ee : \ran(P_2) \to \ran(P_1), \t{\Ee} : \ran(P_2\st) \to
\ran(P_1\st)$ be auxiliary operators for $D_{P_i},D_{P_i}^*$
respectively. Then
\begin{equation}\label{e:wE}
\wE= \begin{pmatrix}
  \Ee &  0\\
  0 &   \t{\Ee}
\end{pmatrix}
:\ran(\wP_2)\to \ran(\wP_1)
\end{equation}
is an auxiliary operator for $\wD_{\wP_1},\wD_{\wP_2}$ and we
have:
\begin{prop}\label{p:Deltascatter}
 $(\D_{P_1},\D_{P_2})$ are
$\z$-comparable with scattering determinant
 $${\rm det}_{\wE}(\wD_{\mu,\wP_1}, \wD_{\mu,\wP_2}) =
  {\rm det}_{F}((\wE S_{\mu}(\wP_2))\ii S_{\mu}(\wP_1))$$ taken on
$H(\wD_{\mu})$.  With $\th = \pi$ and $\la\in\R_+$, one has
\begin{equation}\label{e:relLapdet}
 \frac{{\rm det}_{\z}(\D_{P_1})}{{\rm
det}_{\z}(\D_{P_2})} = {\rm det}_{F}\left(\frac{S(\wP_1)}{\wE S
(\wP_2)}\right) \, . \,
 e^{-{\rm LIM}_{\la\to +\infty}\log\det_F ((\wE S_{-\la} (\wP_2))\ii S_{-\la}
 (\wP_1))}\; .
\end{equation}
If $(P_2,P_1)$ is invertible
\begin{equation}\label{e:relLapdet2}
 \frac{{\rm det}_{\z}(\D_{P_1})}{{\rm
det}_{\z}(\D_{P_2})} = {\rm det}_{F}\left(\frac{S(\wP_1)}{\wP_1 S
(\wP_2)}\right) \, . \,
 e^{-{\rm LIM}_{\la\to +\infty}\log\det_F ((\wP_1 S_{-\la} (\wP_2))\ii S_{-\la}
 (\wP_1))}\; .
\end{equation}
\end{prop}
\begin{proof}
From \cite{Gr99} Cor(9.5) ,\cite{Gr99a} Thm(1), the coefficients
$a_{j,k}, a_j$ in the asymptotic expansion
\eqref{e:Deltaexpansion} are locally determined by the symbols of
$\D$ and $B$, while, provided $P_1-P_2\in \Psi_{l}(H_Y)$ with
$l\geq n$, the expansion coefficients differ only in the $c_j$.
Integrating we hence obtain a resolvent trace expansion in closed
subsectors of $\C\backslash\ol{\R}_{+}$
\begin{equation}\label{e:RelLapAsymp}
\Tr((\D_{P_1} - \mu)\ii-(\D_{P_2} - \mu)\ii)\sim
\end{equation}
$$\sum_{j=1}^{\o}\sum_{k=0}^{1}C_{j,k}(-\mu)^{-j/2- 1} \log(-\mu)
+ \z(\D_{P_1},\D_{P_2},0)(-\mu)\ii$$
 where the coefficients
$C_{j,k}=C_{j,k}(\D,P_1,P_2)$ differ from the $c_j$ by universal
constants.

Since $\wP_1 - \wP_2$ has a smooth kernel we know from
\eqref{e:zrelinv} that so does $\wD_{\la,\wP_1}\ii -
\wD_{\la,\wP_2}\ii$, and from \eqref{e:resolvents} also
$\D_{P_1}\ii - \D_{P_2}\ii$. From \eqref{e:resolvents},
\eqref{e:resolvents2} we have
\begin{eqnarray*}
\Tr\left((\D_{P_1} - \mu)\ii - (\D_{P_2 }- \mu)\ii\right) & = &
\Tr\left(\left[ \wD_{\mu,\wP_1}\ii - \wD_{\mu,\wP_2}\ii
\right]_{(1,1)}\right) \\ & = & \Tr\left(\begin{pmatrix}
    I & 0 \
\end{pmatrix} (\wD_{\mu,\wP_1}\ii -
\wD_{\mu,\wP_2}\ii) \begin{pmatrix}
    I \\
    0 \
\end{pmatrix}\right) \\
& = & \Tr\left(\begin{pmatrix}
    I & 0 \\
    0 & 0
\end{pmatrix} (\wD_{\mu,\wP_1}\ii -
\wD_{\mu,\wP_2}\ii)\right) \\ & = &
-\Tr\left(\frac{\dd}{\dd\mu}(\wD_{\mu})\left(\wD_{\mu,\wP_1}\ii -
\wD_{\mu,\wP_2}\ii\right)\right)  \\ & = &
-\frac{\dd}{\dd\mu}\log{\rm det}_F \left(\frac{S_{\mu}(\wP_1)}{\wE
S_{\mu}(\wP_2)}\right) \ ,
\end{eqnarray*}
where we use \eqref{e:variation2} for the final equality, since
the variation in $U$ is of order $0$.

Hence $(\D_{P_1},\D_{P_2})$ are $\z$-comparable. Since they are
also $\z$-admissible, \eqref{e:relLapdet} is a consequence of
\thmref{t:thmB}.
\end{proof}

\vskip 2mm

\noi $\bb$ {\bf Relation to the right-side of \eqref{e:metrics1}}

\vskip 2mm

\begin{prop}\label{p:alphaD}
Let $S_r(P_i)$ be the boundary integrals defined by a smooth
1-parameter family of Dirac-type operators $D_{r}$. Then
\begin{equation}\label{e:bigScatter3}
 \frac{d}{dr}\log{\rm det}_F\lrb\frac{S_r(\wP_1)}{\wE S_r
(\wP_2)}\rrb =   \frac{d}{dr}\log\frac{{\rm det}_{F}\lrb
S_r(P_1)^* S_r(P_1 )\rrb }{{\rm det}_{F}\lrb S_r(P_2)^* S_r(P_2
)\rrb } \ .
\end{equation}
If $(P_2,P_1):\ran(P_2)\to\ran(P_1)$ is invertible, one has
\begin{equation}\label{e:bigScatter2}
 {\rm det}_{F}\lrb\frac{S(\wP_1)}{\wP_1 S
(\wP_2)}\rrb =  \frac{{\rm det}_{F}\lrb S(P_1)^* S(P_1 )\rrb
}{{\rm det}_{F}\lrb S(P_2)^* S(P_2 )\rrb } \, . \, \frac{1}{{\rm
det}_{F}\lrb P_2 P_1 P_2 \rrb } \ ,
\end{equation}
the determinant in the denominator being taken on $\ran(P_2)$.
\end{prop}

\vskip 1mm

\begin{proof}
Equation \eqref{e:bigScatter2} is a consequence of the following
identity.
\begin{lem}\label{lem:id}
\begin{equation}\label{e:id}
  {\rm det}_{F}\left(\frac{S\pp(P_1\pp)}{P_1\pp S\pp (P_2\pp)}\right)  =
   \ol{{\rm det}_{F}\left(\frac{S(P_1)}{P_1 S (P_2)}\right)}  \ ,
\end{equation}
where the left-hand determinant is taken on $H(D)\pp$.
\end{lem}
\begin{proof} As in \propref{p:stiefel}, we may choose graph Stiefel coordinates
\begin{equation}\label{e:graphs1}
  H(D) = \gr(K:E\to E\pp) \ , \hskip 8mm W_i = \gr(T_i:E\to
  E\pp) \ .
\end{equation}
Using the property
\begin{equation}\label{e:conjugate}
  \ol{{\rm det}_F A} = {\rm det}_F A^*
\end{equation}
for $A:E\to E$ of determinant class and letting ${\rm
det}_{[E]}(A)$ mean the Fredholm determinant taken on $E$,  we
have from \eqref{e:graphscattering2}
\begin{equation*}
 \ol{{\rm det}_{F}\left(\frac{S(P_1)}{P_1 S (P_2)}\right)} =
 \frac{{\rm det}_{[E]}(I + K^* T_1)}{{\rm det}_{[E]}(I + K^* T_2)}
 \frac{{\rm det}_{[E]}(I + T_{2}^* T_2)}{{\rm det}_{[E]}(I + T_{2}^*
 T_1)} \ .
\end{equation*}
From \eqref{e:graphs1} we have
\begin{equation*}
  H(D)\pp = \gr(-K^* :E\pp\to E) \ , \hskip 8mm W_i\pp = \gr(-T_{i}^* :E\pp\to
  E) \ ,
\end{equation*}
and so in Stiefel coordinates
\begin{equation}\label{e:orthogcoords}
 P(D)\pp = \begin{pmatrix}
   K^* \wQ_{K}\ii K &  -K^* \wQ_{K}\ii \\
   -\wQ_{K}\ii K &  \wQ_{K}\ii \
 \end{pmatrix} \ , \hskip 10mm  P_i \pp = \begin{pmatrix}
   T_{i}^* \wQ_{i}\ii T_{i} &  -T_{i}^* \wQ_{i}\ii \\
   -\wQ_{i}\ii T_{i} &  \wQ_{i}\ii \
 \end{pmatrix} \ ,
\end{equation}
where $\wQ_{K} = I +KK^*$, $\wQ_{i} = I +T_i T_{i}^*$. Using these
local representations we compute in a similar fashion to
\eqref{e:graphscattering2}
\begin{equation*}
 {\rm det}_{F}\left(\frac{S\pp(P_1\pp)}{P_1\pp S\pp (P_2\pp)}\right) =
 \frac{{\rm det}_{[E\pp]}(I + T_1 K^*)}{{\rm det}_{[E\pp]}(I + T_1 K^*)}
 \frac{{\rm det}_{[E\pp]}(I + T_2 T_{2}^* )}{{\rm det}_{[E\pp]}(I +
 T_1 T_{2}^*)} \ .
\end{equation*}
Since ${\rm det}_{[E\pp]}(I + ST^*) =
 {\rm det}_{[E]}(I + T^* S)$ for any $S, T :E \to E\pp$ of trace
 class,  we reach the conclusion.
\end{proof}

From \lemref{lem:wCalderon}
\begin{equation*}
  S(\wP_1) = \begin{pmatrix}
    P_1 & 0 \\
    0 & P_{1}^* \
  \end{pmatrix}
  \begin{pmatrix}
    P(D) & \g_0 D_{P(D)}\ii \Kk_* \\
    0 & P(D^*) \
  \end{pmatrix}
  =\begin{pmatrix}
    S(P_1) & P_1\g_0 D_{P(D)}\ii \Kk_* \\
    0 &  \s S\pp (P_1\pp)\s\ii \
  \end{pmatrix} \ ,
\end{equation*}
where we use \eqref{e:star}. Computing $ \wP_1 S(\wP_2) $ in a
similar way, and using \eqref{e:id}, \eqref{e:conjugate} and the
multiplicativity of the Fredholm determinant, we obtain
\begin{equation}\label{e:scatter2}
   {\rm det}_{F}\left(\frac{S(\wP_1)}{\wP_1 S
(\wP_2)}\right)   =   {\rm det}_{F}\left(\frac{S(P_1)}{P_1 S
(P_2)}\right) . {\rm det}_{F}\left(\frac{S\pp(P_1\pp)}{P_1\pp S\pp
(P_2\pp)}\right)
\end{equation}
\begin{eqnarray}& = &  {\rm
det}_{F}\left(\left(\frac{S(P_1)}{P_1 S (P_2)}\right)^*
\frac{S(P_1)}{P_1 S (P_2)}\right) \nonumber\\ & = &  {\rm
det}_{F}\left(\frac{S(P_1)^* S(P_1)}{S (P_2)^* P_2 P_1 P_2 S
(P_2)}\right) \nonumber\\ & = &  \frac{{\rm det}_{F}\lrb S(P_1)^*
S(P_1)\rrb } {{\rm det}_{F}\lrb S(P_2)^* S(P_2 )\rrb }
 \, . \, \frac{1}{{\rm
det}_{F}\left( [S(P_2)^* S(P_2 )]\ii S(P_2)^* P_2 P_1 P_2 S
(P_2)\right)} \nonumber\\ & = &  \frac{{\rm det}_{F}\lrb S(P_1)^*
S(P_1 )\rrb }{{\rm det}_{F}\lrb S(P_2)^* S(P_2 )\rrb }
 \, . \, \frac{1}{{\rm
det}_{F}( P_2 P_1 P_2 )} \nonumber\ ,
\end{eqnarray}
which proves \eqref{e:bigScatter2}.

To see \eqref{e:bigScatter3}, first note that in the same way as
\eqref{e:scatter2} we have
\begin{equation}\label{e:scatter3}
  {\rm det}_{F}\left(\frac{S_r(\wP_1)}{\wE S_r
(\wP_2)}\right)  =   {\rm det}_{F}\left(\frac{S_r(P_1)}{\Ee S_r
(P_2)}\right) . {\rm det}_{F}\left(\frac{S_r\pp(P_1\pp)}{\Ee\pp
S_r\pp (P_2\pp)}\right) \ ,
\end{equation}
where $\Ee\pp := \s\t{\Ee}\s\ii$. Thus we have to show that
\begin{equation*}
\sum_{i=1}^2\Tr(S_r\pp(P_i\pp)\ii \frac{d}{d r}S_r\pp(P_i\pp))  =
 \sum_{i=1}^2\ol{\Tr}(S_r (P_i)\ii \frac{d}{d r}S_r (P_i))  \ ,
\end{equation*}
where the right-side means complex conjugate.  This follows by the
same method used in \lemref{lem:id}, via the Stiefel coordinate
representation for $S(P_i)\ii$ (use \eqref{e:Rinverse} with $T=K,
\widehat{T}=T_i$)  and its analogue for $S\pp(P_i\pp)\ii$ (use
\eqref{e:orthogcoords} ). Or, these coordinate matrices may be
used to prove directly, in a similar way to \propref{p:stiefel},
that \eqref{e:scatter3} differs from the right-side of
\eqref{e:metrics1} by a function independent of $P(D)$.
\end{proof}

\vskip 2mm

\begin{prop}\label{p:metrics6}
Let $D_{r}$, $-\e \leq r \leq \e$, be a 1-parameter family of
Dirac-type operators with product case geometry such that $\wD(r)
= \begin{pmatrix}
  0 & D^{*}_r \\
  D_r & -I
\end{pmatrix}$
satisfies \eqref{e:collarz}. Then $P_1, P_2 \in Gr_{\o}(D_0)$ are
global boundary conditions for $D_r$. If the $D_{r,P_i}$ are
invertible, then, with $$S_r (P_i) = P\circ
P(D_r):H(D_r)\to\ran(P) \ ,$$ one has
\begin{equation}\label{e:metrics6}
\frac{d}{dr}\log\frac{{\rm det}_{\z}(\D_{r,P_1})}{{\rm
det}_{\z}(\D_{r, P_2})} = \frac{d}{dr}\log\frac{{\rm det}_{F}\lrb
S_r(P_1)^* S_r(P_1 )\rrb }{{\rm det}_{F}\lrb S_r(P_2)^* S_r(P_2
)\rrb } \ .
\end{equation}
\end{prop}
\begin{proof} Let $S_r(\wP) = \wP\circ P(\wD(r)):H(\wD(r))\to\ran\wP$.
Then from \propref{e:variation2} we have
\begin{equation}\label{e:leftside2}
\Tr\left(\frac{d}{dr}(\wD(r))\left(\wD(r)_{\wP_1}\ii -
\wD(r)_{\wP_2}\ii\right)\right) = \frac{d}{dr}\log{\rm
det}_{F}\left(\frac{S_r(\wP_1)}{\wE S_r(\wP_2)}\right) \ .
\end{equation}
In view of \eqref{e:bigScatter3}, we need to prove that left-side
\eqref{e:metrics6} = left-side \eqref{e:leftside2}. We show each
expression is equal to
\begin{equation}\label{e:leftside1.1}
\Tr_{L^2}\lrb \dD (D_{P_1}\ii - D_{P_2}\ii)\rrb +  \Tr_{L^2}\lrb
\dD^{*}((D^{*}_{P_1})\ii - (D_{P_2}^*)\ii)\rrb \ ,
\end{equation}
where $\dD = (d/dr)D_r$. We omit $r$ from the operator notation
throughout.

First, for any $P,P_1,P_2\in\Gr_{\o}(D)$ we record the following
identities:
\begin{equation}\label{e:id1}
\dDD = \dD^* D + D^*\dD \ ,
\end{equation}
\begin{equation}\label{e:id2}
D\D_P\ii = (\D_{P}^*)\ii \ , \hskip 10mm D^*\t{\D}_P\ii = D_P\ii
\, (\D_{P}^*)\ii \ ,
\end{equation}
\begin{equation}\label{e:id3}
(D_{P_1}^*)\ii \D_{P_2}\ii = \t{\D}_{P_2}\ii (D_{P_1}^*)\ii \ ,
\end{equation}
\begin{equation}\label{e:id4}
\wD_{\la,\wP}\ii =
\begin{pmatrix}
  (\D_P - \la)\ii & D_{P}^* (\t{\D}_P - \la)\ii \\
 D_P (\D_P - \la)\ii & \la(\t{\D}_{P} - \la)\ii
\end{pmatrix} \ ,
\end{equation}
where $\tilde{\D} =  D^* D$, $\tilde{\D}_{P} = D_P D^{*}_{P}$. To
see \eqref{e:id2}, since $D^* D\D_P\ii = I$ on $L^2 (X,E)$, and
$D\D_P\ii $ has range in $\dom(D^{*}_P)$, one has (using
\eqref{e:master} for $D^{*}_P$)
\begin{equation*}
(D_{P}^*)\ii =  \lrb(D_{P}^*)\ii D^*\rrb D\D_P\ii = D\D_P\ii -
\Kk_{*}(P\st)\g_0 D\D_P\ii = D\D_P\ii \ .
\end{equation*}
The other identities can be checked in a similar fashion. For
brevity let $\D_i = \D_{P_i}, D_i = D_{P_i}, D^{*}_{i} =
D^{*}_{P_i} := D^{*}_{P_i\st}$.

Setting $\la = 0$ in \eqref{e:id4} we have
\begin{equation*}
\dot{\wD} = \begin{pmatrix}
   0 & \dD^* \\
   \dD & 0 \
 \end{pmatrix} \ , \hskip 10mm
 \wD_{\wP_i}\ii = \begin{pmatrix}
   \D_i\ii & D_i\ii \\
  (D^{*}_i)\ii & 0 \
 \end{pmatrix} \ ,
\end{equation*}
and hence $\dot{\wD}\lrb\wD_{\wP_1}\ii - \wD_{\wP_2}\ii\rrb =
\begin{pmatrix}
   \dD^* \lrb(D^{*}_1)\ii - (D^{*}_2)\ii\rrb & 0 \\
  \dD\lrb\D_1\ii - \D_2\ii\rrb &  \dD\lrb D_1\ii - D_2 \ii \rrb \
 \end{pmatrix},$
from which the equality of the left-side of \eqref{e:leftside2}
with \eqref{e:leftside1.1} is clear.

Next, let $\z_{rel}(0) = \z(\D_1,\D_2,0)$. The resolvent trace
\eqref{e:RelLapAsymp} implies a heat trace expansion as $t\to 0$
\begin{equation}\label{e:heattraceasymps}
  \Tr(e^{-t\D_1}- e^{-t\D_2})) \sim
\sum_{j=1}^{\o}\sum_{k=0}^{1}\t{C}_{j,k}t^{j/2} \log^k t +
\z_{rel}(0) \ ,
\end{equation}
while from \eqref{e:detaroundzero2} we have
\begin{equation*}
\log{\rm det}_{\z,\th}(\D_1,\D_2)  =  \lsb \int_{0}^{\o}t^{s-1}
  \Tr(e^{-t\D_1} - e^{-t\D_2}) \ dt -
  \frac{\z_{rel}(0)}{s}\rsb_{|s=0} - \g \z_{rel}(0) \ .
\end{equation*}
Since $\D_1\ii - \D_2\ii$ has a smooth kernel, precisely the same
argument as that leading to \eqref{e:relzetazerovar} yields $
(d/dr)\z_{rel}(0) = 0.$ So the $r$-variation `kills' the pole at
$s=0$. From \eqref{e:heattraceasymps} we therefore have
$(d/dr)\Tr(e^{-t\D_1}- e^{-t\D_2}))=O(t^{1/2})$, and hence
\begin{eqnarray}
\frac{d}{dr}\log{\rm det}_{\z,\th}(\D_1,\D_2)  \, =
&-&\int_{0}^{\o} t^{-1}
  \frac{d}{dr}\lsb \Tr(e^{-t\D_1}) - \Tr(e^{-t\D_2})\rsb \ dt \label{e:one}\\
 = &  & \int_{0}^{\o}
   \Tr(\dDD_1 e^{-t\D_1} - \dDD_2 e^{-t\D_2}) \ dt \nonumber
   \\
= \int_{0}^{\o}
   \Tr(\dD^* D (e^{-t\D_1}  -   e^{-t\D_2})
   & + &   \Tr(D^* \dD (e^{-t\D_1} - e^{-t\D_2})  \ dt
   \ ,\nonumber
\end{eqnarray}
 where in the second equality we use Duhamel's Formula and the
symmetry of the trace. The heat operator
\begin{equation*}
 e^{-t\D_i} = \frac{i}{2\pi}\int_{C_{\pi}}e^{-t\la}(\D_i - \la)\ii
 \, d\la : L^2 (X,E) \too \dom(\D_i) \,
\end{equation*}
has range in $\dom(\D_i)$, and hence it follows that $D^* \dD
e^{-t\D_i} = D_{i}^* \dD_i e^{-t\D_i}$, since $P_i\g_0\psi = 0,
P_i\st\g_0 D \psi = 0$ implies the domain $P_i\g_0\psi = 0,
P_i\st\g_0 \dD \psi = 0$ for $D^* \dD$. Thus, using also
\eqref{e:id3} and the contour integral definition of
$e^{-t\t{\D}_i}$, which imply $e^{-t\D_i}D_{i}^* =
D_{i}^*e^{-t\t{\D}_i}$, we have
\begin{equation}
\Tr\lrb D^* \dD e^{-t\D_i}\rrb =  \Tr\lrb D_{i}^* \dD_i
e^{-t\D_i}\rrb  =  \Tr\lrb D_{i}^* \dD_i
e^{-t\D_i}D_{i}^*(D_{i}^*)\ii\rrb
\end{equation}
$$ =  \Tr\lrb D_{i}^* \dD_i D_{i}^* e^{-t\t{\D}_i}
(D_{i}^{*})\ii\rrb
 =  \Tr\lrb  \dD D^* e^{-t\t{\D}_i}\rrb \ .$$
Hence equation \eqref{e:one} equals
\begin{eqnarray*}
& & \int_{0}^{\o} \Tr\lrb \dD^* D(e^{-t\D_1} - e^{-t\D_2}) \rrb +
\Tr\lrb \dD D^* (e^{-t\t{\D}_1} - e^{-t\t{\D}_2}) \rrb \, dt
\\ & = &  -\int_{0}^{\o} \frac{\dd}{\dd t }\Tr\lrb \dD^* D
(\D_1\ii e^{-t\D_1} - \D_2\ii e^{-t\D_2}) \rrb \\ & & \hskip 10mm
-\int_{0}^{\o} - \frac{\dd}{\dd t }\Tr\lrb \dD D^*(\t{\D}_1\ii
e^{-t\t{\D}_1} - \t{\D}_2\ii e^{-t\t{\D}_2}) \rrb \, dt \\ & = &
-\lim_{\e\to 0} \Tr\lrb \dD^* D (\D_1\ii e^{-t\D_1} - \D_2\ii
e^{-t\D_2}) \rrb |_{\e}^{1/\e} \\ & & \hskip 10mm - \lim_{\e\to 0}
\Tr\lrb \dD D^*(\t{\D}_1\ii e^{-t\t{\D}_1} - \t{\D}_2\ii
e^{-t\t{\D}_2}) \rrb |_{\e}^{1/\e} \\ & = & \Tr\lrb \dD^* D
(\D_1\ii- \D_2\ii) \rrb  + \Tr\lrb \dD D^*(\t{\D}_1\ii -
\t{\D}_2\ii \rrb \\ & = & \Tr\lrb \dD^* ((D_{1}^{*})\ii -
(D_{2}^{*})\ii) \rrb + \Tr\lrb \dD (D_1\ii - D_2\ii \rrb \ ,
\end{eqnarray*}
where we use \eqref{e:id2} for the final equality, and this
completes the proof.
\end{proof}

\begin{rem}
The variational equality \eqref{e:metrics6} also follows from
\eqref{e:variation1} applied to \eqref{e:leftside1.1}, along with
an analogue of \propref{p:alphaD}.
\end{rem}

\begin{cor}\label{cor:N}
For $P_i\in Gr_{\o}(D)$ with $D_{P_i}$ invertible
\begin{equation}\label{e:metrics7}
\frac{{\rm det}_{\z}(\D_{P_1})}{{\rm det}_{\z}(\D_{P_2})} =
\frac{{\rm det}_{F}\lrb S(P_1)^* S(P_1 )\rrb }{{\rm det}_{F}\lrb
S(P_2)^* S(P_2 )\rrb }.N(P_1,P_2) \ ,
\end{equation}
where  $N(P_1,P_2)$ depends only the boundary data. One has
\begin{equation}\label{e:Nmult}
  N(P_1,P_2) . N(P_2,P_3) = N(P_1,P_3) \ .
\end{equation}
\end{cor}

\vskip 2mm

Integrating \eqref{e:metrics6} over $[0,t]\subset(-\e,\e)$,
\eqref{e:metrics7} can be restated
\begin{equation}\label{e:metrics8}
\frac{{\rm det}_{\z}(\D_{t,P_1})}{{\rm det}_{\z}(\D_{t,P_2})}.
\frac{{\rm det}_{\z}(\D_{0,P_2})}{{\rm det}_{\z}(\D_{0,P_1})}=
\frac{{\rm det}_{F}\lrb S_t(P_1)^* S_t(P_1 )\rrb }{{\rm
det}_{F}\lrb S_t(P_2)^* S_t(P_2 )\rrb }.\frac{{\rm det}_{F}\lrb
S_0(P_2)^* S_0(P_2 )\rrb }{{\rm det}_{F}\lrb S_0(P_1)^* S_0(P_1
)\rrb } \ .
\end{equation}

Next, we make use of the homogeneous structure of the
Grassmannian to prove

\vskip 2mm

\noi $\bb$ {\bf $N(P_1,P_2)=1$}.

\vskip 2mm

We use a variational argument generalizing \cite{ScWo99}. Let
$$U_{\o}(H_Y) = U(H_Y)\cap  (I + \Psi_{\o}(H_Y))$$ be the group of
unitary operators on $H_Y = L^2 (Y, (E_1)_{|Y})$ differing from
the identity by a smoothing operator, and let $\t{Gr}_{\o}(D)$ be
the dense open subset of the index zero component of $Gr_{\o}(D)$
\begin{eqnarray}\label{e:tGr}
\t{Gr}_{\o}(D) & = & \{P\in Gr_{\o}(D) \ | \ D_P \ {\rm
invertible}\} \nonumber \\ & = & \{P\in Gr_{\o}(D) \ | \ S(P) :
H(D) \to \ran(P) \ {\rm invertible}\} \nonumber \\ &  = &
U_{H(D)}  \ ,
\end{eqnarray}
where the final equality refers to \eqref{e:Uw1}.

\begin{lem}\label{lem:P1=gP2}
For any $P_1, P_2\in \t{Gr}_{\o}(D)$ there exists a smooth path
\begin{equation}\label{e:pathP1P2}
I = g_0 \leq g_r \leq g_1 = g\ , \hskip 15mm 0\leq r \leq 1 \ ,
\end{equation}
in $U_{\o}(H_Y)$, defining smooth paths of projections  $P_{1,r}
= g_r P_1 g_r\ii$ and $P_{2,r} = g_r P_2 g_r\ii$ in
$\t{Gr}_{\o}(D)$ with $g P_1 g\ii =  P_2  $,
\begin{equation}\label{e:pathP1r}
P_1 \leq P_{1,r} \leq  P_2 \ ,
\end{equation}
and
\begin{equation}\label{e:pathP2r}
g P_1 g\ii  \leq P_{2,r} \leq  g P_2 g\ii \ .
\end{equation}
\end{lem}

\vskip 2mm

We hence obtain a real-valued strictly positive function
$g_r\mtoo N(P_{1,r},P_{2,r})$. The decisive fact is the following:
\begin{lem}\label{lem:dlogN=0}
\begin{equation}\label{e:dlogN=0}
\frac{d}{dr}\log N(P_{1,r},P_{2,r}) = 0 \ .
\end{equation}
\end{lem}

\vskip 2mm

The proofs will be given in a moment. Integrating
\eqref{e:dlogN=0} we have
\begin{equation}\label{e:Nt=N0}
 N(P_{1,1},P_{2,1}) =  N(P_{1,0},P_{2,0}) \ .
\end{equation}
From \eqref{e:pathP1r}, \eqref{e:pathP2r}, \eqref{e:Nt=N0} we
obtain $$N(P_2, g P_2 g\ii) = N(P_1, g P_1 g\ii) \ ,$$ and hence
that $ N(P,g P g\ii)$ depends only on $g\in U_{\o}$ and not on the
basepoint $P$. We define $$N(g):=N(P,g P g\ii)$$ where $P, g
Pg\ii \in \t{Gr}_{\o}(D)$. Then for $g_1, g_2 \in U_{\o}(H_Y)$,
from \eqref{e:Nmult} we have with $P, g_2 Pg_2\ii, g_1 g_2 P
g_2\ii g_1\ii \in \t{Gr}_{\o}(D)$,
\begin{eqnarray*}
  N(g_1 g_2) & = & N(P,g_1 g_2 P g_2\ii g_1\ii)\\
  & = & N(P,g_2 P g_2\ii). N(g_2 P g_2\ii,g_1 (g_2 P g_2\ii)
  g_1\ii)\\ & = & N(g_1) . N(g_2) \ .
\end{eqnarray*}
Thus $g\mto N(g)$ extends to a (Banach) character on
$U_{\o}(H_Y)$. It is a well-known and elementary fact that the
only such characters on $U_{\o}(H_Y)$ are $g\mto \det_F (g)$ ,
$g\mto \det_F (g\ii)$ or the trivial character $g\mto 1$. But $N$
is real-valued positive, while $\det_F$ on $U_{\o}(H_Y)$ takes
values in $U(1)$. Hence $N(g)=1$.

\vskip 2mm

\noi  This completes the proof of \thmref{t:thmA}. It remains to
prove the above Lemmas.

\vskip 2mm

\noi {\bf Proof of \lemref{lem:P1=gP2}} \ \ First, we have from
\eqref{e:tGr} that $\t{Gr}_{\o}(D)$ is path connected, and in fact
contractible. To show that a path of the asserted form exists we
prove that $U_{\o}(H_Y)$ acts transitively on the {\em index
zero} component of $Gr_{\o}(D)$, with non-contractible stabilizer
subgroup $U_{\o}(W) \times U_{\o}(W\pp)$ at $P\in
Gr^{(0)}_{\o}(D)$, $\ran(P) = W$. (The global homogeneous
structure on $Gr_{\o}(D)$ is usually studied via the action of a
restricted linear group \cite{PrSe86}, with contractible
stabilizer $U(W)\times U(W\pp)$, but our purposes on
$Gr^{(0)}_{\o}(D)$ are better suited to the $U_{\o}(H_Y)$
subgroup action.)

It is enough to give the path \eqref{e:pathP1P2} in $GL_{\o}(H_Y)
= GL(H_Y)\cap  (I + \Psi_{\o}(H_Y))$, the group of invertibles
congruent to the identity. For $U_{\o}(H_Y)$ is a retraction of
$GL_{\o}(H_Y)$ via the phase map
\begin{equation*}
GL_{\o}(H_Y) \too U_{\o}(H_Y) \ , \hskip 10mm g \mtoo u_g =
g|g|\ii \ ,
\end{equation*}
where $|g| :=  (g^* g)^{1/2}$. Here
\begin{equation}\label{e:gt}
(g^* g)^{t} = \frac{i}{2\pi}\int_{\g} \mu^{t}(g^* g - \mu)\ii
d\mu \ ,
\end{equation}
with $\g$ a contour surrounding ${\rm sp}(g^* g)$, is a smooth map
$\C\times GL_{\o}(H_Y) \to GL_{\o}(H_Y)$ (Lemma(7.10)
\cite{Qu88}). It follows that if $g_r$ is a path in
$GL_{\o}(H_Y)$ satisfying the properties of \lemref{lem:P1=gP2}
apart from unitarity, then $u_{g_r}$ will be the path required.
To see this, if $P_2 = gP_1 g\ii$ with $g\in GL_{\o}(H_Y)$, then,
since $u_{g}P_1 u_{g}\ii$ is a self-adjoint indempotent, to show
$P_2 = u_g P_1 u_g\ii$ we need only show $\ran(u_{g}P_1 u_{g}\ii)
= \ran(gP_1 g\ii)$. This is equivalent to showing $\ran(|g|P_1
|g|\ii) = \ran(P_1)$, but $gP_1 g\ii = P_2 = P_{2}^* =  P_2
P_{2}^*$ imply $\ran(|g|^2 P_1 (|g|^2)\ii) = \ran(P_1)$, and the
identity then follows from \eqref{e:gt}.

To define the operators $g_r\in GL_{\o}(H_Y)$ we modify an
argument of \cite{BoWo93} $\S$15. To begin with, choose $\e\in
(0,1)$ and suppose $\|P_1 - P_2\| \leq \e < 1$. Let
\begin{equation*}
g_r = I + \frac{r}{\e}(P_2 - P_1)(P_1 - P_1\pp) \ , \hskip 10mm
0\leq r\leq \e \ .
\end{equation*}
Clearly, $P_2 g_{\e}  = g_{\e} P_1$ and since $\|g_r - I\| < 1$
then $g_r$ is invertible. Moreover, since $P_1 - P_2$ is
smoothing, so is $g_r - I$ and hence $g_r\in GL_{\o}(H_Y)$. Since
$D_{P_{i,0}}$ is invertible and invertibility is an open
condition for continuous families of Fredholm operators, then by
taking $\e$ smaller if necessary, $D_{P_{i,r}}$ will be
invertible for $0\leq r\leq \e$. Hence $g_r$ defines locally a
path of the type required. Now $Gr^{(0)}_{\o}(D)$ is path
connected and hence for arbitrary $P, P^{'} \in Gr^{(0)}_{\o}(D)$
we can find a finite sequence $P_1 = P, \ldots, P_m = P^{'},$ in
$Gr^{(0)}_{\o}(D)$ with $\|P_i - P_{i+1}\| \leq \e_i$, $\e_i\in
(0,1)$, for $i=1,\ldots, m-1,$ and a finite sequence of paths
$g_{r_i}$ in $GL_{\o}(H_Y)$ with $$P_i \leq
g_{r_i}P_{i}g_{r_i}\ii\leq P_{i+1}\in \t{Gr}_{\o}(D) \ .$$
Finally, rescaling so that $0\leq r_i \leq 1$ for each path, then
$g_r = g_{r_{m-1}}\ldots g_{r_1}$ is a path in $GL_{\o}(H_Y)$ of
the required form. This completes the proof. \vskip 2mm

\noi {\bf Proof of \lemref{lem:dlogN=0}} \ \ Since we consider the
simultaneous action of $U_{\o}$ on $P_1 , P_2$, we can `gauge'
transform the boundary variation to an order $0$ variation of
$\wD$, and then appeal to \propref{p:metrics6}. First, notice
that the action
\begin{equation*}
  (g_r, P_i) \mtoo g_r P_i g_r\ii = P_{i,r}
\end{equation*}
induces a dual action on the adjoint boundary condition
\begin{equation*}
  (\t{g}_r, P_i\st) \mtoo \t{g}_r P_i\st \t{g}_r\ii = (g_r P_i g_r\ii)\st =
  P_{i,r}\st \ ,
\end{equation*}
where $\t{g}_r = \s g_r \s\ii$. Moreover, with $\wg_r =
\begin{pmatrix}
  g_r & 0 \\
  0 & \t{g}_r
\end{pmatrix}$, we have
\begin{equation*}
\wg_r\wP_{i}\wg_r\ii = P_{i,r} \oplus P_{i,r}\st =
\widehat{P_{i,r}} \ ,
\end{equation*}
and
\begin{equation}\label{e:wswg}
\wg_r\ws = \ws\wg_r \ .
\end{equation}
We can now transform the self-adjoint global boundary problem
$\wD_{\wP_{i,r}}$ to a unitary equivalent operator $\wD(r)_{\wP}$
with constant domain by the method of \cite{Wo99,ScWo99}. Let
$f:[0,1]\to [0,1]$ be a non-decreasing function  with $f(u) = 1$
for $u < 1/4$ and $f(u) = 0$ for $u> 3/4$. Then we extend $g_r$
and $\t{g}_r$ to unitary transformations
\begin{equation*}
U_r= \begin{cases} g_{rf(u)}  \  \  {\rm on}  \ \ \{u\} \times Y
\subset U \\ Id \ \hskip 7mm {\rm on}  \ \ X \setminus U
\end{cases},  \hskip 2 mm \t{U}_r = \begin{cases} \t{g}_{rf(u)} \ \
{\rm on} \ \ \{u\} \times Y  = U
\\ Id \hskip 7mm \ {\rm on} \ \ X \setminus U
\end{cases},
\end{equation*}
on $L^2 (X,E^1)$ and  $L^2 (X,E^2)$, respectively, and hence to a
unitary transformation
\begin{equation*}
\widehat{U}_r= \begin{cases} \wg_{rf(u)}  \  \ {\rm on} \ \ \{u\}
\times Y \subset U \\ Id \hskip 7mm \ {\rm on} \ \ X \setminus U
\end{cases}  =  U_r\oplus \t{U}_r \ ,
\end{equation*}
on $L^2 (X,E^1 \oplus E^2).$ Then
\begin{equation*}
\wD_{\wP_{i,r}} \hskip 5mm {\rm and} \hskip 5mm \wD(r)_{P_i} :=
(\widehat{U}_r\ii \wD\widehat{U}_r)_{\wP_i}
\end{equation*}
are {\em unitarily equivalent}. Moreover, it is easy to check that
\begin{equation}\label{e:wdr}
\wD(r)_{P_i} = (\widehat{\D_r})_{\wP_i} \ ,
\end{equation}
where $\D_r = U_r\ii\D U_r = D_{r}^{*}D_r$ and $D_r$ is the
Dirac-type operator $D_r  = \t{U}_r\ii D U_r.$

Next, since $P(\wD(r)_{\mu}) = \wg_r\ii P(\wD_{\mu})\wg_r$, from
the multiplicativity of $\det_F$ we obtain
\begin{equation*}
{\rm det}_{F}\left(\frac{S_{\mu}(\wP_{1,r})}{\wE S_{\mu}
(\wP_{2,r})}\right) = {\rm
det}_{F}\left(\frac{S_{r,\mu}(\wP_1)}{\wE_r S_{r,\mu}
(\wP_2)}\right) \ ,
\end{equation*}
where $S_{r,\mu}(\wP_1) = \wP_1\circ P(\wD(r)_{\mu})$ and $\Ee_r =
\wg_r\ii\Ee\wg_r$. Hence from \eqref{e:relLapdet} and
\eqref{e:wdr} $$(d/dr)\log{\rm det}_{\z}(\D_{P_{1,r}},
\D_{P_{1,r}}) = (d/dr)\log{\rm det}_{\z}((\D_r)_{P_1},
(\D_r)_{P_2}) \ .$$

Finally, since \eqref{e:wswg} holds, then $\wD(r)_{|U}$ has the
form \eqref{e:collarz} with $\widehat{\Aa}_r = \wg_r\ii
\widehat{\Aa}\wg_r\ii$, and since $\wg_r$ differs from the
identity by a smoothing operator then $\s(\widehat{\Aa}_r)$ is
independent of $r$. Hence we can apply \propref{p:metrics6} and
the identity
\begin{equation*}
{\rm det}_F (S(P_{i,r})^* S(P_{i,r})) =  {\rm det}_F (S_r(P_{i})^*
S_r(P_{i})) \ ,
\end{equation*}
which is a consequence of $P(D_r) = g_r\ii P(D) g_r$, to complete
the proof.

\vskip 2mm

\noi {\bf Proof of \lemref{lem:wCalderon}} \ \  We have
$P(\wD_{\la}) = \g\wKk_{\la}$, where \eqref{e:poisson}
\begin{equation*}
  \wKk_{\la} = \widehat{r}\wD_{\la,d} \ii\wga^*\ws :
  H^{s-1/2}(Y,E_{|Y}^1 \oplus E_{|Y}^2)
  \too \Ker(\wDla,s) \subset H^s (X,E^1 \oplus E^2) \ ,
\end{equation*}
and $\wD_{\la,d}$ is an invertible operator over the double
manifold $\tilde{X}$ with $(\wD_{\la,d})_{|X} = \wD_{\la}$, $\wga
= \g_1 \oplus \g_2$, $\widehat{r} = r_1 \oplus r_2$ with $\g_i
:H^s (\tilde{X},\tilde{E}^i)\to H^{s-1/2}(Y,E_{|Y}^i)$,  $r_i :H^s
(\tilde{X},\tilde{E}^i)\to H^{s}(X,E^i)$ the restriction
operators, and $\ws$ is defined in \eqref{e:ws}.

Let $D_d$ be the double operator of $D$, see for example
\cite{BoWo93}. Then $D_d$ is invertible on the closed double
manifold $\t{X}$ with $(D_d)_{|X} = D$, and hence
\begin{equation*}
  \wD_{\la,d} =\begin{pmatrix}
   -\la  & D^{*}_d \\
    D_d & -I \
  \end{pmatrix} \
\end{equation*}
is invertible on $\t{X}$ with $(\wD_{\la,d})_{|X} = \wDla$. We
compute
\begin{equation*}
  \wD_{\la,d}\ii =\begin{pmatrix}
   (\D_d - \la)\ii  &  D^{*}_d (\t{\D}_d -\la)\ii \\
    D_d(\D_d - \la)\ii & \la (\t{\D}_d -\la)\ii \
  \end{pmatrix} \ ,
\end{equation*}
where $\D_d = D^{*}_d D_d, \t{\D}_d = D_d D^{*}_d$, and hence that
\begin{equation*}
  \wKk_{\la} =\begin{pmatrix}
   r_1 D^{*}_d (\t{\D}_d -\la)\ii \g_{1}^* \s  &  r_1 (\D_d -\la)\ii \g_{0}^* \s\st \\
    \la r_2 (\t{\D}_d -\la)\ii \g_{1}^* \s  &  r_1 D_d (\D_d -\la)\ii \g_{0}^* \s\st
  \end{pmatrix} \ ,
\end{equation*}
with $\s\st := -\s\ii$. Setting $\la = 0$ we obtain
\begin{equation*}
  P(\wD) = \g\begin{pmatrix}
   r_1 D_d\ii \g_{1}^* \s  &  r_1 D_d\ii (D^{*}_d)\ii \g_{0}^* \s\st \\
    0  &  r_1 (D^{*}_d)\ii \g_{0}^* \s\st
  \end{pmatrix} =
  \begin{pmatrix}
   \g r_1 D_d\ii \g_{1}^* \s  &  \g D_X\ii K_{D^*} \\
    0  &  \g r_1 (D^{*}_d)\ii \g_{0}^* \s\st
  \end{pmatrix} \ ,
\end{equation*}
and since $D_X\ii = D_{P(D)}\ii$ we reach the conclusion.


\section{An Application to Ordinary Differential Operators}

In dimension one we  can do better. Because no basepoint is needed
to define the Grassmannian it is possible to apply the method of
\thmref{t:thmB} to obtain formulas for the $\z$-determinant of
individual boundary problems, rather than just relative formulas.

\vskip 2mm

\noi $\bb$ \ {\bf First-order operators.} \ We consider, as in $\S
4$,  a first-order elliptic differential operator $D :
C^{\infty}(X;E) \too C^{\infty}(X;F)$, but where now $X = [0,
\b]$, $\b>0$, and $E , F$ are Hermitian bundles of rank $n$.
Relative to trivializations of $E, F$ one has $D= A(x)d/dx +
B(x)$, where $A(x),B(x)$ are complex $n\times n$ matrices and
$A(x)$ is invertible. The restriction map to the boundary  $Y =
\{0\} \sqcup \{\b\}$ is the map $\gamma:H^{1}(X;E)\too E_{0}\oplus
E_{\beta},$ with $ \gamma(\psi) = (\psi (0),\psi (\beta)),$ and so
global boundary conditions  for $D$ are parameterized by the
Grassmannian $Gr (E_0\oplus E_{\b})$ of the $2n$-dimensional space
of boundary `fields'. For each $P \in Gr (E_0\oplus E_{\b})$ we
have a boundary problem $ D_P: {\rm dom}(D_{P})\too L^{2}(X;E),$
where
\begin{eqnarray*}
{\rm dom}(D_{P}) & = &  \{\psi\in H^{1}(X;E) \ | \ P\gamma\psi =
0\} \\ & = &  \{\psi\in H^{1}(X;E) \ | \ M\psi(0) + N\psi(\b) =
0\}
\end{eqnarray*}
and  $[M \; N ] \in \Hom(E_0\oplus E_{\b},E_0)$ are Stiefel
coordinates for $P$, see \eqref{e:Pstiefel}.

In dimension one, any element of $\Ker (D)$ has the form $K(x)v$
for some $v\in E_{0}$, where $K(x)\in \Hom(E_0,E_x)$ is the
parallel transport operator uniquely solving $D K(x) = 0$ subject
to $K(0)=I$. The isomorphism $\gamma : \Ker (D)\to H(D)$ is clear,
while $H(D) = \gr(K: E_0\to E_{\b}) \subset E_{0}\oplus
E_{\beta}$, where $K := K(\b)$.

Notice that $P(D) \in Gr_n (E_0\oplus E_{\b})$ and hence the
component $Gr_k (E_0\oplus E_{\b})$ of the Grassmannian with
$\tr(P)=k$ is the component with operator index $\ind D_{P} = \ind
\Ss(P) = n - k$. The Poisson operator $\Kk_D : E_0 \oplus E_{\b}
\too \Ci(X,E)$ is the operator $\Kk_D (u)(x) = K(x)p_0 P(D) u$,
where $p_0$ is the projection map $E_0 \oplus E_{\b} \to E_0$. It
is easy to check that $D_{P(D)}\ii D = I - \Kk_D \g$ and hence
that for invertible global boundary problems \eqref{e:relinv}
holds:
\begin{equation}\label{e:relinvdim1}
D_{P_1}\ii - D_{P_2} \ii   =  - \Kk_D (P_1)P_1\g  D_{P_2} \ii \ .
\end{equation}
On the other hand, it is well known from elementary considerations
that in Stiefel coordinates $D_P\ii$ has kernel
 \vskip 1mm
\begin{equation}\label{e:kernel}
k_{P}(x,y) = \left\{\begin{array}{l}
 -K(x)((M + N K)\ii N K)K(y)^{-1} A(y)\ii \hskip 15mm x<y \\
\\[1mm]
 \hskip 3mm K(x)(I-(M + N K)\ii N K)K(y)^{-1}A(y)\ii \hskip 8mm x>y \ ,
\end{array}
\right.
\end{equation}
 \vskip 1mm
Hence, if $P_1, P_2$ are represented by Stiefel coordinates $[M_1
\; N_1], [M_2 \; N_2]$, then the relative inverse $D_{P_1}\ii
-D_{P_2}\ii$ has the  smooth kernel
\begin{equation}\label{e:relkernel}
 -K(x)\left((M_1 + N_1 K)\ii N_1 -(M_2 + N_2 K)\ii
N_2\right)K.K(y)^{-1}A(y)\ii \ ,
\end{equation}
which can also be computed directly from \eqref{e:relinvdim1} by
using \eqref{e:Rinverse} with $\Phi(R):=\Pi(S(P))$.

We assume that $D_P$ is invertible with spectral cut $R_{\th}$.
For $\re(s) >0$ we can then define $D_P\si =
\frac{i}{2\pi}\int_{C}\la_{\th}\si\, (D_P - \la)\ii \,
  d\la \ .$
Let $k_{P,\la}(x,y)$ be the kernel of $(D_P - \la)\ii$. From
\eqref{e:kernel} one has $\lim_{\e\to 0} (k_{P,\la}(x,x+\e) -
k_{P,\la}(x+\e,x) = -A(x)\ii$, and hence for $\re(s) > 1$ the
kernel $p_s (x,y) = \frac{i}{2\pi}\int_{C}\la_{\th}\si\,
k_{P,\la}(x,y)  \, d\la \ $ of $D_P\si$ is continuous,  and
$D_P\si$ is trace class.  Moreover, if $P(0)$ is the projection
onto $E_0$, with Stiefel graph coordinates $[I \; 0]$, then from
\eqref{e:kernel} we have $\Tr(D_{P(0)}\si) = 0$. For $\re(s) > 1$,
\begin{equation}\label{e:zeta1}
\z_{\th}(D_P,s) = \z_{\th}(D_P, D_{P(0)},s) \ ,
\end{equation}
is therefore a relative $\z$-function. Hence
\begin{eqnarray*}
 \z_{\th}(D_P,s) & = & \frac{i}{2\pi}\int_{C}\la\si\,
  \Tr((D_P - \la)\ii - (D_{P(0)} - \la)\ii) \, d\la  \nonumber \\
  & = & -\frac{i}{2\pi}\int_{C}\la\si\,
  \frac{\dd}{\dd\la}\log {\rm det} \lrb \frac{S_{\la}(P)}{PS_{\la}(P(0))}\rrb \, d\la
  \\
  & = & -\frac{i}{2\pi}\int_{C}\la\si\,
  \frac{\dd}{\dd\la}\log {\rm det}(M + N K_{\la}) \, d\la \ ,
\end{eqnarray*}
where $K_{\la}(x)$ is the solution operator for $D-\la$. The
second equality is a restatement of \eqref{e:relresolventtrace},
note $\det_F$ becomes the usual determinant here, and the third
equality follows from \eqref{e:graphscattering2} and
\remref{rem:stiefel} (with $[M_1 \; N_1] = [M \; N], [M_2 \; N_2]
= [I \; 0]$). Alternatively, with $\Mm_{\la} = M + N K_{\la}$, one
can compute directly from \eqref{e:relkernel} $$ \Tr((D_P -
\la)\ii - (D_{P(0)} - \la)\ii)  $$
\begin{eqnarray}
& = & -\int_{0}^{\b}\tr\lsb K_{\la}(x)\Mmla\ii N
K_{\la}K_{\la}(x)\ii A(x)\ii \rsb \, dx   \nonumber \\ & = &
\int_{0}^{\b}\tr\lsb\dla(D-\la)K_{\la}(x)\Mmla\ii N
K_{\la}K_{\la}(x)\ii A(x)\ii \rsb \, dx  \nonumber\\ & = &
-\int_{0}^{\b}\tr\lsb(D-\la)\dla (K_{\la}(x))\Mmla\ii N
K_{\la}K_{\la}(x)\ii A(x)\ii \rsb \, dx \nonumber\\ & = &
-\int_{0}^{\b}\tr\lsb \frac{d}{dx} \lrb K_{\la}(x)\ii \dla
(K_{\la}(x)) \Mmla\ii N K_{\la}\rrb \rsb \, dx \nonumber\\ & = &
-\dla \log\det\,\Mmla \label{e:Mla} \ ,
\end{eqnarray}
using $K_{\la}(x)\ii A(x)\ii (D-\la)K_{\la}(x) = d/dx$ and the
$\la$ derivative of $(D-\la)K_{\la}(x)=0$.

We now assume $D_P$ defines an elliptic boundary problem in the
sense of \cite{Se66}. Then there is an asymptotic expansion as
$\la\to\o$ in a sector $\Lambda_{\th}$
\begin{equation}\label{e:order1asymp1}
\Tr((D_P - \la)\ii - (D_{P(0)} - \la)\ii) \sim
\sum_{j=1}^{\o}b_{j} (-\la)^{-j} \ .
\end{equation}
More precisely, $(D_{P_i} - \la)^{-2}$ is trace class and by
ellipticity $\Tr((D_{P_i} - \la)^{-2}) \sim \sum_{j=1}^{\o}a_{j}
(-\la)^{-j-1}$ as $\la\to\o$, and so applying
\eqref{lem:asympSla} to the relative trace and observing that the
trace class condition on the relative resolvent implies terms
$(-\la)^{-\a_j}\log(-\la)$ with $\a_j \leq 0$ vanish, then
\eqref{e:order1asymp1} follows.

Thus $\z_{\th}(D_P, D_{P(0)},s)$ defines the meromorphic
continuation of $\z_{\th}(D_P,s)$ to $\C$  via the resolvent trace
expansion \eqref{e:order1asymp1}, and this is regular at $s=0$.
Hence we can define ${\rm det}_{\z,\th}(D_P)$. In dimension one
the following stronger variant of \thmref{t:thmD*} ($\S 4.4$)
holds:

\vskip 2mm

\begin{thm}\label{t:thmD**} Let $D_P$ be a
first-order elliptic boundary problem  over $[0,\b]$ and let $[M
\; N]$ be Stiefel coordinates for $P$. Then
\begin{equation}\label{e:stiefelversiondim1}
{\rm det}_{\z,\th}(D_P)= {\rm det}(M + NK)\, . \,
e^{-\LIM^{\th}_{\la\to\infty}\log {\rm det} \lrb M + NK_{\la}\rrb
} \ .
\end{equation}
Invariantly, one has
\begin{equation}\label{e:stiefelversiondim2}
{\rm det}_{\z,\th}(D_P)= {\rm det} \lrb
\frac{S(P)}{PS(P(0))}\rrb\, . \, e^{-\LIM^{\th}_{\la\to\infty}\log
 {\rm det} \lrb \frac{S_{\la}(P)}{PS_{\la}(P(0))}\rrb } \ .
\end{equation}
\end{thm}
\begin{proof}
Immediate from \eqref{e:zeta1}, \eqref{e:Mla},
\eqref{e:order1asymp1} and \propref{p:thmB*} with $\Phi(\la) =
\log\det \Mmla$, or from \thmref{t:thmB} with $\Ss_{\la}$
replaced by $\Mm_{\la}$.
\end{proof}

\vskip 2mm

\noi $\bb$ \ {\bf Operators of order $\geq 2$.} \ There is a
straightforward generalization of these formulas to differential
operators  $D :\Ci(X,E)\to \Ci(X,F)$ of order $r\geq 2$. With
respect to trivializations of $E,F$
\begin{equation*}\label{e:edo}
D = \sum_{k=0}^{r}B_{k}(x)\frac{d^k}{dx^k} : C^{\infty}(X;\Cn)\to
C^{\infty}(X;\Cn) \ ,
\end{equation*}
with complex matrix coefficients and $\det B_{r}(x)\neq 0$. The
restriction map  is
\begin{equation}\label{e:cauchy}
\gamma_{r-1}:H^{r}(X;E)\too \Crn\oplus \Crn \ ,\hskip 5mm
\gamma_{r-1}(\psi) = (\wpsi (0),\wpsi (\beta)) \ ,
\end{equation}
where $\wpsi(x) = (\psi(x),\ldots,\psi^{(r-1)}(x))$.  The form of
$\g_{r-1}$ means that one can study boundary problems for $D$
through the first-order system on $C^{\infty}(X;\Cnr)$
\begin{equation}\label{e:Dequiv}
\wDd = \frac{d}{dx} - \left[\begin{array}{cccc} 0   &   I   &
\ldots & 0\\ 0   &   0   &       \ldots  &   0\\ \vdots  & \vdots
& \vdots & \vdots\\ 0   &   0   &       \ldots  & I\\
-B_{r}(x)^{-1}B_{0}(x) &   -B_{r}(x)^{-1}B_{1}(x)
    &   \ldots  &   -B_{r}(x)^{-1}B_{r-1}(x)
\end{array}\right] \ .
\end{equation}
This is well known \cite{Fo87,LeTo99}. $\wDd$ extends to a
continuous map $H^{1}(X;\Cnr)\to L^{2}(X;\Cnr)$, and with respect
to the inclusion $\hat{\iota}:H^{r}(X;\C^n)\too
H^{1}(X;\Cnr),\;\psi \mapsto \wpsi,$ we have $\gamma_{r-1} =
\gamma \circ \hat{\iota}.$ More precisely, there is an isomorphism
$\Ker(D) \cong \Ker(\wDd)$, for from
\begin{equation}\label{e:ker}
\wDd (\psi,\ldots,\psi^{(r-1)}) = (0,\ldots,0,B_{r}^{-1}D \psi) \
,
\end{equation}
a basis $\{\psi_1,\ldots,\psi_k\}$ for $\Ker(D)$ defines a basis
$\{\wpsi_1,\ldots,\wpsi_k\}$ for $ \Ker (\wDd)$. A solution of
$\wDd$ is characterized by its parallel transport operator
$\wK(x)$,  as before, the columns of which are a preferred basis
for $ \Ker (\wDd)$. One has  $$\ran(P(\wDd)) =
\gr(\wK:\Cnr\to\Cnr)$$ and $\wDd$ has Poisson operator
$\wKk:\Cnr\oplus\Cnr \to \Ci (X;\Cnr)$ defined by $$\wKk(v)(x) =
\wK(x)p_{0}P(\wD)v \ .$$

A global boundary condition for $D$ is defined by a global
boundary condition  $P\in Gr(\Cnr\oplus \Cnr)$ for $\wDd$. That is
\begin{eqnarray*}
{\rm dom}(D_{P}) & = & \{\psi\in H^{r}(X;E) \ | \
P\gamma_{r-1}\psi = 0\} \\ & = & \{\psi\in H^{r}(X;E) \ | \
\wM\wpsi(0) + \wN\wpsi(\b) = 0\} \ ,
\end{eqnarray*}
where $[\wM,\wN]$ are Stiefel coordinates for $P$.

The boundary problem $D_P$ is modeled by the finite-rank operator
on boundary data $\wS(P) := P\circ P(\wDd):K(\wDd)\to \ran(P)$.
From \eqref{e:ker} we have that $D_P$ is invertible if and only if
$\wDd_P$ is invertible, and in that case
\begin{equation*}
D_{P}^{-1} = [\wDd_{P}^{-1}]_{(1,r)}B_{r}^{-1} : L^{2}(X;\Cn) \to
\dom (D_P) \ .
\end{equation*}
Here $[\wDd_{P}^{-1}]_{(1,r)}$ means the integral operator
$\int_{0}^{\b}[\wk(x,y)]_{(1,r)}\psi (y)\, dy$ where $\wk(x,y)$ is
the kernel of $\wDd_{P}^{-1}$, and, as in \cite{LeTo99},
$[T]_{(1,r)}$ is the $n\times n$ matrix in the $(1,r)^{{\rm th}}$
position in an $r\times r$ block matrix $T\in \End (\Cnr)$. For
$D_{P_1}, D_{P_2}$ invertible, this leads to the formula
\begin{equation}\label{e:inverse2}
  D_{P_1}\ii = D_{P_2}\ii -
  [\wKk \wS(P_1)\ii P_1 \gamma\wDd_{P_2}\ii]_{(1,r)}B_{r}\ii \ .
\end{equation}
 In Stiefel
coordinates $D_P\ii$ has kernel
\begin{equation}\label{e:kernelr}
k_{P}(x,y) = \left\{\begin{array}{l}
 -[\wK(x)(\wMm^{-1}\wN\wK) \wK(y)^{-1}]_{(1,r)}B_{r}(y)^{-1} \hskip
15mm x<y \\
\\[1mm]
 \hskip 3mm [\wK(x)(I-\wMm^{-1}\wN\wK) \wK(y)^{-1}]_{(1,r)}B_{r}(y)^{-1} \hskip 8mm x>y,
\end{array}
\right.
\end{equation}
 \vskip 1mm
\noi where $\wMm=\wM+\wN\wK$.

Since we are in dimension one, the resolvent of a differential
operator of order $r\geq 2$ is trace class (as is evident from
\eqref{e:kernelr}). Let $R_{\th}$ be a spectral cut for $D_P$ and
let $P$ now be a local elliptic boundary condition for $D$
 \cite{Se66}. Then Seeley proved
\cite{Se69a} that as $\la\to \o$ in $\Lambda_{\th}$ the resolvent
trace has an asymptotic expansion
\begin{equation}\label{e:asymporderr}
  \Tr((D_P -\la)\ii) \sim \sum_{j=-1}^{\o} b_j (-\la)^{-j/r -1} \ .
\end{equation}
On the other hand, $(D-\la)_{P}\ii =
[\wDd_{\la,P}^{-1}]_{(1,r)}B_{r}^{-1}$ with $\wDd_{\la} =
\widehat{(D-\la)}$,  with $\wK_{\la}(x)$ the parallel transport
operator for $\wDd_{\la}$ and $\wMm_{\la} = \wM + \wN\wK_{\la}$
\begin{eqnarray*}
  \Tr((D_P -\la)\ii) & =  & \int_{0}^{\b}\tr\lsb\dla(\wDd_{\la})\wK_{\la}(x)\wMmla\ii \wN
\wK_{\la}\wK_{\la}(x)\ii \rsb \, dx
\\ & = & -\dla \log\det\,\wMmla \ ,
\end{eqnarray*}
This follows by the same argument as before, using the device
$$\tr([T(x)]_{(1,r)}) = \tr(J T(x)) = -\tr(\dla(\wDd_{\la}) T(x))
\ ,$$ where $J$ is the $n\times n$ block matrix with the identity
in the $(r,1)^{{\rm th}}$ position and zeroes elsewhere.

By \propref{p:thmB*} we therefore obtain the extension of
\thmref{t:thmD**} to higher order operators:

\begin{thm}\label{t:thmE} Let $D_P$ be a local elliptic boundary
 problem of order $r\geq 2$ over $[0,\b]$ and let $[\wM \; \wN]$ be
Stiefel coordinates for $P\in Gr(\Cnr\oplus\Cnr)$. Then
\begin{equation}\label{e:stiefelversiondim1orderr}
{\rm det}_{\z,\th}(D_P)= {\rm det}(\wM + \wN\wK)\, . \,
e^{-\LIM^{\th}_{\la\to\infty}\log {\rm det} \lrb \wM +
\wN\wK_{\la}\rrb } \ .
\end{equation}
Invariantly, one has
\begin{equation*}
{\rm det}_{\z,\th}(D_P)= {\rm det} \lrb
\frac{\wS(P)}{P\wS(P(0))}\rrb\, . \,
e^{-\LIM^{\th}_{\la\to\infty}\log
 {\rm det} \lrb \frac{\wS_{\la}(P)}{P\wS_{\la}(P(0))}\rrb } \ ,
\end{equation*}
where $P(0)$ is the projection to the first factor in
$\Cnr\oplus\Cnr$.
\end{thm}

\vskip 1mm  The formulas \eqref{e:stiefelversiondim1} and
\eqref{e:stiefelversiondim1orderr} were first proved in
\cite{LeTo99}.

\vskip 8mm

\noi {\small King's College, \newline \noi London. \vskip 2mm \noi
Email: sgs@mth.kcl.ac.uk}

\end{document}